\documentclass[11pt]{article}
\usepackage{amsfonts}
\usepackage{latexsym,amssymb,amsfonts}
\usepackage{mathrsfs}
\usepackage{graphicx}
\usepackage{psfrag}
\usepackage{amsmath, amsbsy}
\usepackage{amsopn, amstext}
\usepackage{amsmath,multicol,amsopn}
\usepackage{latexsym,amssymb,color}
\usepackage{authblk}
\usepackage{enumerate} %列表
\usepackage{graphicx}
\usepackage{caption}
\usepackage{booktabs}
\usepackage{framed}
\usepackage{subfig}
\usepackage{float}
\usepackage{caption}
\usepackage{cite}

\usepackage{amsthm}

\def\sqr#1#2{{\vcenter{\vbox{\hrule height.#2pt
              \hbox{\vrule width.#2pt height#1pt \kern#1pt \vrule width.#2pt}
              \hrule height.#2pt}}}}
\def\signed #1{{\unskip\nobreak\hfil\penalty50
              \hskip2em\hbox{}\nobreak\hfil#1
              \parfillskip=0pt \finalhyphendemerits=0 \par}}
\def\endpf{\signed {$\sqr69$}}

\def\dbR{{\mathop{\rm l\negthinspace R}}}

\def\3n{\negthinspace \negthinspace \negthinspace }
\def\2n{\negthinspace \negthinspace }
\def\1n{\negthinspace }

\def\dbF{{\mathbb{F}}}

\def\dbN{{\mathbb{N}}}

\def\dbP{{\mathbb{P}}}

\def\dbR{{\mathbb{R}}}

\def\ds{\displaystyle}

\def\={\buildrel \triangle \over =}

\def\a{\alpha}
\def\b{\beta}
\def\g{\gamma}
\def\d{\delta}
\def\e{\varepsilon}

\def\k{\kappa}
\def\l{\lambda}

 \def\n{\nabla}

\def\t{\times}

\def\th{\theta}
\def\om{\omega}

\def\i{\infty}
\def\ns{\noalign{\ss} }

\def\pa{\partial}

\def\G{\Gamma}
\def\D{\Delta}

\def\Si{\Sigma}

\def\Om{\Omega}

\def\cA{{\cal A}}
\def\cB{{\cal B}}

\def\cF{{\cal F}}

\def\cl{{\cal l}}
\def\mE{{\mathbb{E}}}

\def\no{\noindent}

\def\ss{\smallskip}
\def\ms{\medskip}
\def\bs{\bigskip}
\def\q{\quad}
\def\qq{\qquad}

\def\lan{\big\langle}
\def\ran{\big\rangle}

\def\max{\mathop{\rm max}}

\def\pa{\partial}

\def\wt{\widetilde}
\def\cd{\cdot}
\def\cds{\cdots}

\def\ae{\hbox{\rm a.e.{ }}}
\def\as{\hbox{\rm a.s.{ }}}

\def\cl{\overline}

\def\deq{\mathop{\buildrel\D\over=}}

\def\({\Big (}
\def\){\Big )}
\def\[{\Big[}
\def\]{\Big]}
\def\be{\begin{equation}}
\def\bel{\begin{equation}\label}
\def\ee{\end{equation}}
\def\bt{\begin{theorem}}
\def\bcd{\begin{condition}}
\def\ecd{\end{condition}}
\def\et{\end{theorem}}
\def\bc{\begin{corollary}}
\def\ec{\end{corollary}}
\def\bde{\begin{definition}}
\def\ede{\end{definition}}
\def\bl{\begin{lemma}}
\def\el{\end{lemma}}
\def\bp{\begin{proposition}}
\def\ep{\end{proposition}}
\def\br{\begin{remark}}
\def\er{\end{remark}}
\def\ba{\begin{array}}
\def\ea{\end{array}}
\def\ed{\end{document}}
\def\ns{\noalign{\ms}}
\def\ds{\displaystyle}

\def\square#1{\vbox{\hrule\hbox{\vrule height#1%
     \kern#1\vrule}\hrule}}
\def\rectangle#1#2{\vbox{\hrule\hbox{\vrule height#1%
     \kern#2\vrule}\hrule}}

\font\tenbb=msbm10 \font\sevenbb=msbm7
\font\fivebb=msbm5

\newfam\bbfam
\scriptscriptfont\bbfam=\fivebb
\textfont\bbfam=\tenbb
\scriptfont\bbfam=\sevenbb

\newtheorem{lemma}{Lemma}[section]
\newtheorem{theorem}{Theorem}[section]
\newtheorem{corollary}{Corollary}[section]

\newtheorem{definition}{Definition}[section]
\newtheorem{proposition}{Proposition}[section]
\newtheorem{condition}{Condition}[section]
\theoremstyle{remark}
\newtheorem{remark}{\rm \bf Remark}[section]
\newtheorem{example}{Example}[section]
\makeatletter
   
   \@addtoreset{equation}{section}
\makeatother

\makeatletter
   
   \@addtoreset{figure}{section}
\makeatother

\allowdisplaybreaks

\begin{document}

\title{\bf Stability and regularization for ill-posed Cauchy problem of a stochastic parabolic differential equation}

\author[1]{Fangfang Dou\thanks{Corrsponding author. Email: fangfdou@uestc.edu.cn.}}
\author[1]{Peimin L\"{u}}
\author[2]{Yu Wang}
\affil[1]{\small\it School of Mathematical Sciences, University of Electronic Science and Technology of China, Chengdu, China}
\affil[2]{\small\it School of Mathematics, Southwest Jiaotong University, Chengdu, China.}
\date{}
\maketitle

\begin{abstract}
In this paper, we investigate an ill-posed Cauchy problem involving a stochastic parabolic equation. We first establish a Carleman estimate for this equation. Leveraging this estimate, we derive the conditional stability and convergence rate of the Tikhonov regularization method for the aforementioned ill-posed Cauchy problem. To complement our theoretical analysis, we employ kernel-based learning theory to implement the completed Tikhonov regularization method for several numerical examples.
\end{abstract}

\bs

\no{\bf 2020 Mathematics Subject
Classification}. 35R30,65N21,60H15.

\bs

\no{\bf Key Words}. Carleman estimate,
Cauchy problem of stochastic parabolic differential equation, conditional stability, regularization, numerical approximation.

\section{Introduction }

   %Thanks to their important application requirement, simultaneous identification of several parameters in parabolic equations and slow diffusion processes caught much attention. The uniqueness and stability for determination of the space-dependent reaction coefficient, the initial temperature and the source term from measured temperatures at different instants in \cite{CLL2020,CL2019,JL2008,WCQW2020,YZ2001}, boundary/internal measurement data \cite{YZ2001}, or both boundary/internal measurement and final temperature\cite{CJW2022,YLD2020,ZW2014}. In these works, Gel’fand–Levitan theory, spectral theory, Carleman estimates, regularization methods and iterated numerical algorithms are well applied. We also mention that these problems are also considered for time-fractional different equations, which are used to describe slow diffusion processes, see \cite{LW2019,LYY2015,ZL2021} for example.

%the mild solution of the stochastic parabolic differential equation is given by the fundamental solution of the deterministic parabolic equation. As far as we know, there is no published work mentioned about this kind of practice.
\medskip
To beginning with, we introduce some notations concerning stochastic analysis. More details for these can be found in \cite{LZ2021}.

Let $(\Om,\cF,\dbF,\dbP)$ with $\dbF=\{\cF_t\}_{t\ge 0}$  be
a complete filtered probability space
on which a one-dimensional standard
Brownian motion $\{W(t)\}_{t\ge0}$ is
defined. Let $H$ be a Fr\'echet space. We
denote by $L_{\dbF}^2(0,T;H)$ the Fr\'echet
space consisting of all $H$-valued
$\dbF$-adapted processes
$X(\cd)$ such that
$\mathbb{E}(|X(\cd)|_{L^2(0,T;H)}^2)<\infty$;
by $L_{\dbF}^\i(0,T;H)$ the Fr\'echet space
consisting of all $H$-valued
$\dbF$-adapted bounded
processes; and by
$L_{\dbF}^2(\Om;C([0,T];H))$ the Fr\'echet
space consisting of all $H$-valued
$\dbF$-adapted continuous
processes $X$ such that
$\mathbb{E}(|X|_{C([0,T];H)}^2)<\i$. All of the above spaces
are equipped  with the canonical quasi-norms. Furthermore, all of the above spaces
are Banach spaces equipped  with the canonical  norms, if $H$ is a Banach space.

For simplicity, we  use the notation $
y_i\equiv y_i(x)=\frac{\pa y(x)}{\pa x_i}, $ where
$x_i$ is the $i$-th coordinate of a generic
point $x=(x_1,\cdots,x_n)$ in $\dbR^n$. In a
similar manner, we use the notation $ z_i$,
$v_i$, etc. for the partial derivatives of $z$
and $v$ with respect to $x_i$. Also, we denote
the scalar product in $\dbR^n$ by
$\lan\cd,\cd\ran$, and use $C$ to denote a
generic positive constant independent of the
solution $y$, which may change from line to
line.

Let $T>0$, $G \subset \dbR^{n}$ ($n\in \dbN$) be
a given bounded domain with the $C^4$ boundary
$\pa G$, and $\G$ be a given nonempty open
subset of $\pa G$.

Let
 $a_1\in
L^\infty_\dbF(0,T;L^\infty(G;\dbR^n))$,  $a_2\in
L^\infty_\dbF(0,T;L^\infty(G))$ and   $a_3\in
L^\infty_\dbF(0,T;W^{1,\infty}(G))$. 
We  assume $ f \in L^{2}_{\mathbb{F}}(0,T; L^{2}(G)) $, $g_{1}\in H^{1}(\Gamma)$, $g_{2}\in L^2_\dbF(0,T;L^2(\G))$, and
$a^{ij}: \;\Om\times [0,T]\times \cl{G}\to \dbR^{n\times
n}\;$ ($i,j=1,2,\cdots,n$) satisfies

\medskip
{\bf (H1)} {\it   $a^{ij}\in
L_{\dbF}^2(\Om;C^1([0,T];W^{2,\infty}(G)))$ and
$a^{ij}=a^{ji}$;}

 \ms

 {\bf (H2)} {\it 
%
%\begin{equation}\label{beta0}
$  \sum\limits^{n}_{i,j=1}a^{ij}(\om,t,x)\xi^{i}\xi^{j}
 \geq s_0|\xi|^{2},\ (\om,t,x,\xi)\equiv (\om,t,x,\xi^{1},\cdots,\xi^{n})
 \in \Omega\times (0,T)\t G\t \dbR^{n}$
% \end{equation}
 for some $s_0>0$.}
\ms

Now, the Cauchy problem of the forward stochastic parabolic differential equation can be described as follows.
\begin{equation}\label{system1}
\begin{cases}\ds
dy-\sum^{n}_{i,j=1}(a^{ij}y_i)_jdt=[\lan a_1,\n y\ran+ a_2 y + f]dt+a_3y dW(t)&\mbox{ in } (0,T)\times G,\\
\ns\ds y = g_{1} & \mbox{ on } (0,T)\times \G,
\\
\ns\ds  \frac{\pa y}{\pa\nu} =g_{2} & \mbox{ on } (0,T)\times \G.
\end{cases}
\end{equation}
\medskip

Let $\mathbf G_\G\=\{G'\subset G| \pa G'\cap \pa G\subset \G\}$ and
$$
H_{\G}^2(G)\deq \{\eta\in H^2_{loc}(G): \eta|_{G'}\in  H^2(G'),\q \forall G'\in \mathbf G_\G,\q \eta|_{\G}=g_{1},\; \pa_\nu \eta|_{\G}=g_{2}\}.
$$

The aim of this paper is to study
the Cauchy Problem with the Lateral Data for
the stochastic parabolic differential equation:  

\ms

{\bf Problem (CP)} Find
 a function  $y\in L_{\dbF}^2(\Om;C([0,T];L^2_{loc}(G)))\cap
L_{\dbF}^2(0,T;H_{\G}^2(G))$ 
  that satisfies system  \eqref{system1}. 
\ms

\begin{remark}
  It would be quite interesting to study the general case where $ g_{1} \in L^{2}_{\mathbb{F}}(0,T; H^{1}(\Gamma)) $. 
  However, this will lead to the term $ \mathbb{E} \int_{0}^{T} \int_{\Gamma} \sum\limits_{i,j=1}^{n} \theta \eta a^{i j} v_{i} \nu^{j}  d y d\Gamma $ in \eqref{th1eq4}, which we have no idea  how to handle.
\end{remark}

In certain applications involving diffusive, thermal, and heat transfer problems, measuring some boundary data can be challenging. For instance, in nuclear reactors and steel furnaces in the steel industry, the interior boundary can be difficult to measure. Similarly, the use of liquid crystal thermography in visualization can pose the same problem. To address this issue, engineers attempt to reconstruct the status of the boundary using measurements from the accessible boundary. This, in turn, creates the Cauchy problem for parabolic equations.  Due to the requirements in real applications, numerous researchers have focused on solving the Cauchy problem of parabolic differential equations, particularly on the inverse problems for deterministic parabolic equations, as seen in \cite{IK2000,Ksurvey,Klib2006,KT2007,KZ1998,LYZ2009CPAA,Y,YY2009}, and their respective references. Among these works, some have studied the identification of coefficients in parabolic equations through lateral Cauchy observations, such as uniqueness and stability \cite{IK2000,YY2009,Ksurvey,Y}, and numerical reconstruction \cite{KZ1998}, assuming that the initial value is known. Meanwhile, the determination of the initial value has been considered in \cite{Klib2006,KT2007,LYZ2009CPAA} under the assumption that all coefficients in the governing equation are known.

Stochastic parabolic equations have a wide range of applications for simulating various behaviors of stochastic models which are utilized in numerous fields such as random growth of bacteria populations, propagation of electric potential in neuron models, and physical systems that are subject to thermal fluctuations (e.g.,\cite{Kotelenez,LZ2021}).  In addition, they can be considered as simplified models for complex phenomena such as turbulence, intermittence, and large-scale structure (e.g.,\cite{Glimm}). Given the significant applications of these models, stochastic parabolic equations have been extensively studied in both physics and  mathematics. Therefore,   it is natural to study the Cauchy problem of stochastic parabolic equations in such situations.   However, due to the complexity of stochastic parabolic equations, some tools become ineffective for solving these problems. Thus, the research on inverse problems for stochastic parabolic differential equations is relatively scarce.  In \cite{barbu2}, the author proved backward uniqueness of solutions to stochastic semilinear parabolic equations and tamed Navier-Stokes equations driven by linearly multiplicative Gaussian noises, via the logarithmic convexity property known to hold for solutions to linear evolution equations in Hilbert spaces with self-adjoint principal parts. In \cite{GLWX2021IP,NHZ2020IP}, authors studied inverse random source problems for the stochastic time fractional diffusion equation driven by a spatial Gaussian random field, proving the uniqueness and representation for the inverse problems, as well as proposing a numerical method using Fourier methods and Tikhonov regularization. Carleman estimates play an important role in the study of inverse problems for stochastic parabolic differential equations, such as inverse source problems \cite{Lu1, Zhang2008}, determination of the history for stochastic diffusion processes \cite{DD2022IP,Lu1,wu,yuan2017}, and unique continuation properties \cite{Fu, yin2014}. We refer the reader to \cite{LZ2023} for a survey on some recent advances in Carleman estimates and their applications to inverse problems for stochastic partial differential equations.

In this paper, our objective is to  solve {\bf Problem (CP)}, i.e., we aim to retrieve the solution of equation \eqref{system1} with observed data from the lateral boundary. To this end, we first prove the conditional stability based on a new Carleman estimate for the stochastic parabolic equation \eqref{system1}. Then we construct a Tikhonov functional for the Cauchy problem based on the Tikhonov regularization strategy and prove the uniqueness of the minimizer of the Tikhonov functional, as well as the convergence rate for the optimization problem using variational principles, Riesz theorems, and Carleman estimates established previously.

Generally, the optimization problem for the Tikhonov functional is difficult to solve in the study of inverse problems in stochastic partial differential equations (SPDEs). This is because it involves solving the adjoint problem of the original problem, which is challenging to handle.  In fact, one of the primary differences between stochastic parabolic equations and deterministic parabolic equations is that at least one partial derivative of the solution does not exist, making it impossible to express the solution of that equation explicitly.    Fortunately, we can express the mild solution of the stochastic parabolic equations using the fundamental solution of the corresponding deterministic equation\cite{1986An,T1992}. This idea suggests that we can use kernel-based theory to numerically solve the minimization problem of the Tikhonov functional without computing the adjoint problem. Furthermore, we can solve the problem in one step without iteration, thus reducing the computational cost to some extent. This technique has gained attention in the study of numerical computation for ordinary and partial differential equations, and the use of fundamental solutions as kernels has proven effective for solving inverse problems in deterministic evolution partial differential equations. As far as we know, our work is the first attempt to apply the regularization method combining with kernel-based theory to solve inverse problems for stochastic parabolic differential equations.

\medskip
The strong solution of equation \eqref{system1} is useful in the proof of convergence rate of regularization method. Thus, we recall it here.

\begin{definition}\label{def}
We call $y\in L_{\dbF}^2(\Om;C([0,T];L^2_{loc}(G)))\cap
L_{\dbF}^2(0,T;H_{\G}^2(G))$ a solution to the
equation \eqref{system1} if for any $t\in[0,T]$ and a.e. $x\in G'\in \mathbf G_\G$, it holds
\begin{equation}\label{def eq1}
\begin{array}{ll}\ds
\q y(t,x) - y(0,x) \\
\ns \ds = \int_0^t \Big \{\sum^{n}_{i,j=1}\big(a^{ij}
(s,x)y_i(s,x)\big)_j  +[\lan a_1(s,x),
\n y(s,x)\ran+a_2(s,x)y(s,x)+f(s,x)y(s,x)] \Big \} ds\\
\ns \ds\q +\int_0^t a_3(s,x)y(s,x) dW(s),\qq
\dbP\mbox{-}\as
\end{array}
\end{equation}
\end{definition}
%

%
%\begin{definition}\label{milds}
%For any $t\in[0,T]$ and a.e. $x\in G'\in \mathbf G_\G$, the equation \eqref{system1}  possesses a unique mild solution $y\in L_{\dbF}^2(\Om;C([0,T];L^2_{loc}(G)))\cap
%L_{\dbF}^2(0,T;H_{\G}^2(G))$ %
%\begin{equation}\label{def eq1}
%\begin{array}{ll}\ds
%\q y(t,x) = \int_G S(t,0,x,\xi) y(0,\xi)d\xi +\int_0^t \int_G S(t,s,x,\xi)a_3(s,\xi)y(s,\xi)d\xi dW(s),\qq
%\dbP\mbox{-}\as
%\end{array}
%\end{equation}
%%  
%where $S$ is the strongly continuous two-parameter semi-group on $L^2_{loc}(G)$, which is generated by $\sum^{n}_{i,j=1}(a^{ij}
%y_i)_j  +[\lan a_1,\nabla y\ran+a_2y] $.
%\end{definition}
%%

It should be noted that the assumption regarding the solution, namely that $y\in L_{\dbF}^2(\Om;C([0,T];$ $L^2_{loc}(G)))\cap
L_{\dbF}^2(0,T;H_{\G}^2(G))$, implies a higher degree of smoothness than is strictly necessary for establishing H\"{o}lder stability of the Cauchy problem. However, this additional smoothness is required in order to facilitate the regularization process.

\ms

The remainder of this paper is organized as follows.   Section 2 presents a proof of a H\"{o}lder-type conditional stability result, along with the establishment of the Carleman estimate as a particularly useful tool in this proof. The regularization method, using the Tikhonov regularization strategy, is introduced in Section 3 and we showcase both the uniqueness and convergence rate of the regularization solution. Section 4 provides numerically simulated reconstructions aided by kernel-based learning theory.

%%%%%%%%%%%%%%%%%%%%%%%%%%%%%%%%%%%%%%%%%%%

%%%%%%%%%%%%%%%%%%%%%%%%%%%%%%%%%%%%%%%%%%%%%

\section{Conditional Stability}

In this section, we prove a stability estimate for the Cauchy problem.
\begin{theorem}\label{th1} (H\"{o}lder stability
estimate).  
For any given $G'\subset\subset G$ and
$\e \in \big(0,\frac{T}{2}\big)$,   there exists   $ \delta _{0} \in \left(
0,1\right) $, $\b\in (0,1)$ and a constant
 $C >0$ such
that if for  $\delta \in
\left( 0,\delta _{0}\right) $,
\begin{equation}\label{2.11_1}
  \max\{ 
    |f|_{ L^{2}_{\mathbb{F}}(0,T; L^{2}(G)) }, ~
    |g_{1}|_{H^{1}(\Gamma)}, ~
    \left|g_{2}\right|_{L_{\mathbb{F}}^2\left(0, T ; L^2(\Gamma)\right)}   
\}
\leq \delta,
\end{equation}
then
\begin{equation}\label{2.11_2}
\mE\int_{\e}^{T-\e}\int_{G'} \big( y^2 +  |\nabla
y|^2\big)dxdt\leq C (1 +  |y|^2_{L^2_\dbF((0,T);H^1(G))})\delta ^{2 \beta },
\quad  \forall \, \delta \in (0, \delta_{0}).
\end{equation}
\end{theorem}

We first recall the following exponentially
weighted energy identity, which will play a key
role in the sequel.

\begin{lemma}{\cite[Theorem 3.1]{Tang-Zhang1}}\label{Lemma1}
Let $n$ be a positive integer,
 \begin{equation}\label{1}
 b^{ij}=b^{ji}\in L_{\dbF}^2(\Om;C^1([0,T];W^{2,\infty}(\dbR^n))),\qq i,j=1,2,\cdots,n,
 \end{equation}
and $\ell\in C^{1,4}((0,T)\t\dbR^n),\ \Psi\in
C^{1,2}((0,T)\t\dbR^n)$. Assume $u$ is an
$H^2(\dbR^n)$-valued continuous
semi-martingale. Set
 $\th=e^{\ell }$ and $v=\th u$. Then for $\ae\2n$ $x\in \dbR^n$ and $\dbP$-$\as$
 $\om\in\Om$,
 \begin{align}\label{c1e2a}
  \notag
  & 2\int_0^T\th\[-\sum^{n}_{i,j=1} (b^{ij}v_i)_j+\cA v\]\[du-\sum^{n}_{i,j=1}(b^{ij}u_i)_jdt\]+2\int_0^T\sum^{n}_{i,j=1} (b^{ij}v_idv)_j\\
  \notag
  % \noalign{\ss}
  &\displaystyle\q+ 2\int_0^T\sum^{n}_{i,j=1}\[\sum_{i',j'}\(2b^{ij} b^{i'j'}\ell_{i'}v_iv_{j'}
  -b^{ij}b^{i'j'}\ell_iv_{i'}v_{j'}\)+\Psi
  b^{ij}v_iv- b^{ij}\(\cA\ell_i+\frac{\Psi_i}{2}\)v^2\]_jdt\\
  \notag
  &
  % \noalign{\ss}
  \displaystyle
  =2\int_0^T\sum^{n}_{i,j=1} \Big\{\sum_{i',j'}\[2b^{ij'}\(b^{i'j}\ell_{i'}\)_{j'} -
  \(b^{ij}b^{i'j'}\ell_{i'}\)_{j'}\]-\frac{b_t^{ij}}{2}+\Psi b^{ij}
  \Big\}v_iv_jdt\\
  \notag
  % \noalign{\ss}
  &
  \displaystyle\q
  +\int_0^T \cB v^2dt+2\int_0^T\[-\sum^{n}_{i,j=1} (b^{ij}v_i)_j+\cA v\]\[-\sum^{n}_{i,j=1} (b^{ij}v_i)_j+(\cA-\ell_t)v\]dt\\
  &
  % \noalign{\ss}
  \displaystyle\q+\(\sum^{n}_{i,j=1}
  b^{ij}v_iv_j+\cA v^2\)\Big|_0^T-\int_0^T\th^2\sum^{n}_{i,j=1}
  b^{ij}[(du_i+\ell_idu)(du_j+\ell_jdu)]-\int_0^T\th^2\cA(du)^2,
 \end{align}
where
\begin{equation}\label{c1e3a}
\left\{
 \begin{array}{ll}
 \ds \cA\=-\sum^{n}_{i,j=1}
 (b^{ij}\ell_i\ell_j-b^{ij}_j\ell_i
 -b^{ij}\ell_{ij})-\Psi,\\
  \ns
 \ds
 \cB\=2\[\cA\Psi-
 \sum^{n}_{i,j=1}(\cA b^{ij}\ell_i)_j\] -\cA_t-\sum^{n}_{i,j=1} (b^{ij}\Psi_j)_i.
  \ea
\right.
\end{equation}
\end{lemma}

In the sequel, for a positive integer $p$, denote by $O(\mu^p)$ a function of order
$\mu^p$ for large $\mu$ (which is independent
of $\l$); by $O_{\mu}(\l^p)$ a function of
order $\l^p$ for fixed $\mu$ and for large
$\l$.

\vspace{0.2cm}

\noindent{\it Proof of Theorem \ref{th1}.}\; We borrow
some idea from \cite{LL}. 

Take a bounded domain $J\subset
\dbR^n$  such that $\pa J \cap \overline G = \G$
and that $\widetilde G = J\cup G\cup\G$ enjoys a
$C^4$ boundary $\pa\widetilde G$. Then we have
\begin{equation}\label{wtG}
G\subset \widetilde G,\q \overline{\pa G \cap
\widetilde G} \subset \G,\q \pa G\setminus \G
\subset
 \pa\widetilde G \q\mbox{and }  \widetilde G \setminus G  \mbox{
contains some nonempty open subset. }
\end{equation}

 Let $G_0
\subset\subset \widetilde G \setminus G$ be an
open subdomain. We know that there is a $\psi\in
C^4(\widetilde G)$ satisfying (see \cite[Lemma
5.1]{Tang-Zhang1} for example)
\begin{equation}\label{psi}
\left\{
\begin{array}{ll}
\ds \psi > 0 &\mbox{ in } \widetilde G,\\
\ns\ds \psi = 0 & \mbox{ on } \pa \widetilde G,\\
\ns\ds |\nabla \psi| > 0 & \mbox{ in } G\subset \widetilde G\setminus G_0.
\end{array}
\right.
\end{equation}

Since $G'\subset\subset G$, we can choose a
sufficiently large $N>0$ such that
\begin{equation}\label{G'}
G'\subset \Big\{ x:\,x\in \widetilde
G,\;\psi(x)> \frac{4}{N}|\psi|_{L^\infty(\wt
G)}\Big\} \cap G.
\end{equation}
Further, let
$\rho=\frac{1}{\sqrt{2}}\(\frac{1}{2}-\k\) T>0$, there exists a positive number $c$ such that
\begin{equation}\label{c}
c\rho^2 < |\psi|_{L^\infty(\widetilde G)}  <
2c\rho^2.
\end{equation}
Then, define
\begin{equation}\label{phi}
\phi(t,x)=\psi(x)-c(t-t_0)^2,\q\a(t,x)=e^{\mu\phi(t,x)}
\end{equation}
 for a fix $t_0\in
[\sqrt{2}\rho,T-\sqrt{2}\rho]$, and denote
$\b_k=\b_k(\mu)=e^{\mu\big[\frac{k}{N}|\psi|_{L^\infty(\wt\Om)}
- \frac{c\rho^2}{N} \big]}, k=1,2,3,4$.    Set
\begin{equation}\label{Q1}
Q_k = \{(t,x):\,x\in \overline
G,\;\a(x,t)>\b_k\}, \q  k=1,2,3,4.
\end{equation}
Clearly, $Q_k$ is independent of $\mu$. Moreover, $\psi(x)>\frac{4}{N}|\psi|_{L^\infty(\wt G)}$ for any $(t,x)\in
\(t_0-\frac{\rho}{\sqrt{N}},t_0+\frac{\rho}{\sqrt{N}}\)\times
G'$, 
and thus
$$
\psi(x)-c(t-t_0)^2 >
\frac{4}{N}|\psi|_{L^\infty(\wt G)} -
\frac{c\rho^2}{N}.
$$
Hence, we see $\a(t,x)>\b_4$ in $(t,x)\in
\(t_0-\frac{\rho}{\sqrt{N}},t_0+\frac{\rho}{\sqrt{N}}\)\times
G'$, and thus $Q_4\supset\(t_0-\frac{\rho}{\sqrt{N}},t_0+\frac{\rho}{\sqrt{N}}\)\times
G'$. On the other
hand, for any $(t,x)\in Q_1$,
$$
\psi(x)-c(t-t_0)^2>\frac{1}{N}|\psi|_{L^\infty(\wt
G)} - \frac{c\rho^2}{N}.
$$
This yields
$$
|\psi|_{L^\infty(\wt G)} -
\frac{1}{N}|\psi|_{L^\infty(\wt G)} +
\frac{c\rho^2}{N}>c(t-t_0)^2.
$$
Together with \eqref{c} we have
$$
2\(1-\frac{1}{N}\)c\rho^2 + \frac{c\rho^2}{N} >
c(t-t_0)^2.
$$
Therefore, we conclude
\begin{equation}\label{Qk1}
\(t_0- \frac{\rho}{\sqrt{N}},t_0 +
\frac{\rho}{\sqrt{N}}\)\times G' \subset Q_1
\subset (t_0-\sqrt{2}\rho,t_0+\sqrt{2}\rho)\times
\overline G.
\end{equation}

Next,  for any $(t,x)\in\pa Q_1$, we know $x\in
\overline G$ and $\a(t,x)\geq\b_1$. In above set, if $x\in G$, $\a(t,x)=\b_1$, and if $x\in
\pa G$, it must hold $x\in \G$.
Indeed, if $x\in \pa G\setminus\G$, from $ \partial G\setminus\G\subset
\pa\wt G$ we get $\psi(x)=0$. On
the other hand, since $\a(t, x)\geq\b_1$, we know
$$
\psi(x)-c(t-t_0)^2 = -c(t-t_0)^2\geq
\frac{1}{N}|\psi|_{L^\infty(\wt G)} -
\frac{c\rho^2}{N}.
$$
And thus
$$
0\leq c(t-t_0)^2\leq \frac{1}{N}(c\rho^2 -
|\psi|_{L^\infty(\wt G)}),
$$
which contradicts \eqref{c}. Therefore, we have
\begin{equation}\label{C1}
\pa Q_1 = \Si_1\cup\Si_2,
\end{equation}
where
\begin{equation}\label{C2}
\begin{array}{ll}\ds
\Si_1 \subset [0,T]\times \G,\qq \Si_2 = \{(t,x) \in [0,T] \times G:\, x\in
G,\,\a(t,x)=\b_1\}.
\end{array}
\end{equation}

Let $\eta\in C_0^\infty(Q_2)$ such that
 $0\leq \eta
\leq 1$ and $\eta = 1$ in $Q_3$. For any $y$
solving \eqref{system1}, let $z = \eta y$, then
$z$ solves
\begin{equation}\label{system2}
\left\{
\begin{array}{ll}\ds
dz - \sum^{n}_{i,j=1}(a^{ij}z_i)_j dt = \big( \lan a_1, \nabla z \ran + a_2 z + \tilde{f} + \eta f\big)dt + a_3 z dB(t) &\mbox{ in } Q_1,\\
\ns\ds z=\frac{\pa z}{\pa\nu}=0 &\mbox{ on }
\Si_2.
\end{array}
\right.
\end{equation}
Here $\ds \tilde{f} = -\sum^{n}_{i,j=1} (a^{ij}_j\eta_i y +
a^{ij}\eta_{ij}y + 2a^{ij}\eta_i y_j) - \lan
a_1,\nabla \eta\ran y + \eta_t y  \in
L^2_{\mathbb{F}}(0,T;L^2(G))$ and $\tilde{f}$ is supported in
$Q_2\setminus Q_3$.

Applying  Lemma \ref{Lemma1} to
\eqref{system2} with $u=z$, $b^{ij}=a^{ij}$,
$\ell = \l\a$ and
$\ds\Psi=2\sum^{n}_{i,j=1}a^{ij}\ell_{ij}$,
integrating it in $G$ and taking mean value,
noting that $z$ is supported in $Q_1$, we see
\begin{equation}\label{th1eq0}
\begin{array}{ll}
\displaystyle
 2\mE\int_{Q_1} \th\[-\sum^{n}_{i,j=1} (a^{ij}v_i)_j+\cA v\]\[dz-\sum^{n}_{i,j=1}(a^{ij}z_i)_jdt\]dx +2\int_{Q_1}\sum^{n}_{i,j=1} (a^{ij}v_idv)_jdx \\
\ns \ds\q+
2\mE\int_{Q_1}\sum^{n}_{i,j=1}\[\sum_{i',j'}\(2a^{ij}
a^{i'j'}\ell_{i'}v_iv_{j'}
-a^{ij}a^{i'j'}\ell_iv_{i'}v_{j'}\)+\Psi
a^{ij}v_iv- a^{ij}\(\cA\ell_i+\frac{\Psi_i}{2}\)v^2\]_jdxdt\\
\ns \ds \ge 2\sum^{n}_{i,j=1}\mE\int_{Q_1}
c^{ij}v_iv_j dxdt
+\int_0^T \cB v^2 dxdt+\mE\int_{Q_1} \Big|-\sum^{n}_{i,j=1} (a^{ij}v_i)_j+\cA v\Big|^2dxdt\\
\noalign{\ss}
\displaystyle\q-\mE\int_{Q_1}\th^2\sum^{n}_{i,j=1}
a^{ij}(dz_i+\ell_idz)(dz_j+\ell_jdz)dx-\mE\int_{Q_1}\th^2\cA(dz)^2dx,
\end{array}
\end{equation}
where
\begin{equation}\label{th1eq01}
\left\{
 \begin{array}{ll}
 \ds \cA =-\sum^{n}_{i,j=1}
 \big(a^{ij}\ell_i\ell_j-a^{ij}_j\ell_i
 -a^{ij}\ell_{ij}\big)-\Psi,\\
  \ns
 \ds
 \cB =2\[A\Psi-
 \sum^{n}_{i,j=1}\big(Aa^{ij}\ell_i\big)_j\] -A_t-\sum^{n}_{i,j=1} \big(a^{ij}\Psi_j\big)_i-\ell_t^2,\\
 \ns
 \ds c^{ij} =\sum_{i',j'}\[2a^{ij'}\big(a^{i'j}\ell_{i'}\big)_{j'} -
 \big(a^{ij}a^{i'j'}\ell_{i'}\big)_{j'}\]-\frac{a_t^{ij}}{2}+\Psi a^{ij}.
  \end{array}
\right.
\end{equation}
Now, we estimate $\cA$, $\cB$ and $c^{ij}$.
\begin{equation}\label{A}
\begin{array}{ll}\ds
\cA \!\!\!&\ds =-\sum^{n}_{i,j=1}
 \big(a^{ij}\ell_i\ell_j-a^{ij}_j\ell_i
 -a^{ij}\ell_{ij}\big)-\Psi\\
 \ns& \ds = -\l^2\mu^2\a^2\sum^{n}_{i,j=1}a^{ij}\psi_i\psi_j
  + \l\mu\a\sum^{n}_{i,j=1}a^{ij}_j\psi_i -
  \l\mu^2\a\sum^{n}_{i,j=1}a^{ij}\psi_i\psi_j-\l\mu\a\sum^{n}_{i,j=1}a^{ij}\psi_{ij}\\
  \ns&\ds =-\l^2\mu^2\a^2\sum^{n}_{i,j=1}a^{ij}\psi_i\psi_j
  +\l\a O(\mu^2),
\end{array}
\end{equation}
In order to estimate $\cB$, we first do some
computations. For $\Psi$, we have
\begin{equation}\label{Psi}
\Psi=2\sum^{n}_{i,j=1}a^{ij}\big(\l\mu^2\a
\psi_i\psi_j+\l\mu\a \psi_{ij}\big) =2\l\mu^2\a
\sum^{n}_{i,j=1}a^{ij}\psi_i\psi_j+\l\a O(\mu).
\end{equation}
Next, we have
$$
\ell_{i'j'j}=\l\mu^3\a\psi_{i'}\psi_{j'}\psi_j+\l\a
O(\mu^2),\q
\ell_{i'j'ij}=\l\mu^4\a\psi_{i'}\psi_{j'}\psi_i\psi_j+\l\a
O(\mu^3).
$$
Therefore, we get
$$
\Psi_j=2\sum_{i',j'=1}^n\big(a^{i'j'}\ell_{i'j'}\big)_j=2\sum_{i',j'}\big(a^{i'j'}_j\ell_{i'j'}+a^{i'j'}\ell_{i'j'j}\big)
=2\l\mu^3\a\sum_{i',j'=1}^na^{i'j'}\psi_{i'}\psi_{j'}\psi_j+\l\a
O(\mu^2),
$$
and
$$
\Psi_{ij}=2\sum_{i',j'=1}^n\big(a^{i'j'}_{ij}\ell_{i'j'}+a^{i'j'}\ell_{i'j'ij}+2a^{i'j'}_j\ell_{i'j'i}\big)
=2\l\mu^4\a\sum_{i',j'=1}^n a^{i'j'}\psi_{i'}\psi_{j'}\psi_i\psi_j+\l\a
O(\mu^3).
$$
Hence, we find
\begin{equation}\label{th1eq001}
-\sum^{n}_{i,j=1}\big(a^{ij}\Psi_j\big)_i=-\sum^{n}_{i,j=1}\big(a^{ij}_i\Psi_j+a^{ij}\Psi_{ji}\big)
=-2\l\mu^4\a\(\sum^{n}_{i,j=1}a^{ij}\psi_i\psi_j\)^2
+\l\a O(\mu^3).
\end{equation}
Further, from \eqref{A} and \eqref{Psi}, we have
\begin{equation}\label{th1eq02}
\cA\Psi=-2\l^3\mu^4\a^3\(\sum^{n}_{i,j=1}a^{ij}\psi_i\psi_j\)^2+\l^3\a^3O(\mu^3)+\l^2\a^2O(\mu^4).
\end{equation}
From the definition of $A$, we find
$$
\begin{array}{ll}\ds
\cA_j\!\!\!&\ds
=-\sum_{i',j'=1}^n\big(a_j^{i'j'}\ell_{i'}\ell_{j'}+2a^{i'j'}\ell_{i'}\ell_{j'j}-a_{j'j}^{i'j'}\ell_{i'}
-a_{j'}^{i'j'}\ell_{i'j}
+a_j^{i'j'}\ell_{i'j'}+a^{i'j'}\ell_{i'j'j}\big)\\
\ns & \ds
=-\sum_{i',j'=1}^n\big(2a^{i'j'}\ell_{i'}\ell_{j'j}+a^{i'j'}\ell_{i'j'j}\big)+\big(\l\a+\l^2\a^2
\big)O(\mu^2)\\
\ns & \ds
=-2\l^2\mu^3\a^2\sum_{i',j'=1}^na^{i'j'}\psi_{i'}\psi_{j'}\psi_j+
O_\mu(\l)+\l^2\a^2 O(\mu^2).
\end{array}
$$
Hence, we see
$$
\sum^{n}_{i,j=1}\cA_ja^{ij}\ell_i=-2\l^3\mu^4\a^3\Big(\sum^{n}_{i,j=1}a^{ij}\psi_{i}\psi_{j}\Big)^2+
O_\mu(\l^2)+\l^3\a^3 O(\mu^3),
$$
which leads to
\begin{equation}\label{th1eq03}
\ba{ll} \ds
\sum^{n}_{i,j=1}\big(\cA a^{ij}\ell_i\big)_j \3n & \ds=\sum^{n}_{i,j=1}\cA_ja^{ij}\ell_i+\cA\sum^{n}_{i,j=1}\big(a_j^{ij}\ell_i+a^{ij}\ell_{ij}\big)\\
\ns &\ds
=-3\l^3\mu^4\a^3\Big(\sum^{n}_{i,j=1}a^{ij}\psi_{i}\psi_{j}\Big)^2
+O_\mu(\l^2) +\l^3\a^3 O(\mu^3). \ea
\end{equation}
Further, we have
\begin{equation}\label{th1eq04}
\begin{array}{ll} \ds
\cA_t\3n&\ds =-\sum^{n}_{i,j=1}
\(a^{ij}\ell_i\ell_j-a^{ij}_j\ell_i
+a^{ij}\ell_{ij}\)_t\\
\ns  & \ds=-\sum^{n}_{i,j=1}
\[a^{ij}(\ell_i\ell_j)_t-a^{ij}_j\ell_{it}
+a^{ij}\ell_{ijt}\]+\l^2\a^2 O(\mu^2)+\l\a O(\mu^2)+\l\a TO(\mu^3)\\
\ns & \ds  =\l^2\a^2 T O(\mu^3) +\l^2\a^2
O(\mu^2)+\l\a O(\mu^2)+\l\a TO(\mu^3).
\end{array}
\end{equation}
Finally, it holds
\begin{equation}\label{th1eq05}
\ell_t^2=\l^2\mu^2\a^2\phi_t^2= O_\mu(\l^2).
\end{equation}

From the definition of $B$ (see
\eqref{th1eq01}), and combing
\eqref{th1eq001}--\eqref{th1eq5}, we have
$$
\ba{ll} \cB\3n&=\ds
-4\l^3\mu^4\a^3\(\sum^{n}_{i,j=1}a^{ij}\psi_i\psi_j\)^2
+6\l^3\mu^4\a^3\(\sum^{n}_{i,j=1}a^{ij}\psi_{i}\psi_{j}\)^2
-2\l\mu^4\a\(\sum^{n}_{i,j=1}a^{ij}\psi_i\psi_j\)^2\\
\ns&\ds\q + \l^3\a^3 O(\mu^3) + O_\mu(\l^2)
  \\
\ns &=\ds 2\l^3\mu^4\a ^3\(\sum^{n}_{i,j=1}
a^{ij}\psi_i\psi_j\)^2+\l^3\a^3O(\mu^3) +
O_\mu(\l^2).
\end{array}
$$
Hence, we know
\begin{equation}\label{B}
\cB\ge 2s_0^2\l^3\mu^4\a ^3|\n \psi|^4
+\l^3\a^3O(\mu^3) + O_\mu(\l^2).
\end{equation}
Now we estimate $c^{ij}$. By direct computation,
we have
\begin{equation}\label{cij}
\3n\begin{array}{ll} \ds \sum^{n}_{i,j=1}
c^{ij}v_iv_j\3n&
\ds=\sum^{n}_{i,j=1}\Big\{\sum_{i',j'=1}^n\[2a^{ij'}a^{i'j}\ell_{i'j'}
+
a^{ij}a^{i'j'}\ell_{i'j'}+2a^{ij'}a_{j'}^{i'j}\ell_{i'}
-(a^{ij}a^{i'j'})_{j'}\ell_{i'}\]-\frac{a_t^{ij}}{2}\Big\}v_iv_j\\
\ns &\ds
=\sum^{n}_{i,j=1}\Big\{\sum_{i',j'=1}^n\[2\l\mu^2\a
a^{ij'}a^{i'j}\psi_{i'}\psi_{j'}  + \l\mu^2\a
a^{ij}a^{i'j'}\psi_{i'}\psi_{j'} + \l\a  O(\mu)\]+O(1)\Big\}v_iv_j\\
\ns &\ds=2\l\mu^2\a
\(\sum^{n}_{i,j=1}a^{ij}\psi_iv_j\)^2\!+\!\l\mu^2\a
\(\sum^{n}_{i,j=1}
a^{ij}\psi_i\psi_j\)\(\sum^{n}_{i,j=1}a^{ij}v_iv_j\)\\
\ns &\ds\q + \l\a |\n v|^2 O(\mu)+O(1)|\n v|^2\\
\ns &\ds\ge [s_0^2\l\mu^2\a |\n\psi|^2+\l\a
O(\mu)+O(1)]|\n v|^2.
\end{array}
\end{equation}

Owing to that $z$ is a solution to equation \eqref{system2}, we find
\begin{equation}\label{th1eq2}
\begin{array}{ll}\ds
\q 2\mE\int_{Q_1}\th\[-\sum^{n}_{i,j=1} \big(a^{ij}v_i\big)_j+\cA v\]\[dz-\sum^{n}_{i,j=1}\big(a^{ij}z_i\big)_jdt\]dx\\
\ns\ds   \leq \mE\int_{Q_1} \[-\sum^{n}_{i,j=1} \big(a^{ij}v_i\big)_j+\cA v\]^2dxdt + \mE\int_{Q_1} \th^2\( \lan a_1, \nabla z \ran + a_2 z + \tilde{f} + \eta f  \)^2dxdt \\
\ns\ds   \leq \mE\int_{Q_1} \[-\sum^{n}_{i,j=1} \big(a^{ij}v_i\big)_j+\cA v\]^2dxdt + 3 |a_1|^2_{L^\infty(0,T;L^\infty(G;\dbR^n))}\mE\int_{Q_1}\th^2 |\nabla z|^2 dxdt  \\
\ns\ds \q   +  3
|a_2|^2_{L^\infty(0,T;L^\infty(G))}\mE\int_{Q_1}\th^2
z^2 dxdt + 6\mE\int_{Q_1}\th^2 (\tilde{f}^2 + \eta^{2} f^{2}) dxdt.
\end{array}
\end{equation}
It is clear that $\overline{Q_{2}} \cap \overline{\Sigma_{2}} = \emptyset$. Hence, $\eta = 0$ on $\Sigma_{2}$.
Recall that $ y = g_{1}  $ on $ \Sigma_{1} $ and $ g_{1} \in H^{1}(\Gamma)$.
Then it holds that 
\begin{equation}
  \label{th1eq4}
  \begin{aligned}
    \mE\int_{Q_1}\sum^{n}_{i,j=1}
  \big(a^{ij}v_idv\big)_jdx 
  & = 
  \mathbb{E} \int_{\Sigma_{1}} \sum_{i,j=1}^{n} a^{i j} v_{i} \nu^{i} d v d \Gamma
  \\
  &
  \leq C \lambda  \mu \mathbb{E} \int_{0}^{T} \int_{\Gamma} \alpha  \theta^{2} (|g_{2}|^{2} + |\nabla_{\Gamma} g_{1}|^{2} + \lambda^{2} \mu^{2} \alpha^{2} |g_{1}|^{2}) d\Gamma d t
  ,
  \end{aligned}
\end{equation}
where $ \nabla_{\Gamma} $ denotes the tangential gradient on $ \Gamma $.

By means of  $z=\frac{\pa
z}{\pa\nu}=0$ on $\Si_2$, we find
\begin{equation}\label{th1eq5}
\begin{array}{ll}\ds
\Big|\mE\int_{Q_1}\sum^{n}_{i,j=1}\[\sum_{i',j'=1}^{n}\(2a^{ij}
a^{i'j'}\ell_{i'}v_iv_{j'}
-a^{ij}a^{i'j'}\ell_iv_{i'}v_{j'}\)+\Psi
a^{ij}v_iv- a^{ij}\(\cA\ell_i+\frac{\Psi_i}{2}\)v^2\]_jdxdt\Big| \\
 \ns\ds 
 \leq C
 \lambda \mu \mathbb{E} \int_{0}^{T} \int_{\Gamma} \alpha \theta^{2} (|g_{2}|^{2} + |\nabla_{\Gamma} g_{1}|^{2} + \lambda^{2} \mu^{2} \alpha^{2} |g_{1}|^{2}) d \Gamma d t.
 \end{array}
\end{equation}
From \eqref{cij} and \eqref{B}, we know that
there is a $\mu_0>0$ such that for every
$\mu\geq\mu_0$, one can find a $\l_0(\mu)>0$ so
that for all $\l\geq\l_0(\mu)$, it holds that
\begin{equation}\label{th1eq6}
\mE\sum^{n}_{i,j=1}\int_{Q_1} c^{ij}v_iv_j dxdt \geq
C\l\mu^2 \mE \int_{Q_1}\a|\n
 v|^2dxdt,
\end{equation}
and
\begin{equation}\label{th1eq7}
\mE\int_{Q_1} Bv^2 dxdt \geq C\l^3\mu^4
\mE\int_{Q_1}\a^3 v^2 dxdt.
\end{equation}
Utilizing the fact that $z$ solves
\eqref{system2} again, we get
\begin{equation}\label{th1eq8}
\begin{array}{ll}\ds
\q\mE\int_{Q_1}\th^2\sum^{n}_{i,j=1}
 a^{ij}(dz_i+\ell_idz)(dz_j+\ell_jdz) dx \\
 \ns\ds  =  \mE\int_{Q_1}\th^2\sum^{n}_{i,j=1}
 a^{ij} \[(a_3 z)_i (a_3 z)_j + 2 \l\mu\a\psi_i (a_3 z)_j a_3 z + \l^2\mu^2\a^2 \psi_i\psi_j a_3^2 z^2  \]dxdt \\
 \ns\ds \leq C|a_3|^2_{L^\infty_\dbF(0,T;W^{1,\infty}(G))}\(\mE\int_{Q_1}\th^2\big(z^2 + |\nabla z|^2 \big)dxdt +  \l^2\mu^2 \mE\int_{Q_1} \th^2 \a^2|z|^2dxdt \).
 \end{array}
\end{equation}
By $ z_i=\th^{-1} (v_i-\ell_i v)=\th^{-1}
(v_i-\l\mu\a\psi_i v)$ and $v_i=\th (z_i+\ell_i
z)=\th (z_i+\l\mu\a\psi_i z)$, we get
\begin{equation}\label{th1eq9}
\frac{1}{C}\th^2\big(|\n
z|^2+\l^2\mu^2\a^2z^2\big)\le |\n
v|^2+\l^2\mu^2\a^2v^2\le C\th^2\big(|\n
z|^2+\l^2\mu^2\a^2z^2\big).
\end{equation}

Combing \eqref{th1eq0},
\eqref{th1eq2}--\eqref{th1eq9}, we know that
there is a $\mu_1\geq \mu_0$ such that for all
$\mu\geq\mu_1$, we can find a $\l_1(\mu_1)\geq
\l_0(\mu_0)$ so that for every $\l\geq
\l_1(\mu_1)$, it holds that
\begin{equation}\label{th1eq10}
\begin{array}{ll}\ds
\q \l^3\mu^4 \mE\int_{Q_1} \a^3 \theta^2z^2 dxdt + \l\mu^2  \mE\int_{Q_1}\a\theta^2|\nabla z|^2 dxdt \\
\ns\ds \leq C \mE\int_{Q_1} \theta^2 (\tilde{f}^2 + f^{2}) dxdt +
C
 \lambda \mu \mathbb{E} \int_{0}^{T} \int_{\Gamma} \alpha \theta^{2} |g_{2}|^{2} d \Gamma d t + C \lambda \mu  \int_{\Gamma} \alpha \theta^{2}  (|\nabla_{\Gamma}  g_{1}|^{2} + \lambda^{2} \mu^{2} \alpha^{2} |g_{1}|^{2}) d \Gamma.
\end{array}
\end{equation}
Recalling that $\ds \tilde{f} = -\sum^{n}_{i,j=1}
(a^{ij}_j\eta_i y + a^{ij}\eta_{ij}y +
2a^{ij}\eta_i y_j) - \lan a_1,\nabla \eta\ran y
+ \eta_t y $, from  \eqref{th1eq10}, we find
that
\begin{equation}\label{th1eq11}
\begin{array}{ll}\ds
\q\l^3\mu^4 \mE\int_{Q_1} \a^3 \theta^2z^2 dxdt + \l\mu^2  \mE\int_{Q_1}\a\theta^2|\nabla z|^2 dxdt\\
\ns\ds   \leq C \mE\int_{Q_2\setminus Q_3}
\theta^2 \big(|y|^2 + |\nabla y|^2 \big) dxdt +
C
 \lambda \mu \mathbb{E} \int_{0}^{T} \int_{\Gamma} \alpha \theta^{2} |g_{2}|^{2} d \Gamma d t 
 \\
 \ns\ds \quad 
 + C \lambda \mu  \int_{\Gamma} \alpha \theta^{2}  (|\nabla_{\Gamma}  g_{1}|^{2} + \lambda^{2} \mu^{2} \alpha^{2} |g_{1}|^{2}) d \Gamma
 + C \mathbb{E} \int_{Q} \theta^{2} f^{2} d x d t
 ,
\end{array}
\end{equation}
which, together with $z=y$ in $Q_3$, implies that
\begin{equation}\label{th1eq12}
\begin{array}{ll}\ds
\q \l^3\mu^4 \mE\int_{Q_3} \a^3\theta^2 y^2 dxdt + \l\mu^2  \mE\int_{Q_3} \a\theta^2|\nabla y|^2 dxdt \\
\ns\ds  \leq C \mE\int_{Q_2\setminus Q_3}
\theta^2 \big(|y|^2 + |\nabla y|^2 \big) dxdt +
C
 \lambda \mu \mathbb{E} \int_{0}^{T} \int_{\Gamma} \alpha \theta^{2} |g_{2}|^{2} d \Gamma d t 
 \\
 \ns\ds \quad 
 + C \lambda \mu  \int_{\Gamma} \alpha \theta^{2}  (|\nabla_{\Gamma}  g_{1}|^{2} + \lambda^{2} \mu^{2} \alpha^{2} |g_{1}|^{2}) d \Gamma
 + C \mathbb{E} \int_{Q} \theta^{2} f^{2} d x d t
\end{array}
\end{equation}
From the choice of $Q_3$ and \eqref{G'}, we
know that
\begin{equation}\label{th1eq14.1}
\begin{array}{ll}\ds
\q\l^3\mu^4 \mE\int_{Q_3} \a^3\theta^2 y^2 dxdt
+ \l\mu^2  \mE\int_{Q_3} \a\theta^2|\nabla y|^2
dxdt\\
\ns\ds \geq  e^{2\l\b_4}
\mE\int_{t_0-\frac{\rho}{\sqrt{N}}}^{t_0+\frac{\rho}{\sqrt{N}}}\int_{G\,'}
\big(\l^3 y^2 +  \l|\nabla y|^2\big)dxdt.
\end{array}
\end{equation}

From the definition of $Q_2$ and $Q_3$, we find
\begin{equation}\label{th1eq14.2}
\begin{array}{ll}\ds
\mE\int_{Q_2\setminus Q_3}
\theta^2 \big(|y|^2 + |\nabla y|^2 \big) dxdt \leq e^{2\l\b_3} \mE\int_{Q_1}
\big(|y|^2 + |\nabla y|^2 \big) dxdt.
\end{array}
\end{equation}
Now, we fix $\mu=\mu_1$ in $\b_3$, $\b_4$ and
$\th$. Let $\g=\max_{(t,x)\in \overline{G}\times (0,T)}\a(t,x)$. From \eqref{th1eq12}--\eqref{th1eq14.2},
we obtain that
$$
\begin{array}{ll}\ds
\q e^{2\l\b_4}
\mE\int_{t_0-\frac{\rho}{\sqrt{N}}}^{t_0+\frac{\rho}{\sqrt{N}}}\int_{G'}
\big(\l^3 y^2 +  \l|\nabla y|^2\big)dxdt\\
\ns\ds \leq Ce^{2\l\b_3} \mE\int_{Q_1}
\big(|y|^2 + |\nabla y|^2 \big) dxdt + Ce^{2\l\g}
\big (
  |f|_{ L^{2}_{\mathbb{F}}(0,T; L^{2}(G)) }^{2} +
  |g_{1}|_{H^{1}(\Gamma)}^{2} + 
  \left|g_{2}\right|_{L_{\mathbb{F}}^2\left(0, T ; L^2(\Gamma)\right)}^{2}  
  \big)
  ,
\end{array}
$$
which implies that for all $\l\geq \l_1$, it
holds that
\begin{equation}\label{th1eq14.3}
\begin{array}{ll}\ds
\q\mE\int_{t_0-\frac{\rho}{\sqrt{N}}}^{t_0+\frac{\rho}{\sqrt{N}}}\int_{G\,'}
\big( y^2 + |\nabla y|^2\big)dxdt \\
\ns\ds \leq  Ce^{-2\l(\b_4-\b_3)} \mE\int_{Q_1}
\big(|y|^2 + |\nabla y|^2  \big) dxdt + Ce^{2\l\g}
\big (
  |f|_{ L^{2}_{\mathbb{F}}(0,T; L^{2}(G)) }^{2} +
  |g_{1}|_{H^{1}(\Gamma)}^{2} + 
  \left|g_{2}\right|_{L_{\mathbb{F}}^2\left(0, T ; L^2(\Gamma)\right)}^{2}  
  \big)
  .
\end{array}
\end{equation}
By taking $t_0
= \sqrt{2}\rho + \frac{2i\rho}{\sqrt{N}}$,
$i=0,1,\cds,m$ such that
$$
\sqrt{2}\rho + \frac{2m\rho}{\sqrt{N}}\leq
T-\sqrt{2}\rho \leq \sqrt{2}\rho +
\frac{(2m+1)\rho}{\sqrt{N}},
$$
and $\e=\sqrt{2}\rho$,
we get
\begin{align}\label{C3}\notag
  & \q \mE\int_{\e}^{T-\e}\int_{G'} \big( y^2 + |\nabla
  y|^2\big)dxdt\\\notag
  & =\mE\int_{\sqrt{2}\rho
  }^{T-\sqrt{2}\rho}\int_{G'}
  \big( y^2 + |\nabla y|^2\big)dxdt \\\notag
  &  \leq \sum_{i=1}^m\mE\int_{\sqrt{2}\rho +
  \frac{(i-1)\rho}{\sqrt{N}}}^{\sqrt{2}\rho +
  \frac{(i+1)\rho}{\sqrt{N}}}\int_{G'} \big( y^2
  + |\nabla y|^2\big)dxdt
  \\ 
  &  \leq  Ce^{-2\l(\b_4-\b_3)} \mE\int_{Q_1}
  \big(|y|^2 + |\nabla y|^2  \big) dxdt +Ce^{2\l\g}
  \big (
    |f|_{ L^{2}_{\mathbb{F}}(0,T; L^{2}(G)) }^{2} +
    |g_{1}|_{H^{1}(\Gamma)}^{2} + 
    \left|g_{2}\right|_{L_{\mathbb{F}}^2\left(0, T ; L^2(\Gamma)\right)}^{2}  
    \big)
    .
\end{align}
We now balance the terms in the right hand side of \eqref{C3} via choosing $\l= \l(\d)$ such
that $ e^{\l\g}\d =Ce^{-\l(\b_4-\b_3)}$. This implies that
$$
\l = \frac{-\ln \d}{\g+\b_4-\b_3}.
$$
Hence, for $\d \in (0, \d_0)$, where the number $\d_0 $ is so small
that
$\frac{-\ln \d}{\g+\b_4-\b_3}\geq \l_2$, we have  \eqref{2.11_2} with $\b=\frac{\b_4-\b_3}{\g+\b_4-\b_3}$.
\endpf

\begin{remark}
The inequality \eqref{th1eq10} is called the Carleman estimate for \eqref{system1}.
\end{remark}

\section{Regularization Method}
\label{sec:4}

For $u\in L_{\dbF}^2(\Om;C([0,T];L^2(G)))\cap
L_{\dbF}^2(0,T;H^2(G))$,  set
$$
\begin{array}{ll}\ds
(Pu)(t,x) \3n&\ds \= u(t,x) - u(0,x)\\
\ns&\ds \q - \int_0^t \Big \{\sum^{n}_{i,j=1}\big(a^{ij}
(s,x)u_i(s,x)\big)_j  +[\lan a_1(s,x),
\n u(s,x)\ran+a_2(s,x)u(s,x)] \Big \} ds\\
\ns&\ds \q -\int_0^t a_3(s,x)u(s,x) dW(s),\qq
\dbP\mbox{-}\as
\end{array}
$$

Let
\begin{align*}
  \mathcal{H}  \deq \{ u \in L^{2}_{\mathbb{F}}(0,T;H^{2}(G)) 
  & : 
  P u  \in L^{2}_{\mathbb{F}}(\Omega; H^{1}(0,T;L^{2}(G))), 
  \\
  & \quad  
  u|_{(0,T) \times \Gamma} = g_{1}, ~ \partial_{\nu} u |_{ (0,T) \times \Gamma} = g_{2} \}
  .
\end{align*}

Denote $ \hat{f}(t) = \int_{0}^{t} f(s) d s $. Clearly, we have 
\begin{align}
  \label{eqfEstimate}
  |\hat{f}|_{L^{2}_{\mathbb{F}}(\Omega;H^{1}(0,T;L^{2}(G)))} \leq C |f|_{L^{2}_{\mathbb{F}}(0,T;L^{2}(G))}
  .
\end{align}
Given a function $F\in  \mathcal{H}$, the
Tikhonov functional is now constructed as
\begin{eqnarray}\label{4.10}
J_{\gamma }\left( u\right)\3n &=\3n&|
Pu - \hat{f}|_{L^{2}_{\mathbb{F}}(\Omega;H^{1}(0,T;L^{2}(G)))}^{2}+\gamma |u-F|
_{L_{\dbF}^2(0,T;H^2(G))}^{2},
\end{eqnarray}
where $ u \in \mathcal{H} $ and $ \gamma \in (0,1) $.

We have the following result.
\begin{theorem}\label{th4.1} 
For every  $\gamma \in \left( 0,1\right) $, 
there exists a unique minimizer $u_{\gamma }\in
\mathcal{H}$ of the functional  $J_{\gamma}\left( u\right) $ in \eqref{4.10}.
% and with a constant $C 
% >0$  the following
% estimate holds 
% %
% \begin{equation}\label{4.12}
% |u_{\gamma }|_{L_{\dbF}^2(0,T;H^2(G))}\leq
% \frac{C}{\sqrt{\gamma }} |
% F|_{L_{\dbF}^2(0,T;H^2(G))
% }.  
% \end{equation}
%
\end{theorem}
\textbf{Proof}.  Let
\begin{align*}
\mathcal{H}_{0} \=\{ u \in L^{2}_{\mathbb{F}}(0,T;H^{2}(G)) 
& : 
P u  \in L^{2}_{\mathbb{F}}(\Omega; H^{1}(0,T;L^{2}(G))), 
\\
& \quad  
u|_{(0,T) \times \Gamma} = 0, ~ \partial_{\nu} u |_{ (0,T) \times \Gamma} = 0 \}
.
\end{align*}
Fix $ \gamma \in (0,1) $.
We define the inner products as follows:
\begin{align}
    \label{eqDou20231}
    \langle  \varphi, \psi \rangle_{\mathcal{H}_{0}} = 
    \langle P \varphi, P \psi \rangle_{L^{2}_{\mathbb{F}}(\Omega; H^{1}(0,T;L^{2}(G)))} 
    + 
    \gamma \langle \varphi, \psi \rangle_{L^{2}_{\mathbb{F}}(0,T;H^{2}(G))} 
\end{align}
where $ \varphi, \psi \in \mathcal{H}_{0} $. 
Let $\overline{\mathcal{H}}_{0}$ be the completion of $\mathcal{H}_{0}$ with respect to the inner product $\langle \cdot, \cdot \rangle_{\mathcal{H}_{0}}$, and still denoted by $\mathcal{H}_{0}$ for the sake of simplicity.

For $ u \in \mathcal{H} $, let $v=u-F$. Then $v\in \mathcal{H}_{0}$. 
By \eqref{4.10}, we
should minimize the following
functional 
\begin{equation}\label{2.108}
\overline{J}_{\gamma }\left( v\right)
=|Pv+ PF - \hat{f}|_{L^{2}_{\mathbb{F}}(\Omega; H^{1}(0,T;L^{2}(G)))}^{2}+\gamma
|v|_{L_{\dbF}^2(0,T;H^2(G))}^{2},\q v\in \mathcal{H}_{0}.  
\end{equation}
If $v_{\gamma }\in  \mathcal{H}_{0}$ is a minimizer of the functional
\eqref{2.108}, then $u_{\gamma }=v_{\gamma }+F$
is a
minimizer of the functional \eqref{4.10}. On the other hand, if $u_{\gamma }$ is a minimizer of the functional \eqref{4.10},
then $v_{\gamma }=u_{\gamma }-F\in
\mathcal{H}_{0}$ is a
minimizer of the functional \eqref{2.108}.

By the variational principle, any minimizer $v_{\gamma }$ of the functional \eqref{2.108} should satisfy the following condition
\begin{equation}\label{2.109}
\begin{array}{ll}\ds
\lan Pv_{\gamma },Ph\ran_{L_{\dbF}^2(0,T;L^2(G))} +\gamma \lan
v_{\gamma },h\ran_{L_{\dbF}^2(0,T;H^2(G))}  =\lan  Ph, \hat{f} -PF\ran_{L^{2}_{\mathbb{F}}(\Omega; H^{1}(0,T;L^{2}(G)))}
,\q\forall \, h\in \mathcal{H}_{0}
,
\end{array}
\end{equation}
%
%  Denote
% %
% \begin{equation}\label{2.110}
% \lan v,h\ran_{\gamma }=\lan 
% Pv,Ph\ran_{L_{\dbF}^2(0,T;L^2(G))} +\gamma \lan  v,h\ran_{L_{\dbF}^2(0,T;H^2(G))},\qq \forall
% v,h\in L_{\dbF}^2(0,T;H^2_{\G,0}(G)).
% \end{equation}
%
% Hence, $\lan v,h\ran_{\gamma }$ defines a new
% scalar product in the Hilbert
% space $L_{\dbF}^2(0,T;H^2_{\G,0}(G))$ and the corresponding norm $|v|_\g$ satisfies
% %
% \begin{equation}\label{2.111}
% \sqrt{\gamma }|v|_{L_{\dbF}^2(0,T;H^2(G))}\leq |
% v|_{\gamma }\leq C|
% v|_{L_{\dbF}^2(0,T;H^2(G))}.
% \end{equation}
% %
% Thus, the scalar product \eqref{2.110} generates
% the new norm $|v|_{\gamma}$,
% which is equivalent with the norm $|
% v|_{L_{\dbF}^2(0,T;H^2(G))
% }$. 
From \eqref{eqDou20231}, it can
be rewritten as
\begin{equation}\label{2.112}
\lan v_{\gamma },h\ran_{\mathcal{H}_{0}}=\lan
Ph, \hat{f}-PF\ran_{L^{2}_{\mathbb{F}}(\Omega; H^{1}(0,T;L^{2}(G)))},\q \forall \, h\in \mathcal{H}_{0}.  
\end{equation}
%
% It follows from (\ref{2.111}) that
From \eqref{eqfEstimate}, it holds that
\begin{equation}\label{2.113}
\big|\lan Ph, \hat{f}-PF\ran_{L^{2}_{\mathbb{F}}(\Omega; H^{1}(0,T;L^{2}(G)))}\big|
\leq C (| P F|_{L^{2}_{\mathbb{F}}(\Omega; H^{1}(0,T;L^{2}(G)))} + |f|_{L^{2}_{\mathbb{F}}(0,T;L^{2}(G))} ) |h|_{\mathcal{H}_{0} }.
\end{equation}
Hence, the right hand side of (\ref{2.112}) is a bounded linear functional
  on $\mathcal{H}_{0}$. By Riesz representation theorem, there
exists an element $w_{\gamma }\in \mathcal{H}_{0}$ such that 
$$
\lan Ph, \hat{f}-PF\ran_{L^{2}_{\mathbb{F}}(\Omega; H^{1}(0,T;L^{2}(G)))}
=\lan w_{\gamma },h\ran_{\mathcal{H}_{0}},\qq\forall \, h\in  \mathcal{H}_{0}.
$$ 
This and (\ref{2.112}) imply that
$$
\lan v_{\gamma },h\ran_{\mathcal{H}_{0}
}=\lan w_{\gamma },h\ran_{\mathcal{H}_{0}
},\qq\forall \, h\in  \mathcal{H}_{0}.
$$ 
Hence, the minimizer  $ 
v_{\gamma }=w_{\gamma }.$ 
% Also, by Riesz theorem
% and (\ref{2.113}),
% %
% $$
% |v_{\gamma }|
% _{\gamma }\leq C |
% F|_{L_{\dbF}^2(0,T;H^2(G))
% }.
% $$ 
% %
% Hence, the
% minimizer $v_{\gamma }$ is unique and the left
% inequality \eqref{2.111} implies \eqref{4.12}.
\endpf

\begin{remark}
In the proof of Theorem \ref{th4.1}, we utilized solely the variational principle and Riesz's theorem, without invoking the Carleman estimate. However, we make use of this estimate in Theorem \ref{th4.2}, where we establish the rate of convergence of minimizers $u_{\gamma}$ to the precise solution, provided that certain conditions are met.
\end{remark}

Assume that there exists an exact solution $y^{\ast }$ of the problem \eqref{system1} with the exact data
$$
f^{*} \in L^{2}_{\mathbb{F}}(0,T;L^{2}(G)), \quad 
y^{\ast }| _{\Gamma}=g_{1}^{\ast }\in
 H^1(\Gamma),\q\partial
_{\nu}y^{\ast }|_{\Gamma}=g_{2}^{\ast
}\in L_\dbF^{2}(0,T;L^2(\Gamma)).
$$ 
By
Theorem \ref{th1}, the exact solution $y^{\ast }$ is
unique. Because of the existence of $y^{\ast },
$ there also exists an exact function $F^{\ast
}\in \mathcal{H}$ such that
\begin{align}
  \label{eqFstar}
  | F^{*} |_{L^{2}_{\mathbb{F}}(0,T;H^{2}(G))} \leq C |y^{*}|_{L^{2}_{\mathbb{F}}(0,T;H^{2}(G))}  
    .
\end{align}
Here is
an example of such a function $F^{\ast }$. Let
the
function $\rho \in C^{2}\left(\overline{Q}\right) $ be such that $ 
\rho \left( x\right) =1$ in a small neighborhood
$$
\rho \left(t,x\right)=\begin{cases} \ds 1, & (t,x)\in 
N_{\sigma }\left((0,T)\times\Gamma\right) =\left\{
(t,x)\in  Q:\, \operatorname*{dist}\left((t,x),(0,T)\times\Gamma\right)
<\sigma \right\},\\
\ns\ds 0, & x\in  Q \setminus N_{2\sigma
}\left((0,T)\times\Gamma\right),
\end{cases}
$$ 
where $\sigma >0$
is a sufficiently small number. Then $F^{\ast }$
can be constructed as $F^{\ast
}\left(t, x\right) =\rho \left(t, x\right) y^{\ast }\left( t, x\right) $. Let $\delta >0$ be a sufficiently small number
characterizing the error in the data. We assume
that
\begin{equation}\label{2.114}
\begin{array}{ll}\ds
  |f^{*} - f|_{L^{2}_{\mathbb{F}}(0,T;L^{2}(G))} \leq \delta, \quad 
 |g_{1}^{\ast
}-g_{1}|_{ H^1(\Gamma) }\leq \delta,\q |g_{2}^{\ast
}-g_{2}|_{L^2_\dbF(0,T; L^2(\Gamma))}\leq \delta,
\end{array}
\end{equation}
and 
\begin{align}
  \label{eqFStarF}
   |PF^{*} - P F|_{L^{2}_{\mathbb{F}}(\Omega; H^{1}(0,T;L^{2}(G)))}  
   + |F^{\ast}-F|_{L_{\dbF}^2(0,T;H^2(G))}\leq \delta.  
\end{align}

\begin{theorem}\label{th4.2} (Convergence rate).
 Assume   \eqref{2.114} and \eqref{eqFStarF},   and let the regularization parameter
$\gamma   =\delta
^{2\alpha }$, where $\alpha  \in
\left(
0,1\right] $ is a constant.  Let $ G', \varepsilon, \delta_{0} $, and    $\beta $ be the
same as in Theorem \ref{th1}. Then there exists a
sufficiently small number $\delta _{1} \in \left( 0,1\right) $ and a constant  $ 
C >0$  such that if  $\delta_{1} \in
( 0,\delta _{0}^{1/\alpha })
$,  then  
\begin{equation}\label{2.1140}
|y_{\gamma }-y^{\ast }|
_{L_{\dbF}^2(\e,T-\e;H^1(G'))
}\leq C \big( 1+|y^{\ast
}|_{L_{\dbF}^2(0,T;H^2(G))
}\big) \delta^{\alpha \beta },\q \forall \, \delta
\in \left( 0,\delta _{1}\right),
\end{equation}
where $y_{\gamma \left( \delta \right)
}$ is the minimizer of the functional
\eqref{4.10}.
\end{theorem}

 \textbf{Proof}.  Let $v^{\ast }=y^{\ast }-F^{\ast }$.
Then $v^{\ast }\in \mathcal{H}_{0}$ and $Pv^{\ast }= \hat{f^{*}} -PF^{\ast
}$. Hence,
\begin{equation}\label{2.115}
\begin{array}{ll}\ds
\lan Pv^{\ast },Ph\ran_{L^{2}_{\mathbb{F}}(\Omega; H^{1}(0,T;L^{2}(G)))} +\gamma \lan
v^{\ast },h\ran_{L_{\dbF}^2(0,T;H^2(G))}\\
\ns\ds =\lan Ph, \hat{f^{*}} -PF^{\ast
}\ran_{L^{2}_{\mathbb{F}}(\Omega; H^{1}(0,T;L^{2}(G)))} +\gamma \lan v^{\ast },h\ran_{L_{\dbF}^2(0,T;H^2(G))}
,\q \forall \, h\in \mathcal{H}_{0}.  
\end{array}
\end{equation}
Subtract identity (\ref{2.109}) from identity (\ref{2.115}) and denote $ 
\widetilde{v}_{\gamma }=v^{\ast }-v_{\gamma }$, $ \widetilde{f} = \hat{f^{*}} -\hat{f} $ and 
$\widetilde{F}=F^{\ast }-F$. 
Then
\begin{equation*}
\begin{array}{ll}\ds
\lan P\widetilde{v}_{\gamma },Ph\ran_{L^{2}_{\mathbb{F}}(\Omega; H^{1}(0,T;L^{2}(G)))} +\gamma \lan \widetilde{v} 
_{\gamma },h\ran_{L_{\dbF}^2(0,T;H^2(G))}\\
\ns\ds =\lan  Ph, \widetilde{f} -P\widetilde{F}\ran_{L^{2}_{\mathbb{F}}(\Omega; H^{1}(0,T;L^{2}(G)))} +\gamma 
\lan v^{\ast },h\ran_{L_{\dbF}^2(0,T;H^2(G))},\q  \forall \, h\in \mathcal{H}_{0}. 
\end{array}
\end{equation*}
Setting here $h\=\widetilde{v}_{\gamma }$, we
obtain
\begin{equation}\label{2.116}
\begin{array}{ll}\ds
| P\widetilde{v}_{\gamma }|
_{L^{2}_{\mathbb{F}}(\Omega; H^{1}(0,T;L^{2}(G)))}^{2}+\gamma
| \widetilde{v}_{\gamma }|
_{L_{\dbF}^2(0,T;H^2(G))}^{2}\\
\ns\ds =\lan P\widetilde{v}_{\gamma }, 
\widetilde{f}-P\widetilde{F}\ran_{L^{2}_{\mathbb{F}}(\Omega; H^{1}(0,T;L^{2}(G)))} +\gamma \lan v^{\ast },\widetilde{v} 
_{\gamma }\ran_{L_{\dbF}^2(0,T;H^2(G))}.  
\end{array}
\end{equation}
Applying the Cauchy-Schwarz inequality to (\ref{2.116}), we obtain
\begin{equation}\label{2.1160}
\begin{array}{ll}\ds
|P\widetilde{v}_{\gamma }|
_{L^{2}_{\mathbb{F}}(\Omega; H^{1}(0,T;L^{2}(G)))}^{2}+\gamma
|\widetilde{v}_{\gamma }|
_{L_{\dbF}^2(0,T;H^2(G))}^{2}
\\
\ns\ds
\leq \frac{1}{2}|
P\widetilde{v}_{\gamma }|_{L^{2}_{\mathbb{F}}(\Omega; H^{1}(0,T;L^{2}(G)))}^{2}
+\frac{1}{2}\big|\widetilde{f}
- P\widetilde{F} \big|_{L^{2}_{\mathbb{F}}(\Omega; H^{1}(0,T;L^{2}(G)))}^{2}
+\frac{\gamma }{2} |v^{\ast }|_{L_{\dbF}^2(0,T;H^2(G))}^{2}
\\
\ns\ds \quad 
+\frac{ \gamma }{2}| \widetilde{v}_{\gamma
}|_{L_{\dbF}^2(0,T;H^2(G))
}^{2}.
\end{array}
\end{equation}
Hence, by \eqref{2.114} and \eqref{eqFStarF}
\begin{equation}\label{2.117}
|P\widetilde{v}_{\gamma }|
_{L^{2}_{\mathbb{F}}(\Omega; H^{1}(0,T;L^{2}(G)))}^{2}+\gamma
|\widetilde{v}_{\gamma }|
_{L_{\dbF}^2(0,T;H^2(G))}^{2}\leq
C \delta ^{2}+\gamma |v^{\ast }|_{L_{\dbF}^2(0,T;H^2(G))}^{2}. 
\end{equation}
Since $\gamma =\delta ^{2\alpha }$, where
$\alpha \in \left( 0,1\right] $,
then $\delta ^{2}\leq \gamma $. Hence, (\ref{2.117}) implies that
\begin{equation}\label{2.120}
| \widetilde{v}_{\gamma }|^{2}_{L_{\dbF}^2(0,T;H^2(G))}\leq
C \big( 1+| v^{\ast }|^{2}_{L_{\dbF}^2(0,T;H^2(G))}\big),
\end{equation}
\begin{equation}\label{2.121}
|P\widetilde{v}_{\gamma }|
_{L^{2}_{\mathbb{F}}(\Omega; H^{1}(0,T;L^{2}(G)))}^{2}\leq
C \big( 1+| v^{\ast }|
_{L_{\dbF}^2(0,T;H^2(G))}^{2}\big)
\delta ^{2\alpha }. 
\end{equation}
Let $w_{\gamma }=\widetilde{v}_{\gamma }\big(
1+|v^{\ast}| _{L_{\dbF}^2(0,T;H^2(G))}\big)^{-1}$. Then (\ref{2.120}), (\ref{2.121}) and Theorem \ref{th1} imply
that
$$| w_{\gamma }| _{L_{\dbF}^2(\e,T-\e;H^1(G'))}\leq
C\delta ^{\alpha \beta },\q\forall \, \delta \in
\left(
0,\delta _{1}\right).
$$ 
Therefore,
\begin{equation}\label{2.118}
| \widetilde{v}_{\gamma }|
_{L_{\dbF}^2(\e,T-\e;H^1(G'))
}\leq C \left( 1+| v^{\ast
}| _{L_{\dbF}^2(0,T;H^2(G))
}\right) \delta ^{\alpha \beta },\q\forall \,
\delta \in \left( 0,\delta _{1}\right).
\end{equation}
Next, since $\widetilde{v}_{\gamma }=\left(
  y^{\ast } - y_{\gamma }\right) +\left(F - F^{\ast
}\right) $,  by \eqref{2.114} and \eqref{eqFStarF},
$$
| F^{\ast }-F| _{L_{\dbF}^2(\e,T-\e;H^1(G'))}\leq \delta,
$$
then  
\begin{equation}\label{2.119}
\begin{array}{ll}\ds
|\widetilde{v}_{\gamma }|_{L_{\dbF}^2(\e,T-\e;H^1(G'))
}\3n&\ds\geq |y_{\gamma }-y^{\ast
}|_{L_{\dbF}^2(\e,T-\e;H^1(G'))}-| F^{\ast
}-F|_{L_{\dbF}^2(\e,T-\e;H^1(G'))}\\
\ns&\ds\geq |
y_{\gamma }-y^{\ast }|_{L_{\dbF}^2(\e,T-\e;H^1(G'))}-\delta.
\end{array}
\end{equation}
Since numbers $\beta ,\delta \in \left(
0,1\right) $ and since $\alpha \in
\left( 0,1\right] $, then $\delta ^{\alpha \beta }>\delta $. Thus, using  \eqref{eqFstar},
\eqref{2.118} and \eqref{2.119}, we obtain
\eqref{2.1140}. \endpf

\section{Numerical Approximations}

In this section, we aim to numerically solve the ill-posed Cauchy problem of the stochastic parabolic differential equation given by \eqref{system1}. For the sake of simplicity, we set $a_1=0$, $a_2=0$, $a_3=1$, $ f = 0 $  and $T=1$ in the system for all numerical tests to follow. Since an explicit expression for the exact solution is unavailable, we resort to numerically solving the initial-boundary value problem
\begin{equation}\label{spde}
\left\{
\begin{array}{ll}
\ds dy-\sum_{i,j=1}^n(a^{ij}y_i)_jdt=\lbrack\langle a_1,\nabla y\rangle+a_2y\rbrack dt+a_3ydW(t)\quad \mbox{in}\quad(0,1)\times G,\\
\ns\ds  y(0,x)=y_{0}(x) \quad\mbox{in}\quad G,\\
\ns\ds  y(t,x)=g_1(t,x) \quad \mbox{on}\quad (0,1)\times\partial G,
\end{array}\right.
\end{equation}
by employing the finite difference method with time discretized via the Euler-Maruyama method \cite{KP1992} to obtain the Cauchy data on $(0,1)\times\Gamma$. We then construct numerical approximations for the Cauchy problem \eqref{system1} via Tikhonov regularization, with the aid of kernel-based learning theory, for which we have established a convergence rate guaranteed in Section \ref{sec:4}.

%
%\subsection{Method of Fundamental Solutions}

%Now we solve the Cauchy problem of stochastic parabolic differential equation
%\begin{equation}\label{spde}
%du-\sum_{i,j=1}^n(a^{ij}u_i)_jdt=\lbrack\langle a_1,\nabla u\rangle+a_2u\rbrack dt+a_3udW(t)\quad \mbox{in}\quad(0,T)\times G
%\end{equation} 
%
%Then we should consider the following initial value problem of deterministic parabolic equation f
%\begin{equation}\label{pde}
%du-\sum_{i,j=1}^n(a^{ij}u_i)_jdt=\lbrack\langle a_1,\nabla u\rangle+a_2u\rbrack dt\quad \mbox{in}\quad(0,T)\times G
%\end{equation}

%We consider the problem for $G=(0,1)\subset\dbR^1$, $G$ is unit square and unit disk in $\dbR^2$, respectively, for different noisy levels. 
%Let $u(x,t)$ and $u_c(x,t)$ be the solution of problem which is obtained by solving \eqref{spde} numerically and its regularized approximation by the method given in section \ref{sec:4}, respectively. The error $E(x)$ is defined 
%by \begin{equation}
%E(x)=\frac{\|u_\lambda(x,t)-u(x,t)\|_{L^2(G)}}{\|u(x,t)\|_{L^2(G)}}
%\end{equation}
%for $x\in G$.
We verify the proposed method by using following three examples.

\begin{example}\label{1dex} Let $G=(0,1)$ and $\Gamma={1}$. Suppose that
\begin{enumerate}
\item[(a)] 
$
y_{0}(x)=x(1-x),\quad g_1(t,x)=0.
$
\item[(b)]
$
y_{0}(x)=\begin{cases}
4x,&x\in\lbrack 0,0.25)\\
-\frac{4}{3}x+\frac{4}{3},&x\in\lbrack 0.25,1\rbrack
\end{cases},\quad g_1(t,x)=0.
$
\end{enumerate}
\end{example}

%Firstly, we apply the finite difference method to solve the direct initial-boundary value problem numerically
%%
%\begin{equation}\label{forwardpro1}
%\begin{cases}
%y_t=u''+udW(t),&\quad (x,t)\in (0,1)\times (0,1),\\
%u(x,0)=f(x), & x\in(0,1),\\
%u(0,t)=g(t), & t\in(0,1),\\
%u(1,t)=g_1(t),& t\in(0,1).
%\end{cases}
%\end{equation}
%%
%to obtain the Cauchy data on the lateral boundary $x=1$.

The simplest time-discrete approximation is the stochastic version of the Euler approximation, also known as the Euler-Maruyama method \cite{KP1992}. 
We simply describe the process of solving the initial-boundary value problem \eqref{spde} in Example \ref{1dex} by the Euler-Maruyama method in the following. 
Let $y_{j,k}=y(x_j,t_k)$  with $x_j=jh$, $t_k=k\tau$, $j=1,\cdots,m+1$, $k=1,\cdots,n+1$, where $h=1/n$ and $\tau=1/m$ denote the spatial and temporal step sizes, respectively. It was shown that not all heuristic time-discrete approximations of stochastic differential equations converge to the solution process in a useful sense as the time step size $\tau$ tends to zero \cite{KP1992}. %Although this approximation is stable for computing the numerical solution of a stochastic parabolic equation, the results of multiple attempts may not differ significantly. %$$
%M=\lbrace (x_i,t_j):x_i=ih,t_j=j\tau,i\in \lbrack 1,m+1\rbrack,j\in \lbrack 1,n+1\rbrack\rbrace,
%$$
%where $h=1/n$ and $\tau =1/m$ are the spatial and temporal stepsizes, respectively.  
Since the numerical approximation of observations obtained by direct application of backward finite difference scheme is not adapted to the filtration $\dbF$, we can only solve the initial-boundary value problem \eqref{spde} by the explicit finite difference scheme, i.e.,
$$
\frac{y_{j,k+1}-y_{j,k}}{\tau}=\frac{y_{j-1,k}-2y_{j,k}+y_{j+1,k}}{h^2}+\frac{W(t_{k+1})-W(t_k)}{\tau}y_{j,k},\ 2\leq j\leq m,1\leq k\leq n.
$$
%Let $r=\frac{\tau}{h^2}$ and $W_k=\frac{W(t_{k+1})-W(t_k)}{\tau}$, then it can be rewritten as
%$$
%y_{j,k+1}=(1-2r+\tau W_k)y_{j,k}+r(y_{j-1,k}+y_{j+1,k}),2\leq j\leq m,1\leq k\leq n$$
where the initial and boundary value are given by
$$
y_{j,1}=y_{0}(x_j),\ 1\leq j\leq m+1,\quad
y_{1,k}=g_1(0,t_k),\ y_{m+1,k}=g_1(1,t_k),\ 2\leq k\leq n+1.
$$
To ensure the numerical scheme is stable, we choose $m=15$ and $n=450$ in the computation.
By solving above algebraic system, we obtain distribution of $y$ at meshed grids of the initial-boundary value problem \eqref{spde}. In the process of solving Cauchy problem \eqref{system1} numerically, we using $y_{m,k}$ and $y_{m+1,k}$ instead of the Cauchy data at $x=1$.
%$$
%\begin{cases}
%u(x_i,t_i)=u(x_{m+1},t_i)=g_{1,i},&i=1,2,\cdots,n,\\
%u(x_i,t_i)=u(x_m,t_i)=g_{1,i}-hg_{2,comp,i},&i=n+1,\cdots,2n.
%\end{cases}
%$$
%in the discrete case.

%Let $\sigma\in \lbrack -1,1\rbrack$ be the random variable representing the white noise, $g_1^{(m)}=\max_j|g_{1,j}|$ and $q^{(m)}=\max_j|q_{comp,j}|$, then the noisy data are 
%
%$$
%\begin{cases}
%\tilde{u}_{m+1i}=g_{1,i}+\delta g_1^{(m)}\sigma_i,&i=1,2,\cdots,n\\
%\tilde{u}_{mi}=g_{1,i}-h(g_{2,comp,i}+\delta g_2^{(m)}\sigma_i),&i=n+1,\cdots,2n.
%\end{cases}
%$$
%

In the following, we numerically solve the optimization problem which is given in \eqref{4.10} by kernel-based learning theory.

%One of the principal features which distinguish stochastic parabolic differential equations from deterministic differential equations, is that at least one of partial derivatives of the solution doesn't exist, which cause that the solution of any stochastic partial equation can not be expressed explicitly.  This yields one of the difficult obstacles to deal with the stochastic partial differential equations.  However, one may rewrite it as an integral equation, and then show that in this form there is a solution which is a continuous, though non-differentiable function. In multi-dimensional domains, the solution turns out to be a distribution, not a function. 
Suppose that $\varphi(x,t)$  is the fundamental solution of parabolic equation 
$$
%\left\{
%\begin{array}{ll}&\ds 
y_t-\sum_{i,j=1}^n(a^{ij}y_i)_j-\langle a_1,\nabla y\rangle+a_2y=0\quad \mbox{in}\quad(0,T)\times \mathbb{R}^d,
%\\
%\ns\ds& u(x,0)=y_{0}(x), \quad x\in\mathbb{R}^n,\end{array}\right.
$$
then the mild solution $y(x,t)$ of 
\begin{equation}\label{spde1}\left\{
\begin{array}{ll}  \ds dy-\sum_{i,j=1}^n(a^{ij}y_i)_jdt=\lbrack\langle a_1,\nabla y\rangle+a_2y\rbrack dt+a_3ydW(t)\quad \mbox{in}\quad(0,T)\times \mathbb{R}^n,\\
\ns\ds  y(x,0)=y_{0}(x), \quad x\in\mathbb{R}^n,\end{array}\right.
\end{equation} 
can be written as
\begin{equation}
y(x,t)=\int_0^t\int_{\mathbb{R}^{d}} \varphi(x,t,y,s)a_3(y)dydW(s)+\int_{\mathbb{R}^{d}}\varphi(x,t,y,0)y_{0}(y)dy.
\label{mildsolution}
\end{equation}
%where $G(x,t,y,s)$ is the fundamental solution of the above equation and $y_{0}(x)=u(x,0)$. 
%$G(x,t)=\frac{1}{(4\pi t)^\frac d 2}e^{-\frac{|x|^2}{4\pi t}}$ 
Since the fundamental solution $\varphi(x,t)$ is deterministic, so the mild solution $y(x,t)$ given by \eqref{mildsolution} is well-known.

%The mild solution \eqref{mildsolution} of \eqref{spde1} gives us the idea that regularization solution be given in \eqref{4.10}--\eqref{4.11} can be numerically solved by the optimization problem with aid of  the kernel-based learning theory. Recently, the theory has drawn attention in the studies of numerical computation for ordinary and partial differential \cite{Schaback2007}. In particular, the use of fundamental solutions as kernels has proven to be effective for solving inverse problems of deterministic evolution partial differential equations \cite{DH2012, HW2004}.  It's the first time to apply regularization method combining with the kernel-based theory to solving inverse problem for stochastic parabolic differential equations.
%In the following, we will illustrate the effective of the proposed method. Mention that one of the advantage of this method is that we can solve the problem in one step and without iteration, thus save computation cost to some extent. 

%In computation, we can rewrite the mild solution as
Assume that %$\tilde{u}(x,t)$
\begin{equation}\label{kbam}
%u(x,t)\thickapprox
%\tilde{u}(x,t)=\sum_{j=1}^N \lambda_j\phi_j(x,t)=\sum_{j=1}^N \lambda_j G(x-\eta_j,t-\tau_j),\quad j=1,\cdots,N,
\tilde{y}(x,t)=\sum_{l=1}^N \lambda_l\Phi_l(x,t),\quad l=1,\cdots,N. 
\end{equation}
Here $\Phi_l$ are the basis functions, $N$ is the number of source points, and $\lambda_l$ are unknown coefficients to be determined. 

From \eqref{mildsolution}, we can let $\Phi_l(x,t)=\varphi(x-\xi_l,t-\tau_l), l=1,\cdots,N $, where $\varphi$ is the fundamental solution of the deterministic parabolic equation, $(\xi_l,\tau_l)$ are source points.
The suitable choice of the source points can ensure that $\tilde{y}$ is analytic in the domain $(0,T)\times G$ and make the algorithm effective and accurate.  
However, the optimal rule for the location of source points is still open. In recent work, we occupy uniformly distributed points on $DT\times \lbrack -R,R+1\rbrack$ below $t=0$, where $DT$ and $R$ are post-prior determined, as the source points. This choice is given in \cite{DH2012} and related works, and performs better in comparison with other schemes. 
%\textcolor{blue}{In this part, we assume that the spatial domain $\Omega$ is the interval $\lbrack 0,1\rbrack$ and we fix $m=20$ and $n=40$ in the following computation. }

Since $\tilde{y}$ be given in \eqref{kbam} has already satisfied the equation in system \eqref{spde}, coefficients $\lambda_l,l=1,\cdots,N$ should be determined by Cauchy conditions. 
%From the process of solving the initial boundary value problem, we obtain the distribution of $y$ at meshed grids, and use $y'(1,t)\approx g_{2,comp}(t)=\frac{y(x_{m+1},t_i)-y(x_{m},t_i)}h$ as the Cauchy data. 
From the process of solving the initial boundary value problem we know that, $(x_{m+1},t_k)$ and $(x_{m},t_k), k=1,\cdots,n$ can be set as the collocation points, and the problem converts to solving unknowns 
$\lambda=[\lambda_1,\cdots,\lambda_N]^\mathsf{T}$ from the linear system 
%\begin{equation}
%\tilde{u}(1,t_i)=\sum_{j=1}^N \lambda_j\varphi(1-\xi_j,t-\tau_j)=
%\begin{cases}
%g(1,t_i)&i=1,2,\cdots,n,\\
%g(1,t_i)-g_2(x_i,t_i)\Delta x&i=n+1,\cdots,2n.
%\end{cases}
%\end{equation}
\begin{equation}\label{matrixequation}
A\lambda=b,
\end{equation}
where $$A=\[
\begin{array}{c}
\varphi(x_{m+1}-\xi_l,t-\tau_l)\\
\varphi(x_m-\xi_l,t-\tau_l)
\end{array}
\]_{2n\times N},$$
and
$$
b=
[y(x_{m+1},t_k),
y(x_m,t_k)]_{2n}^\mathsf{T},
$$ 
%Here $(x_{m+1},t_i), (x_{m},t_i)$ are the collocation points obtain by solving direct initial-boundary problem \eqref{spde}.
%Rewrite the above linear system into matrix form, i.e., 

%$$
%b=
%\begin{bmatrix}
%g_1(1,t_i)\\
%g_2(1,t_i)-g_2(x_i,t_i)\Delta x
%\end{bmatrix}_{n+m},
%$$ 
%and $\lambda$ is the unknown $n+m$ vector.
%$$
%\lambda=
%\begin{bmatrix}
%\lambda_1\\
%\vdots\\
%\lambda_{n+m}
%\end{bmatrix}.
%$$ 

By comparing the Tikhonov functional be given in \eqref{4.10} we know, $Pu$ is calculated by $A\lambda-b$. 
%Due to  severe ill-posedness of problem \eqref{system1} and ill-condition of the deducing matrix $A$, most of the standard numerical methods cannot achieve good accuracy in solving the matrix equation \eqref{matrixequation}. 
%Thus, regularization strategy is needed. In the computation we adapt the standard Tikhonov regularization to solve the matrix equation \eqref{matrixequation}, it means that, solving the least square problem 
%\begin{equation}
%\min_{\lambda}\lbrace \|A\lambda-b\|^2+\alpha^2\|\lambda\|^2\rbrace
%\end{equation}
%to obtain the Tikhonov-regularized solution $\lambda_\alpha$, where $\|\cdot\|$ denotes the usual Euclidean norm, $\alpha$ is the regularization parameter. 
Choosing the regularization parameter $\gamma$ by using the L-curve method, that is, a regularized solution near the ``corner" of the L-curve \cite{Hansen}
$$
L=\lbrace (\log{(\|\lambda_\gamma\|^2)},\log{(\|A\lambda_\gamma-b\|^2)}),\gamma>0\rbrace,
$$
and denoting the regularized solution of linear system \eqref{matrixequation} by $\lambda_\gamma^*$, leads to the approximated solution 
\begin{equation}
\tilde{y}_\gamma^*(x,t)=\sum_{l=1}^N\lambda_{\gamma,l}^*\Phi_l(x,t)
\end{equation}
of problem \eqref{system1}.

To illustrate the comparison of the exact solution and its approximation $\tilde{y}_\gamma^*$ and in order to avoid ``inverse crime", we choose $M=m+n$ and compute $y(x,t)$ by finite difference method again such that $y(x,t)$ and $\tilde{y}_\gamma^*$ be defined on the same grid. 

\begin{figure}[htbp]
  \centering
\subfloat[{$\delta=1\%$}]{ \includegraphics[width=0.3\textwidth]{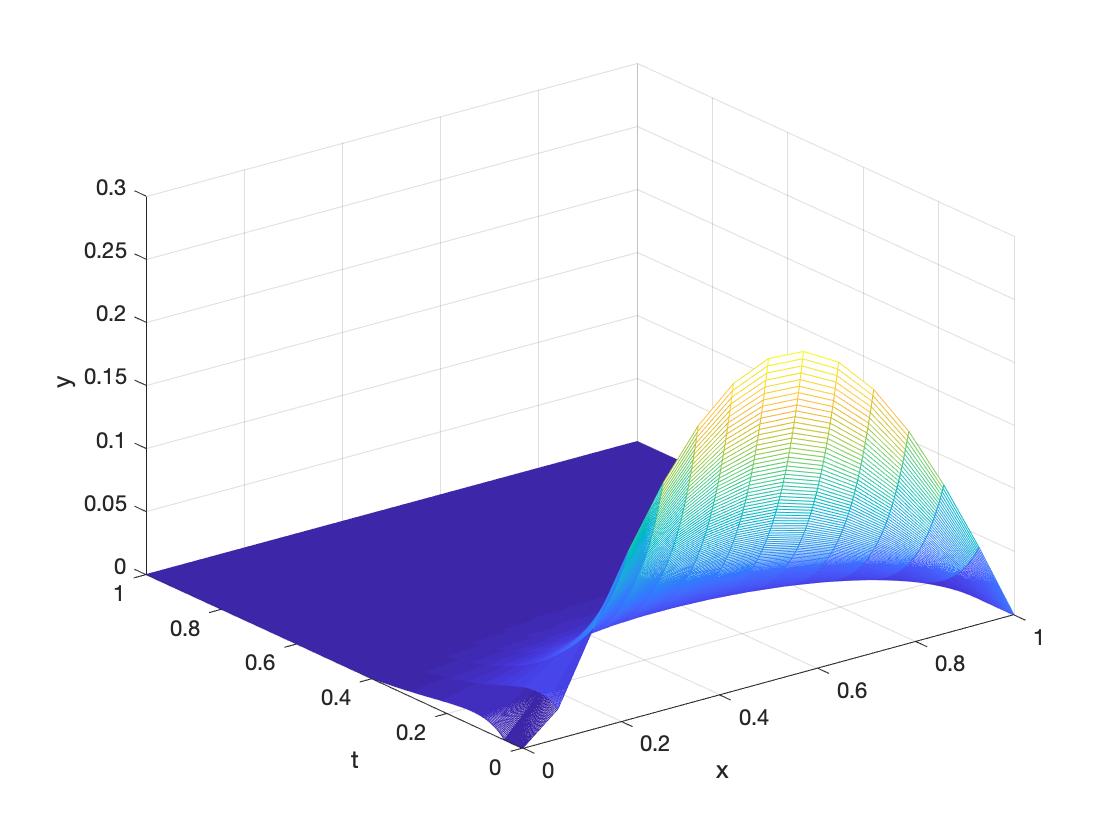}}
\subfloat[{$\delta=5\%$}]{\includegraphics[width=0.3\textwidth]{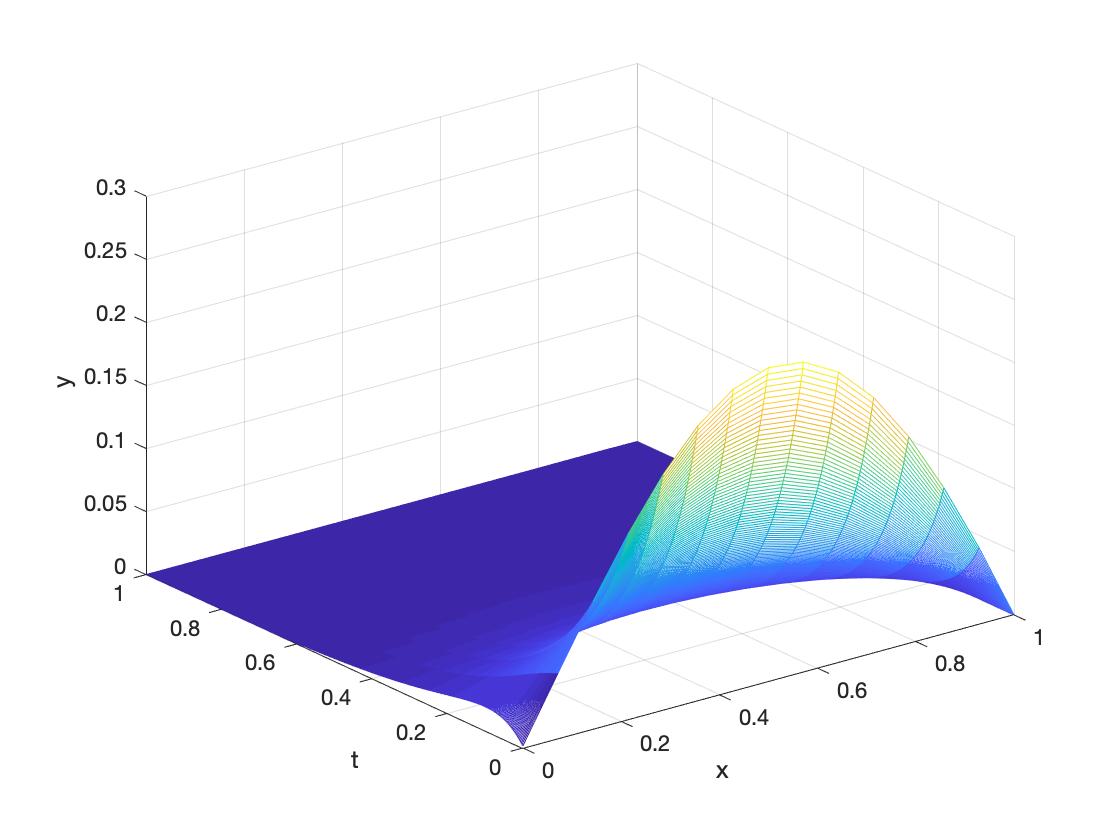}}
\subfloat[{$\delta=10\%$}]{ \includegraphics[width=0.3\textwidth]{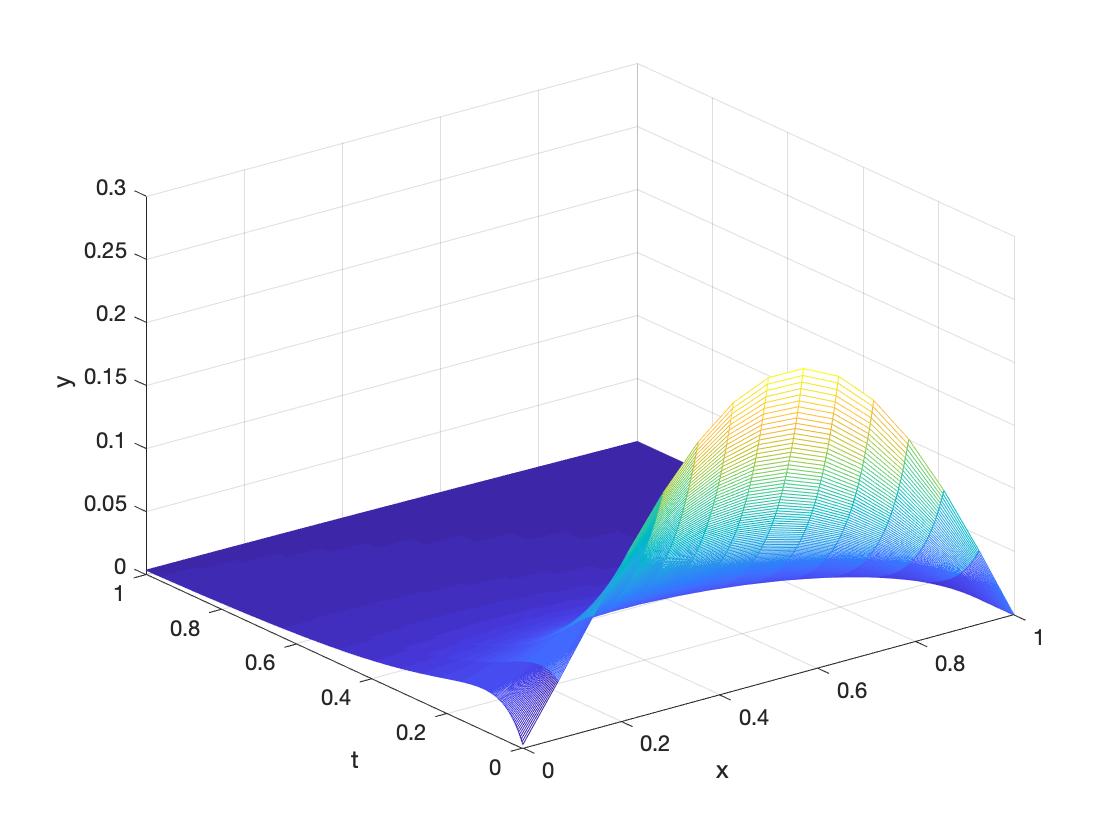}}
\caption{Approximated solution with different noisy levels of Example \ref{1dex}(a).}\label{ex1f1}
\end{figure}

Figure \ref{ex1f1} shows the numerical solution of the initial-boundary problem and the approximation solutions of the Cauchy problem of Example \ref{1dex}(a) with different noisy levels. Since the data is given at $x=1$, the numerical solution seems worse as $x$  tends to $0$, this is consistence with the convergence estimation in section \ref{sec:4} because no information is given at partial boundary $\{x=0\}$. Furthermore, the convergence rate always holds in the temporal interval $[\varepsilon, T-\epsilon], \epsilon>0$ in section \ref{sec:4}. However, from this figure, we find that the proposed method also works well at $t=0$. 
 %One can observe that the proposed numerical scheme performs well when $\delta$ is small, even when $\delta=10\%$, the approximation solution is still acceptable. 
%\begin{figure}[htbp]
%  \begin{center}
%  \tiny{(a)}\includegraphics[width=0.4\textwidth,natwidth=1200,natheight=1000]{delta=0.01exact1dcount=1.jpg}
% \tiny{(b)}\includegraphics[width=0.4\textwidth,natwidth=1200,natheight=1000]{delta=0.011dcount=1.jpg}\\
%\tiny{(c)}\includegraphics[width=0.4\textwidth,natwidth=1200,natheight=1000]{delta=0.051dcount=1.jpg}
%\tiny{(d)} \includegraphics[width=0.4\textwidth,natwidth=1200,natheight=1000]{delta=0.11dcount=1.jpg}
%\end{center}
%\caption{The exact solution and its approximations with $\delta=1\%,5\%,10\%$. (a) Exact solution. (b) $\delta=1\%$. (c) $\delta=5\%$. (d) $\delta=10\%$. }
%\end{figure}\label{ex1f1}

Denote the relative error by
\begin{equation}\label{error1}
E(x)=\frac{\|\tilde{y}_\gamma^*(x,t)-y(x,t)\|_{L^2(G)}}{\|y(x,t)\|_{L^2(G)}}.
\end{equation} 

\begin{figure}[htbp]
\centering
\subfloat[{Relative errors $E(x)$}]{\includegraphics[width=0.34\textwidth]{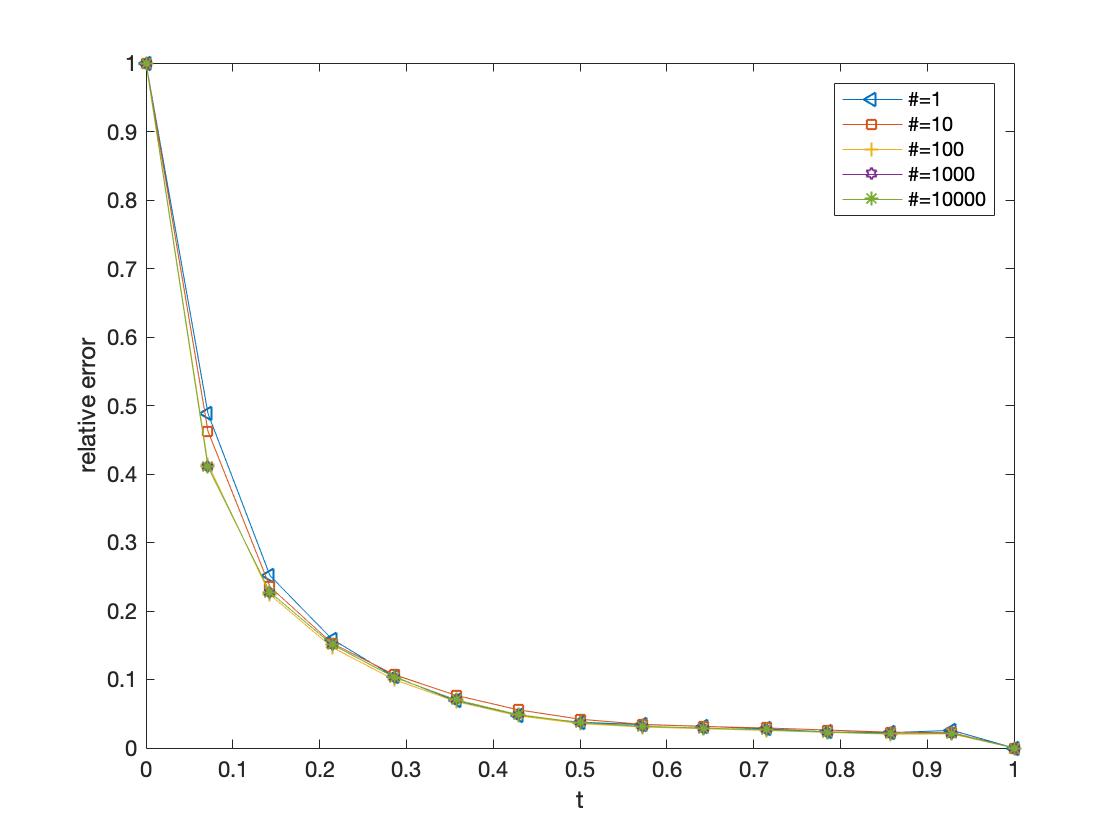}}
\subfloat[{Numerical approximations of $y(0,x)$}]{\includegraphics[width=0.34\textwidth]{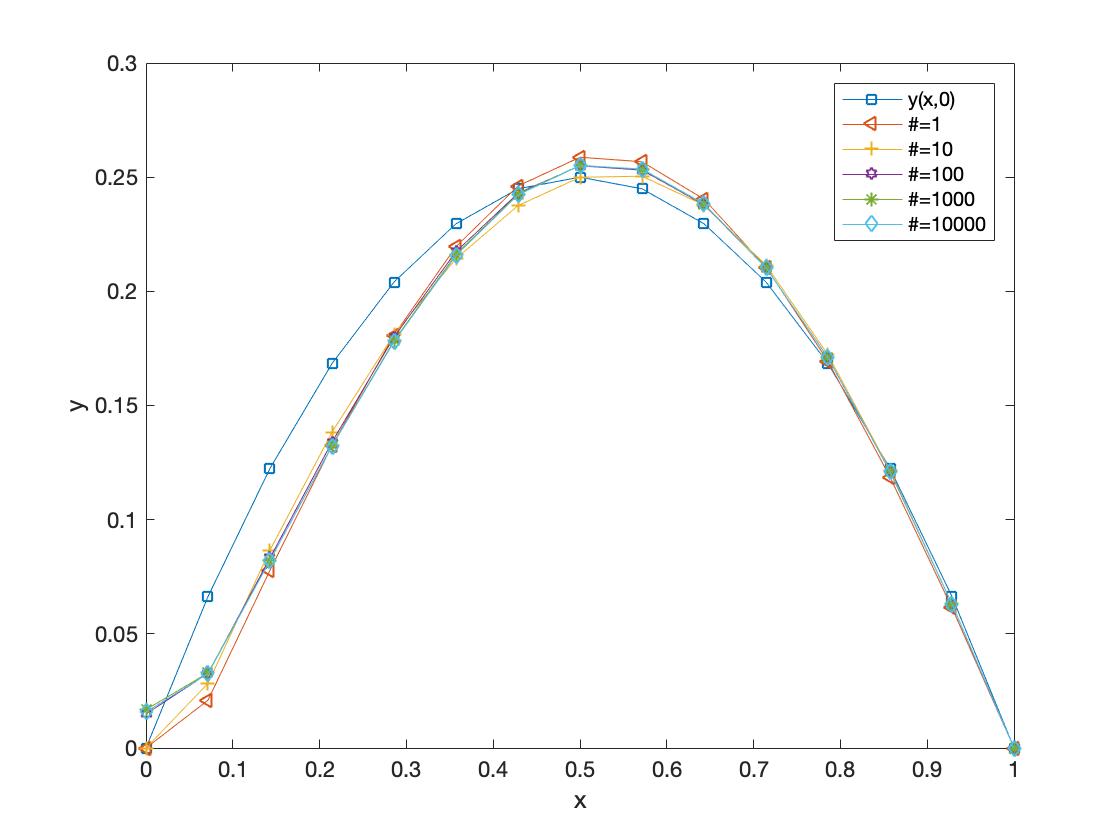}}
\caption{Numerical results by expectation of solutions of different numbers of sample paths.}\label{ex1f3}
\end{figure}

As in the study of stochastic equations, large numbers of sample paths should be checked in the numerical experiments for simulating the expectation of the solution. Thus, we do the tests with different sample paths. It is interesting from Figure \ref{ex1f3} that when the number of sample paths is great than 10, the results seems no better off. Thus in the following, we only consider do the experiments with number of sample paths $\#=10.$

\begin{figure}[htbp]
\centering
\subfloat[{Approximations of $y(t,0)$}]{\includegraphics[width=0.34\textwidth]{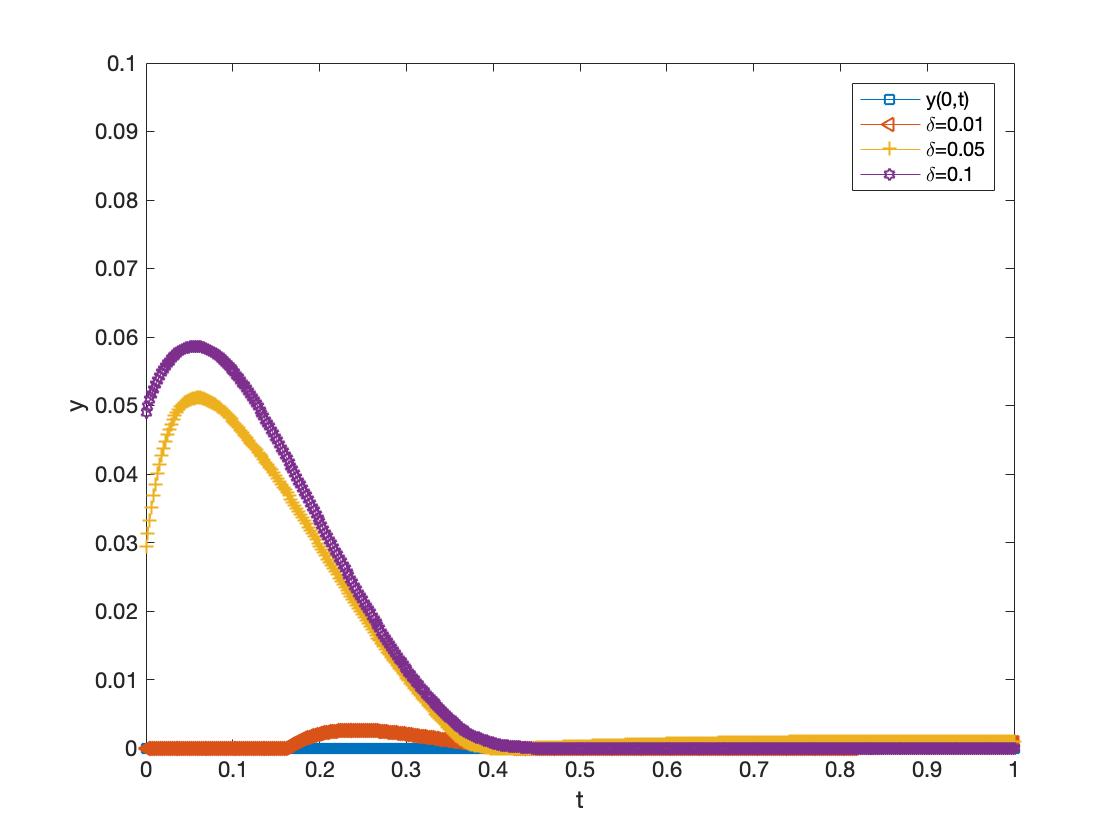}}
\subfloat[{Approximations of $y(0,x)$}]{\includegraphics[width=0.34\textwidth]{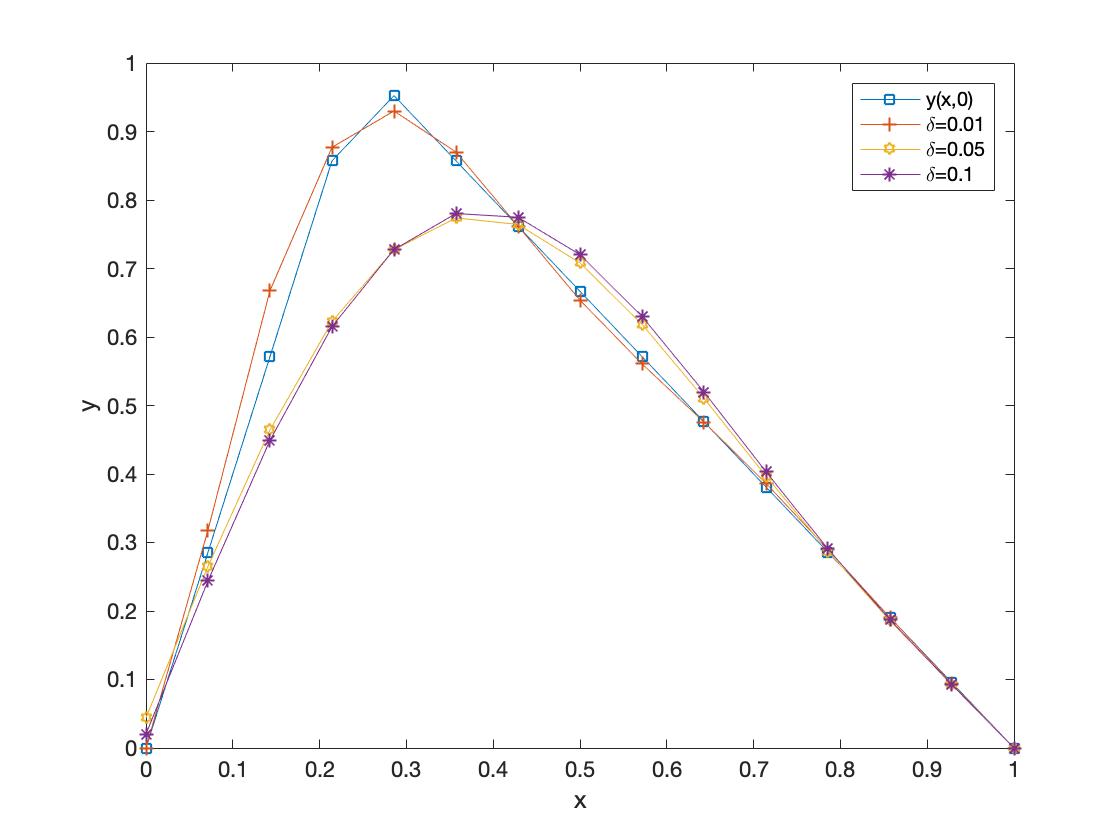}}\\
\subfloat[{Relative errors $E(x)$}]{\includegraphics[width=0.34\textwidth]{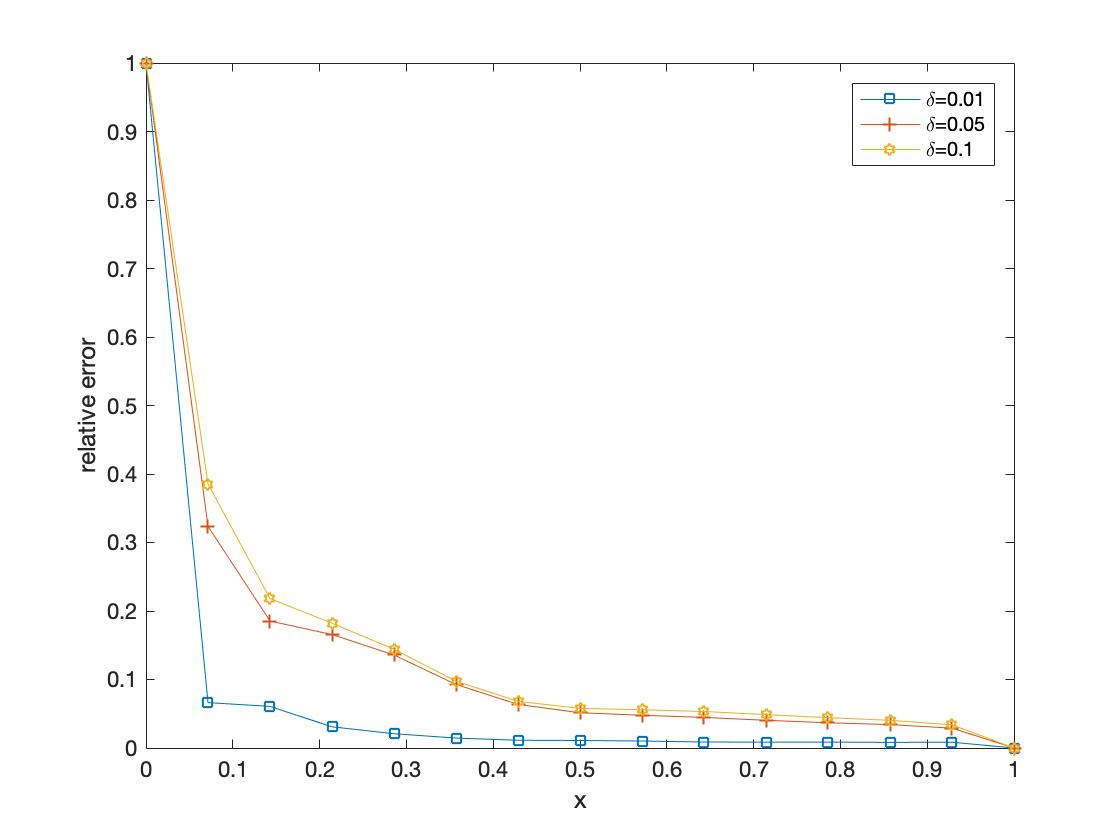}}
\subfloat[{Relative errors $E(t)$}]{\includegraphics[width=0.34\textwidth]{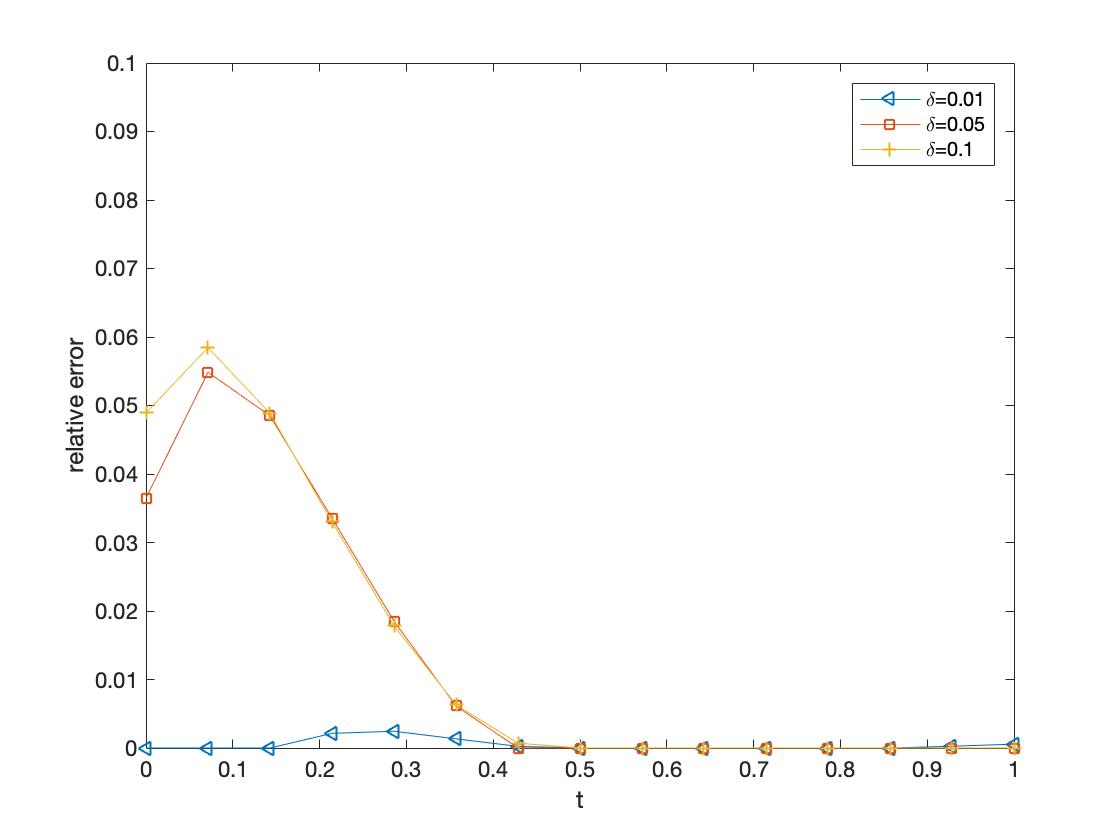}}
\caption{Numerical results of Example \ref{1dex}(b) with different noisy level data $\delta$.}
\label{ex2f1}
\end{figure}

%Firstly, we consider $\Gamma=\lbrace x=1\rbrace$ for different noise level of one single experiment.
%\begin{figure}[htbp]
%  \centering
%\subfloat[{numerical solution of initial-boundary value problem}]{\includegraphics[width=0.4\textwidth,natwidth=1200,natheight=1000]{delta=0.011d2exactu.jpg}}
%\subfloat[{approximation with $\delta=1\%$ and count=1}]{ \includegraphics[width=0.4\textwidth,natwidth=1200,natheight=1000]{delta=0.011d2count=1.jpg}}\\
%\subfloat[{approximation with $\delta=5\%$ and count=1}]{\includegraphics[width=0.4\textwidth,natwidth=1200,natheight=1000]{delta=0.051d2count=1.jpg}}
%\subfloat[{approximation with $\delta=10\%$ and count=1}]{ \includegraphics[width=0.4\textwidth,natwidth=1200,natheight=1000]{delta=0.11d2cpunt=1.jpg}}
%\caption{}
%\end{figure}
%And the relative error is given by fugure().
Now we consider the problem in Example \ref{1dex}(b). In this case, $y_{0}$ is a piecewise smooth function. 
Figure \ref{ex2f1} shows the approximations of boundary value $y(t,0)$ in (a) and the approximations of initial value $y(0,x)$ with different noise levels in (b). Moreover, to illustrate the effectiveness of the proposed method, the change of relative error along lines $x$ and $t$ are also given in Figure \ref{ex2f1}(c) and (d), respectively.  
%However, we can show in Figure \ref{ex2f1}(b) that the numerical approximations also work well at $t=0$.
%changes with the change of $x$, we set $N=2n$ to make $u(x,t)$ and $u_\lambda(x,t)$ be defined on the same grid and norms used below. For each $x$ from this grid we define the 

We mention that the Cauchy problem of the deterministic parabolic equation  has been considered in \cite{KKLY2017}, that the numerical algorithm based on Carleman weight function proposed can stably invert the solution near the known boundary conditions, more precisely, solutions for $t\in[0.5,1]$ for 1-D example. With comparison of the results of Example \ref{1dex}, our method works well for recovering the solution for $t\in [0.1,1]$ for the stochastic parabolic equation. We believe this algorithm also works well for deterministic parabolic equation, and can be an improvement of previous works. Furthermore, the proposed method also works well for 2-D examples in both rectangular and discal domains.

%\section{Numerical solution for the two-dimensional situation}
We explore the effectiveness of the proposed numerical methods in following two examples in bounded domains in $\mathbb{R}^2$, as well as the influence of parameter in the kernel-based learning method. We always let $\delta=3\%$ unless particularly stated.

\begin{figure}[htbp]
%\centering
\vspace{1cm}
\includegraphics[width=\textwidth]{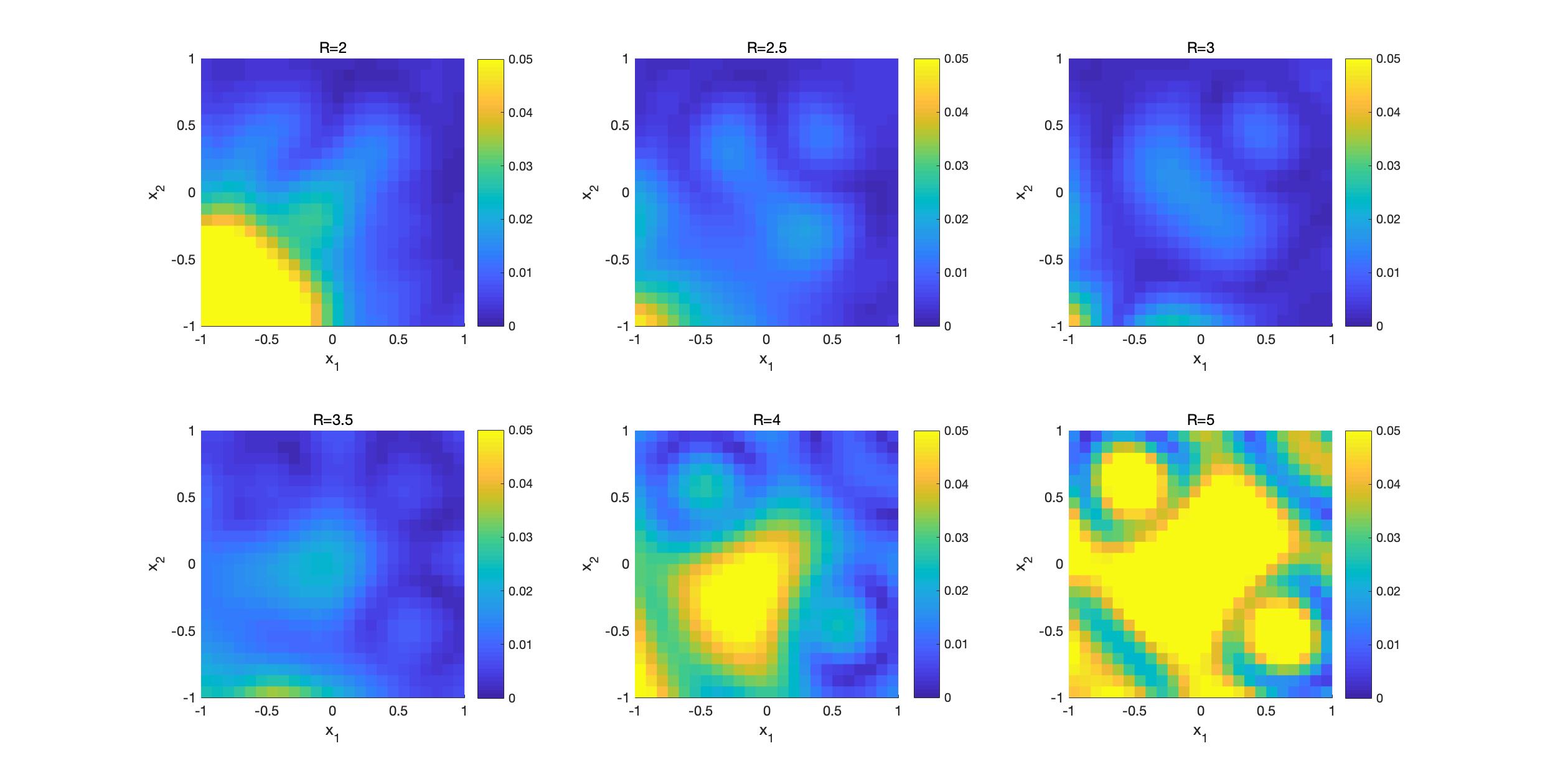}
\caption{The relative error $E(x_1,x_2)$ for different $R$.}\label{ex3R}
\end{figure}

\begin{figure}[htbp]
%\centering
\includegraphics[width=\textwidth]{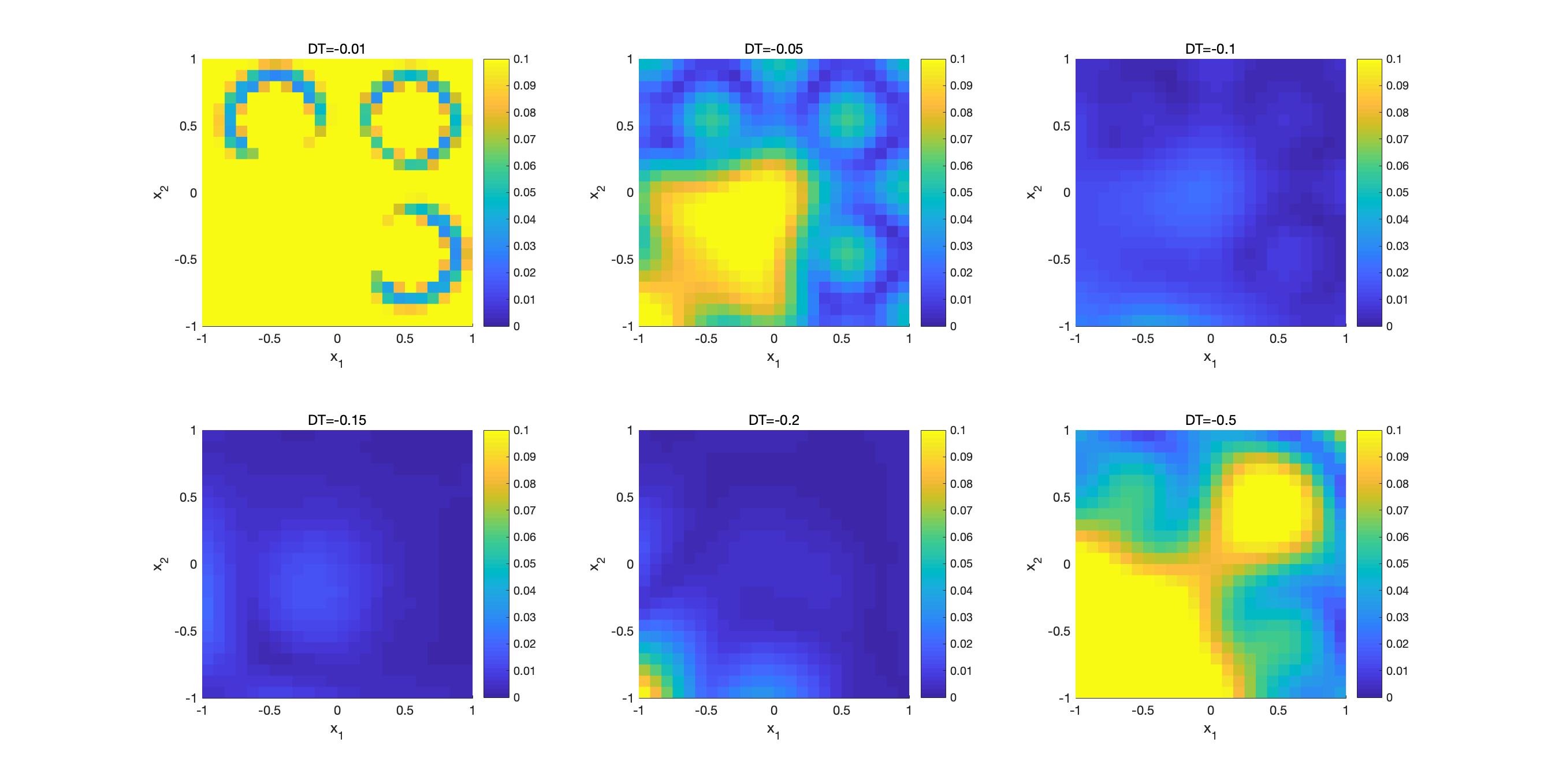}
\caption{The relative error $E(x_1,x_2)$ for different $DT$.}\label{ex3dt}
\end{figure}

\begin{figure}[htbp]
  \centering
%\subfloat[{numerical solution $y(x_1,x_2,0)$ of initial-boundary value problem}]{\includegraphics[width=0.4\textwidth,natwidth=1200,natheight=800]{sinexactut=0.jpg}}
\subfloat[{approximation with $\delta=1\%$ }]{ \includegraphics[width=0.3\textwidth]{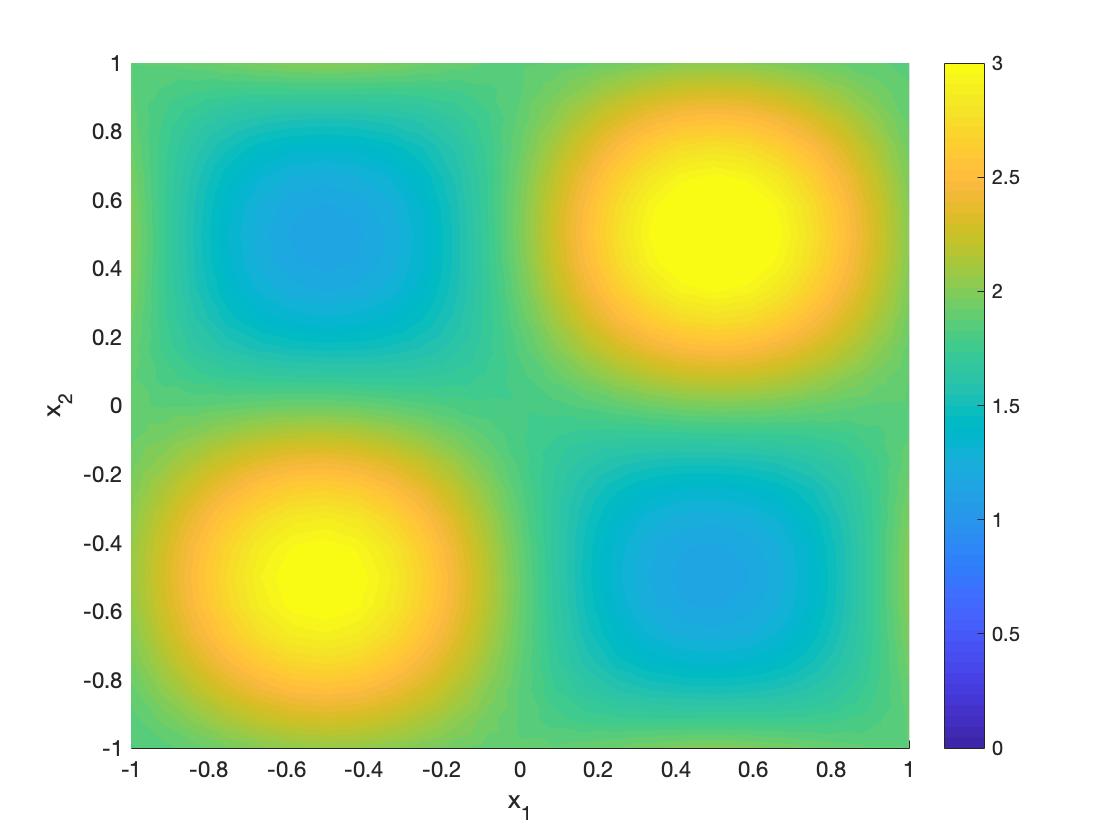}}
\subfloat[{approximation with $\delta=5\%$ }]{\includegraphics[width=0.3\textwidth]{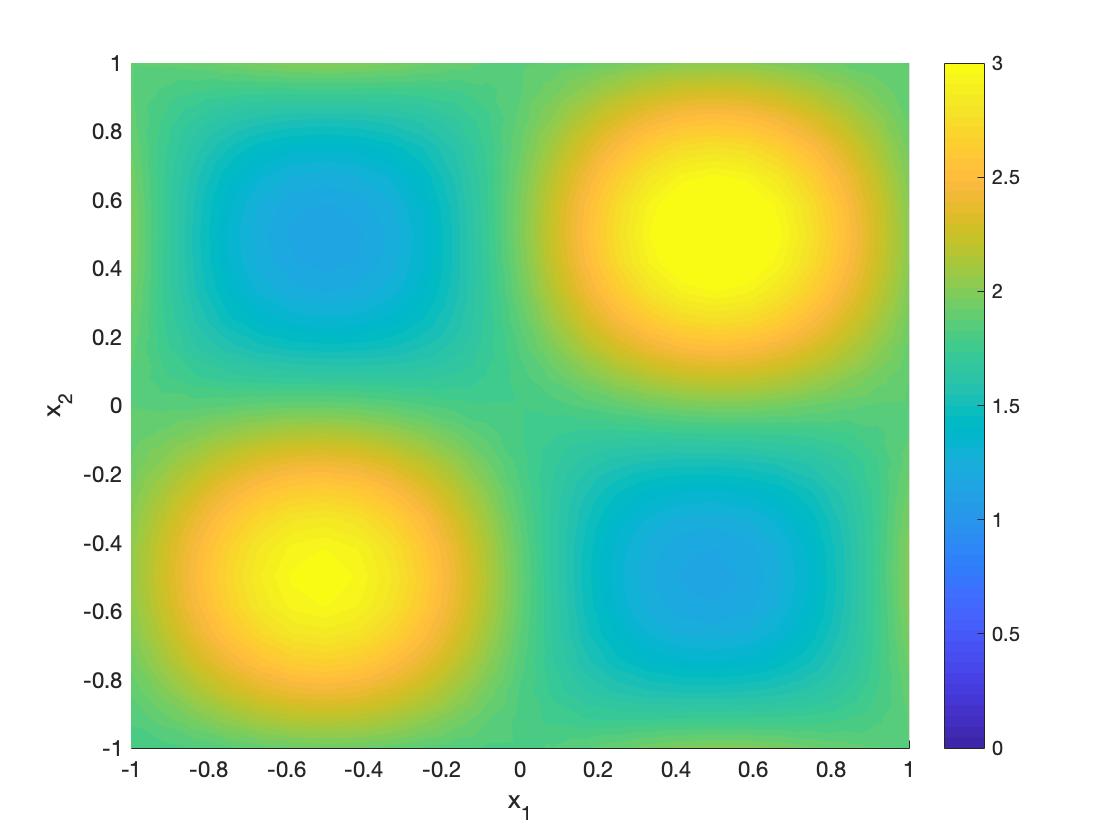}}
\subfloat[{approximation with $\delta=10\%$ }]{ \includegraphics[width=0.3\textwidth]{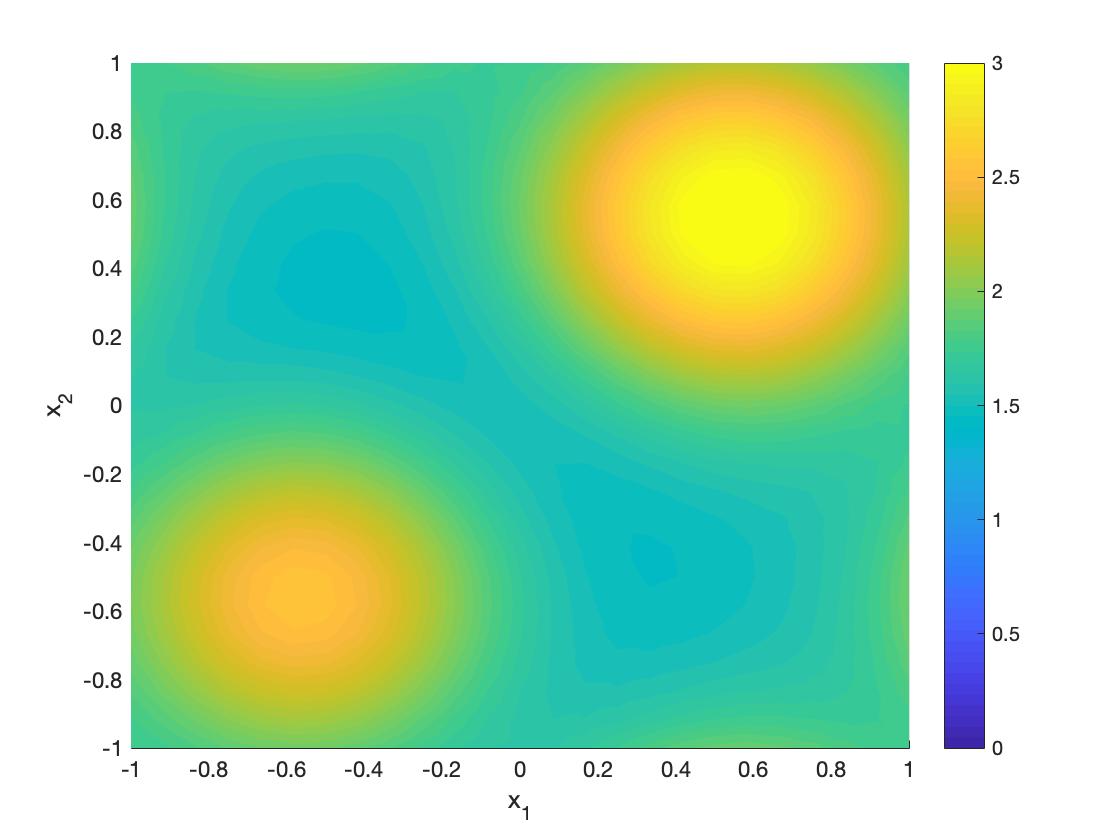}}\\
\subfloat[{relative error with $\delta=1\%$ }]{ \includegraphics[width=0.3\textwidth]{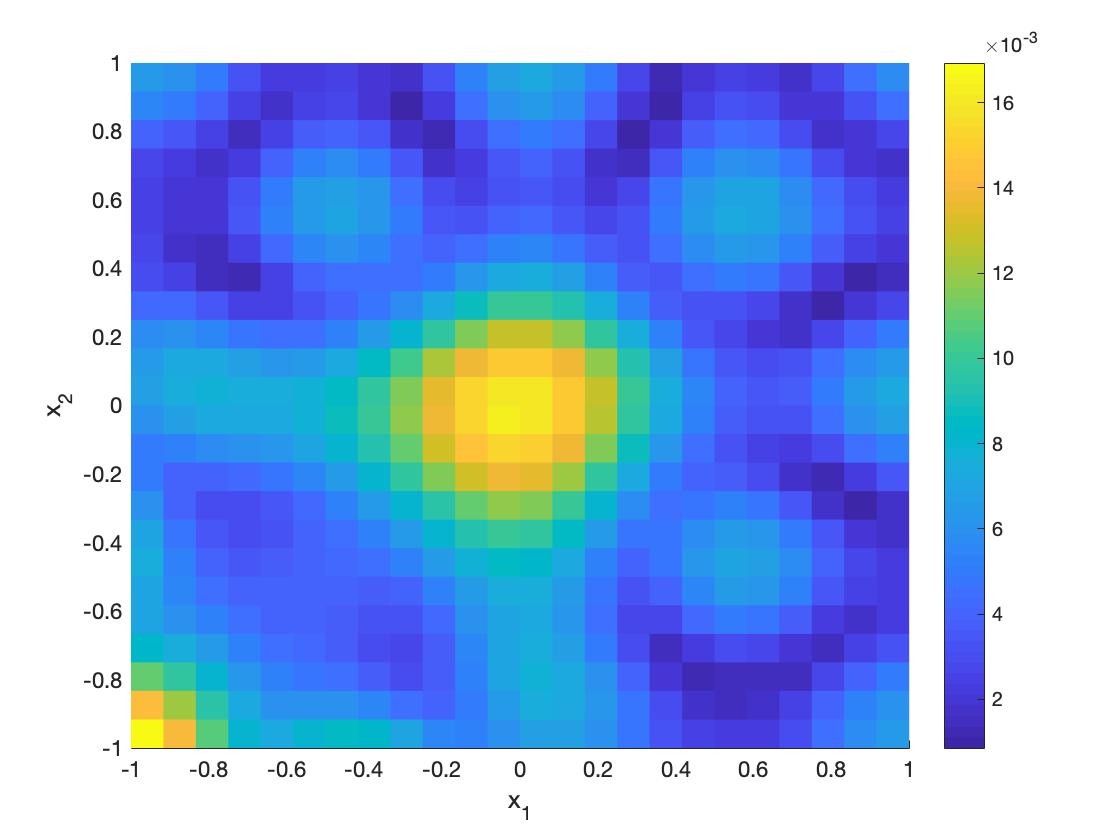}}
\subfloat[{relative error with $\delta=5\%$ }]{\includegraphics[width=0.3\textwidth]{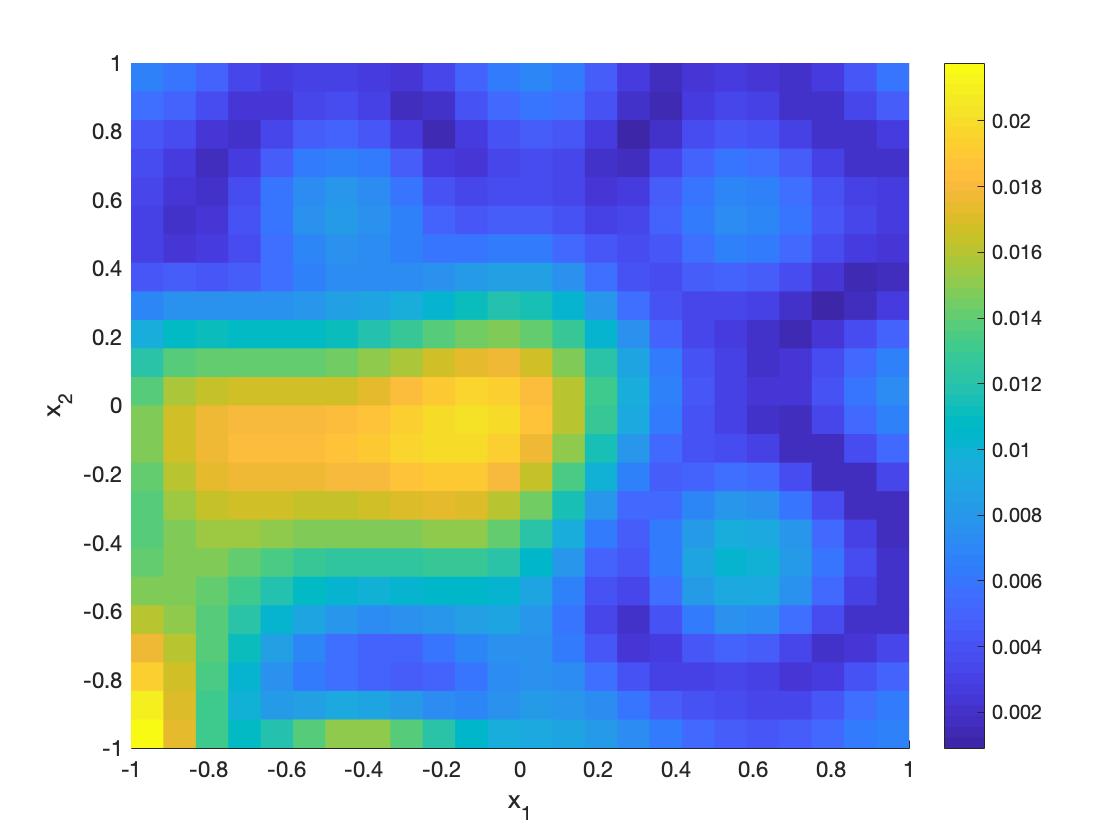}}
\subfloat[{relative error with $\delta=10\%$ }]{ \includegraphics[width=0.3\textwidth]{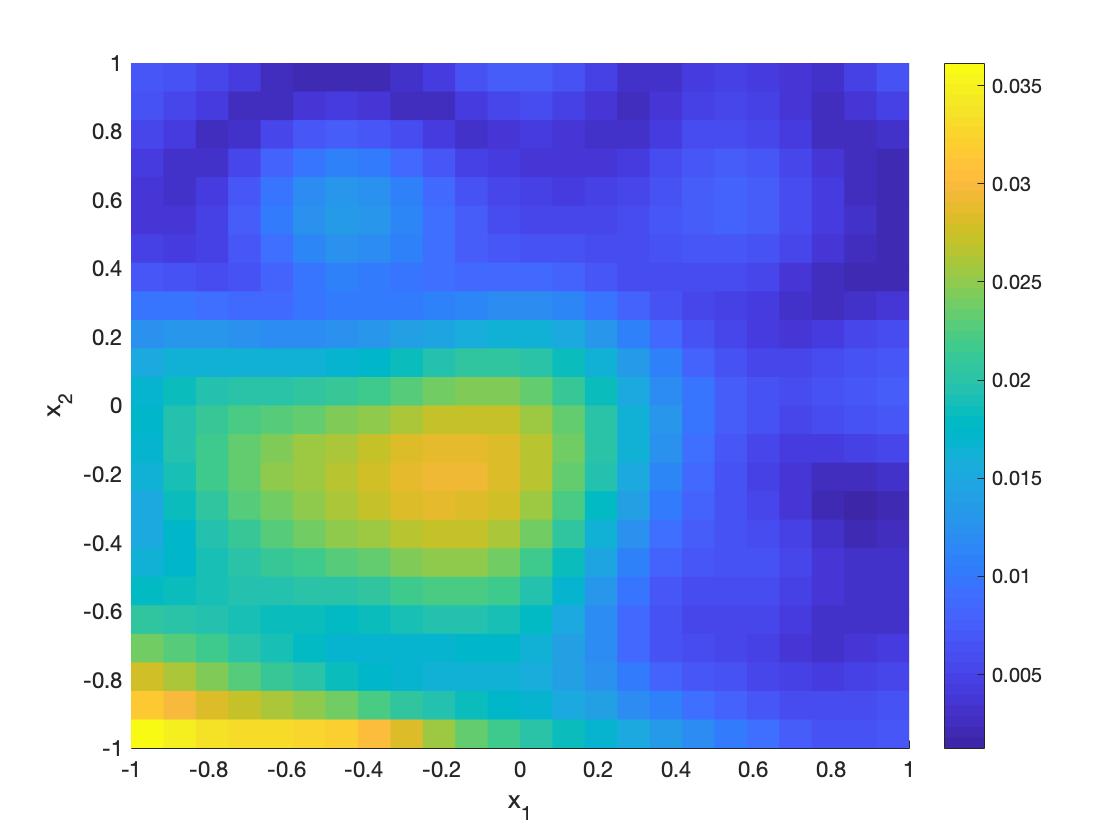}}
\caption{The approximations for $y(0,x_1,x_2)$ and relative errors $E(x_1,x_2)$ for different noise levels for Example \ref{ex4}(a).}\label{ex3err}
\end{figure}

\begin{example}\label{ex4}
Suppose that $G=\lbrack -1,1\rbrack \times\lbrack -1,1\rbrack$, and $\Gamma_1=\{x_1=1,x_2\in[-1,1]\}, \Gamma_2=\{x_1\in[-1,1]\},\Gamma_3=\partial G\backslash\{x_1=-1,x_2\in[-1,1]\}$. Let
\begin{enumerate}
\item[(a)]
$
y_{0}(x_1,x_2)=\sin{(\pi x_1)}\sin{(\pi x_2)}+2,\quad g_1(t,x_1,x_2)=2.
$
\item[(b)]
$
y_{0}(x_1,x_2)=
\begin{cases}
3,\quad (x_1-0.5)^2+(x_2-0.5)^2\leq 0.15^2\\
3,\quad (x_1-0.5)^2+(x_2+0.5)^2\leq 0.15^2\\
3,\quad (x_1+0.5)^2+(x_2-0.5)^2\leq 0.15^2\\
3,\quad (x_1+0.5)^2+(x_2+0.5)^2\leq 0.15^2\\
1,\quad otherwise.
\end{cases},\quad g_1(t,x_1,x_2)=1.
$
\end{enumerate}
\end{example}
We first consider the optimal choices of parameters $R$ and $DT$ in case that the Cauchy data are given on $\Gamma_2$ in the kernel-based approximation processes. We set  $DT = -0.1$ and observe relative errors $E$ with the change of $R$. It can be seen from figure \ref{ex3R} that $R \in \lbrack 2.5,3.5\rbrack$, the methods performs well. Numerical results in Figure \ref{ex3dt} illustrate that any $DT \in \lbrack-0.2,-0.1\rbrack$ is a reasonable choice. Thus we fix $R = 3.5,DT=-0.1$ in the following computing. 
The numerical approximations and relative errors for different noise level are shown in Figure \ref{ex3err}.

%\begin{figure}[htbp]
%  \centering{
%%\subfloat[{numerical solution $y(x_1,x_2,0)$ of initial-boundary value problem}]{\includegraphics[width=0.4\textwidth,natwidth=1200,natheight=800]{sinexactut=0.jpg}}
%\subfloat[{approximation with $\delta=1\%$ and count=1}]{ \includegraphics[width=0.3\textwidth,natwidth=1200,natheight=800]{sindelta=0.01t=0.jpg}}
%\subfloat[{approximation with $\delta=5\%$ and count=1}]{\includegraphics[width=0.3\textwidth,natwidth=1200,natheight=1000]{sindelta=0.05t=0.jpg}}
%\subfloat[{approximation with $\delta=10\%$ and count=1}]{ \includegraphics[width=0.3\textwidth,natwidth=1200,natheight=1000]{sindelta=0.1t=0.jpg}}}
%\end{figure}
%And the relative error is given by fugure().
%\begin{figure}[h]
%\centering{
%\includegraphics[width=0.8\textwidth,natwidth=1800,natheight=600]{sindeltachange.jpg}}
%\caption{The relative error $E(x_1,x_2)$ for different noise level of single experiment.}
%\end{figure}
%
%
%
%Secondly, we try to obtain the result when $\delta=0.03$ by averaging the number of experiments.

%Finally, we discuss the boundary condition.
%\begin{figure}[h]
%\centering
%\includegraphics[width=0.9\textwidth,natwidth=1800,natheight=600]{sindelta=0.03.jpg}
%\end{figure}

\begin{figure}[htbp]
  \centering
\subfloat[{approximation with data on $\Gamma_1$}]{\includegraphics[width=0.3\textwidth]{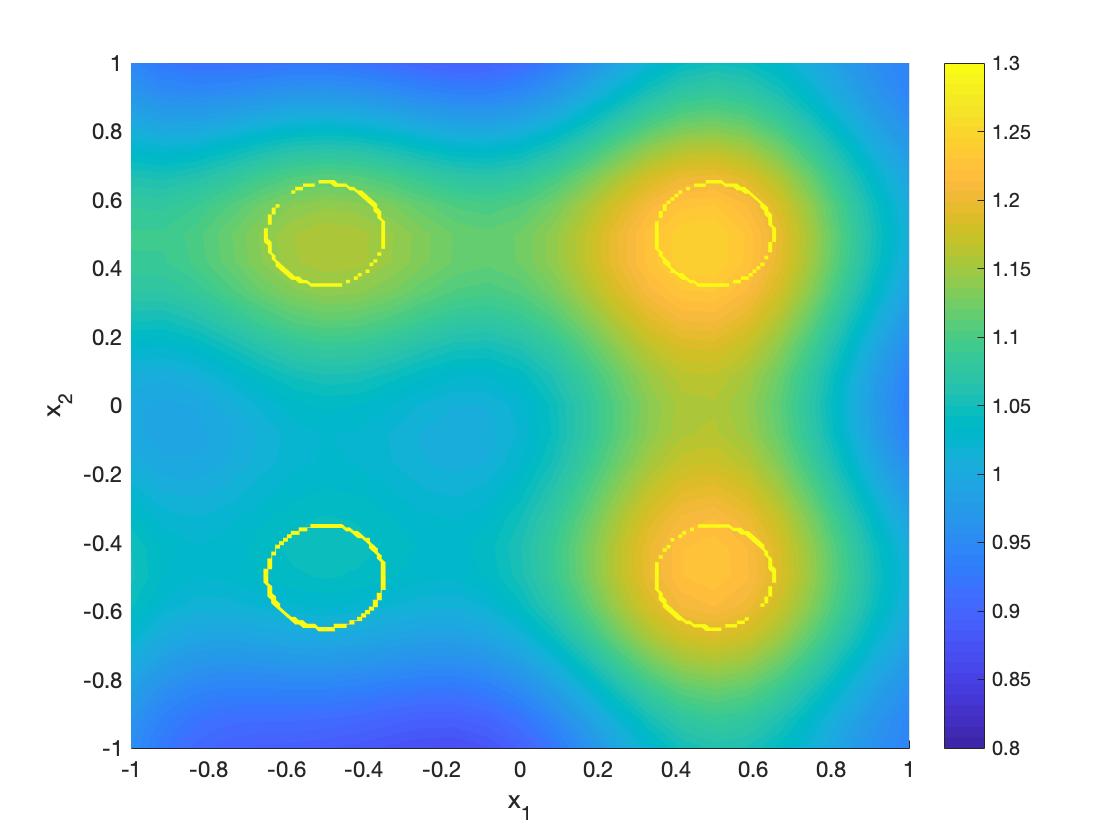}}
\subfloat[{approximation with data on $\Gamma_2$}]{ \includegraphics[width=0.3\textwidth]{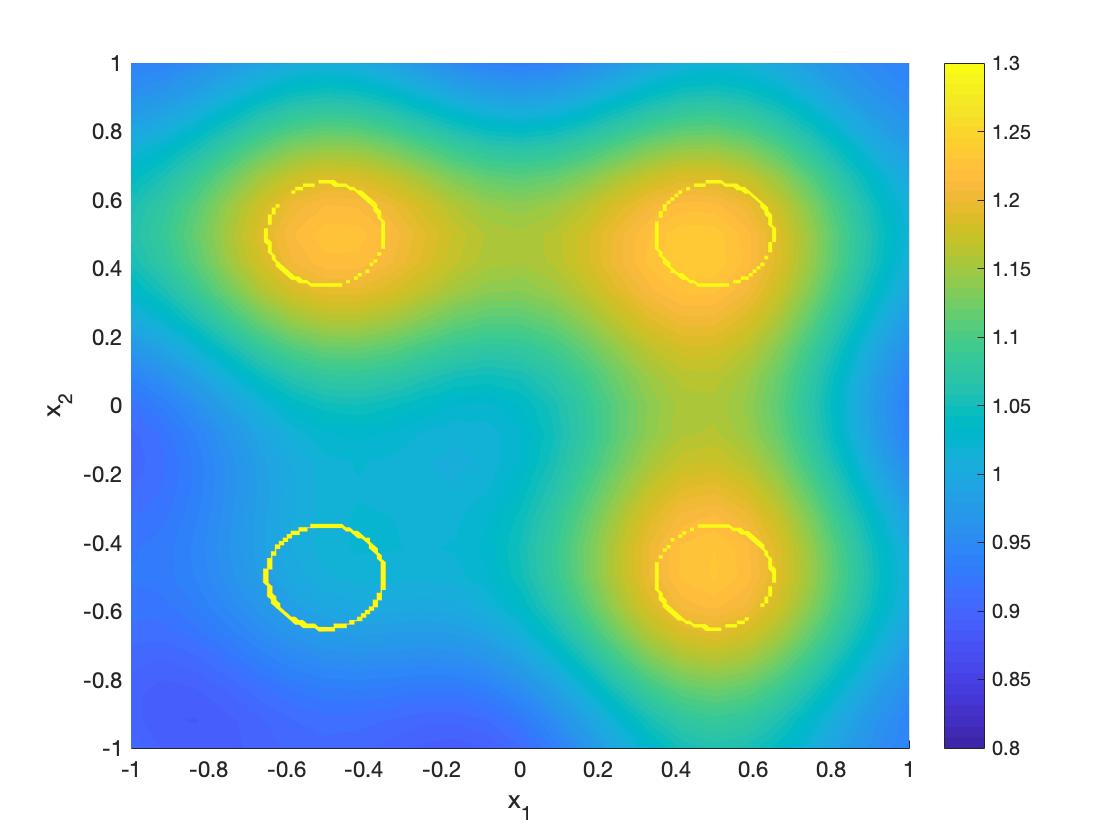}}
\subfloat[{approximation with data on $\Gamma_3$}]{ \includegraphics[width=0.3\textwidth]{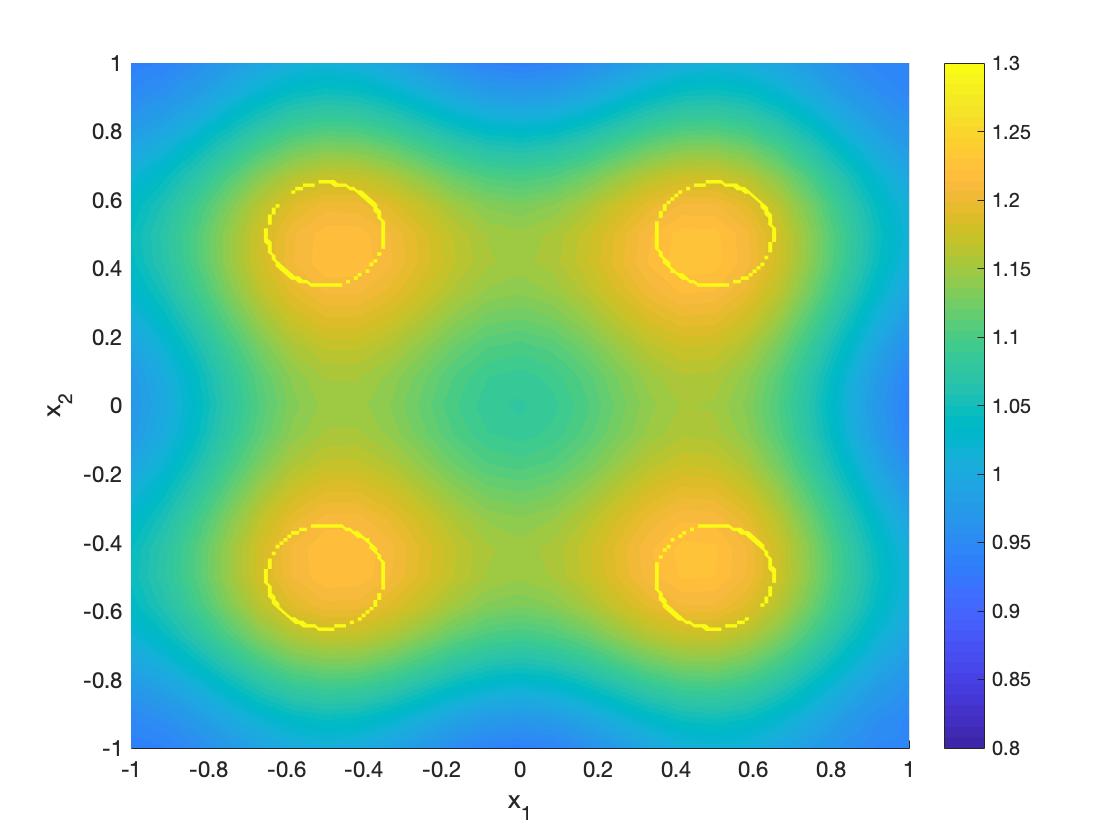}}\\
\subfloat[{relative error with data on $\Gamma_1$}]{\includegraphics[width=0.3\textwidth]{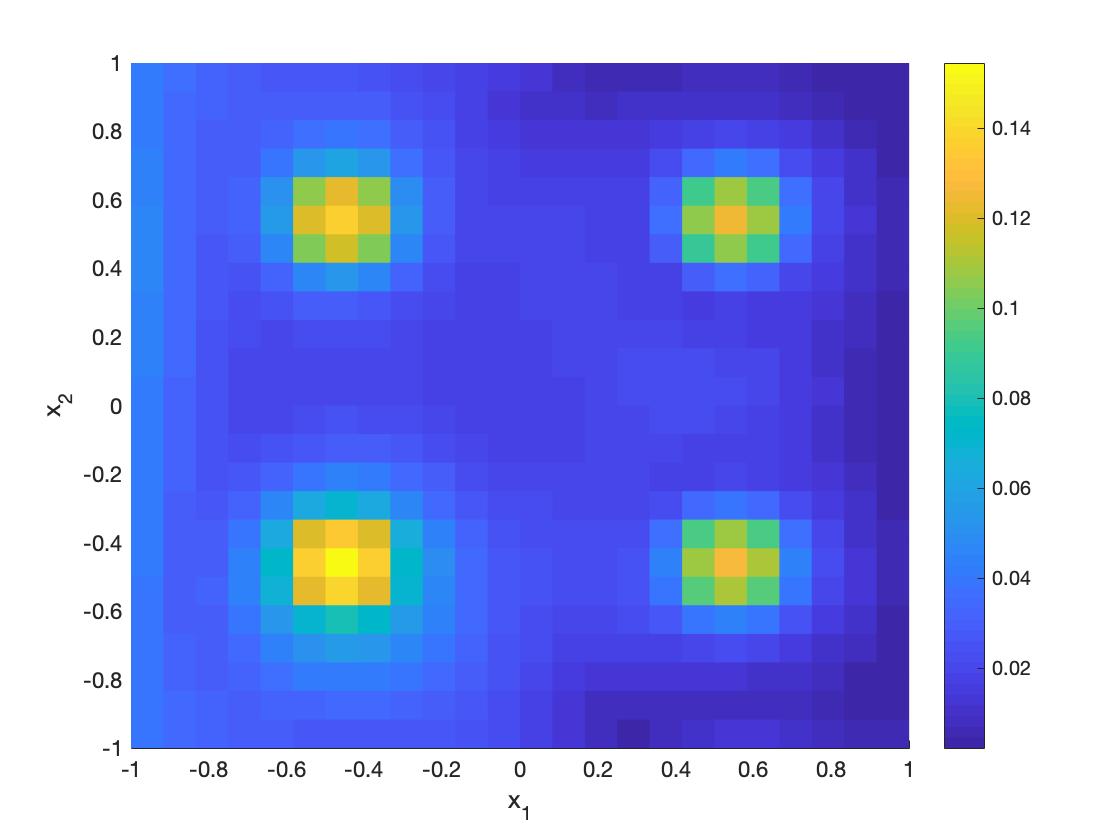}}
\subfloat[{relative error with data on $\Gamma_2$}]{ \includegraphics[width=0.3\textwidth]{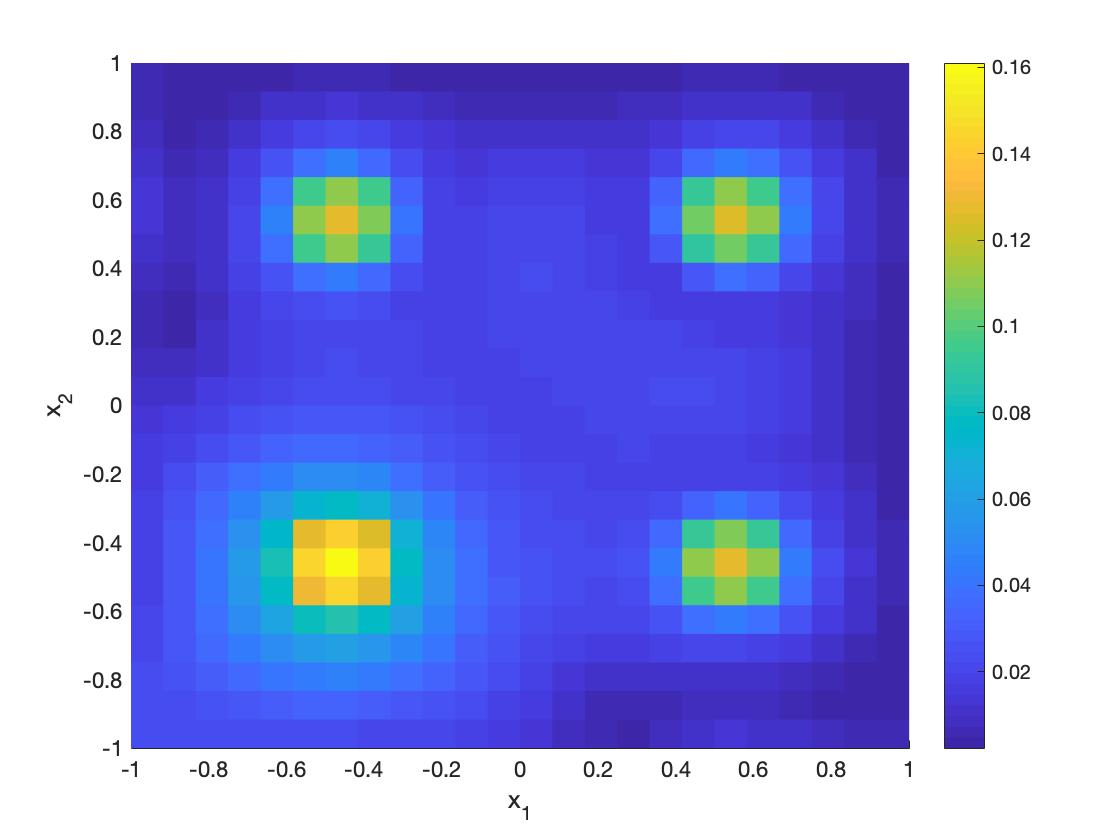}}
\subfloat[{relative error with data on $\Gamma_3$}]{ \includegraphics[width=0.3\textwidth]{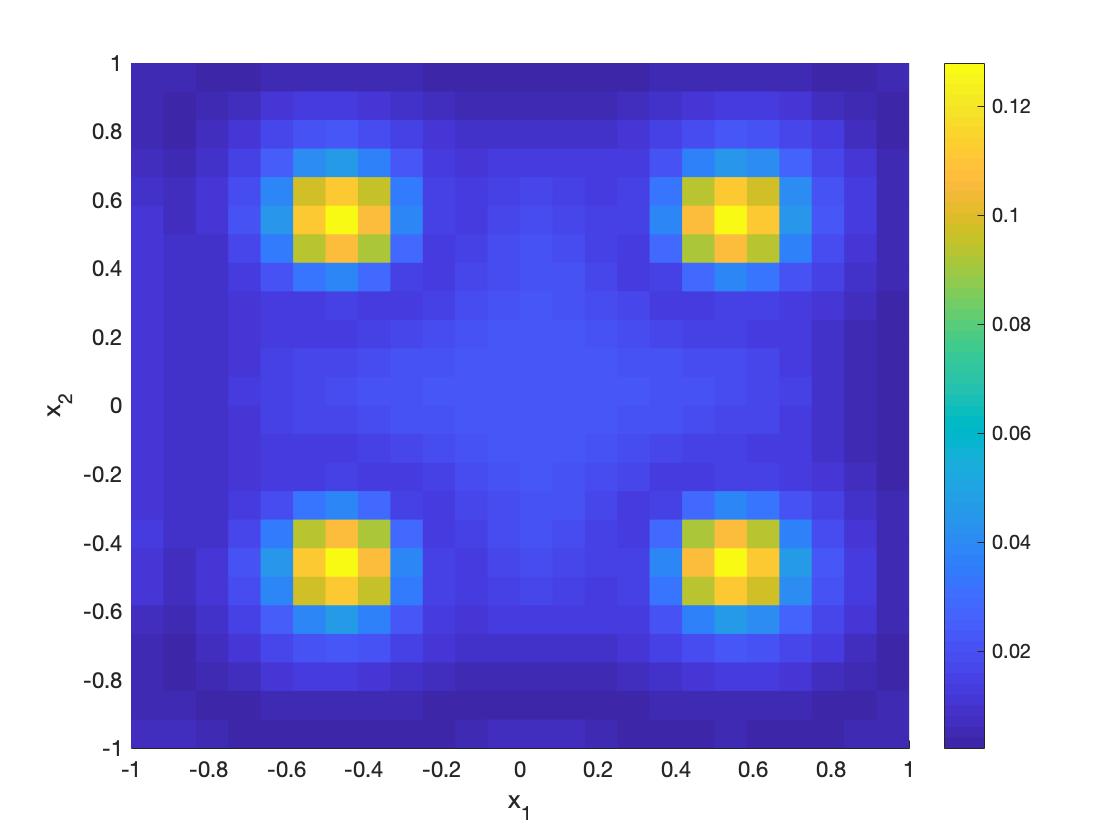}}
\caption{Approximations and relative errors for Example \ref{ex4}(b) with Cauchy data be given on different partial of boundary.}\label{ex4err}
\end{figure}

According to the conditional stability and convergence estimate we analyzed in section 2 and 3, the length of partial of boundary for which Cauchy data be given, will also effect on the approximation. Thus, we verify the proposed method for this example with Cauchy data be given in $\Gamma_j,j=1,2,3$ by Figure \ref{ex4err}.
%\section{circle domain}

%\subsection{forward problem}

%In this section, we consider the problem in circle domain. We use radial basis function method to slove the forward problem.
%\begin{equation}
%\begin{cases}
%du-(u_{xx}+u_{yy})dt=udW(t)&\quad \mbox{in}\quad (0,1)\times G\\
%u(x,y,0)=y_{0}(x,y)&\quad \mbox{in}\quad \lbrace 0\rbrace\times G\\
%u(x,y,t)=g(x,y,t)&\quad \mbox{in}\quad (0,1)\times\partial G
%\end{cases}\label{circlrfor}
%\end{equation}
In Example \ref{ex4}, $G$ is a rectangular domain that the initial-boundary problem could be computed well by the finite difference method. However, in case that $G$ is a general bounded domain, we should improve the finite difference method by heterogeneous grids. To avoid the complicated analysis for the numerical solution of  initial-boundary problem, we solve the problem by kernel-based learning theory, by treating the fundamental solution as kernels.

\begin{figure}[htbp]
  \centering
\subfloat[{approximation with $\delta=1\%$ }]{\includegraphics[width=0.3\textwidth]{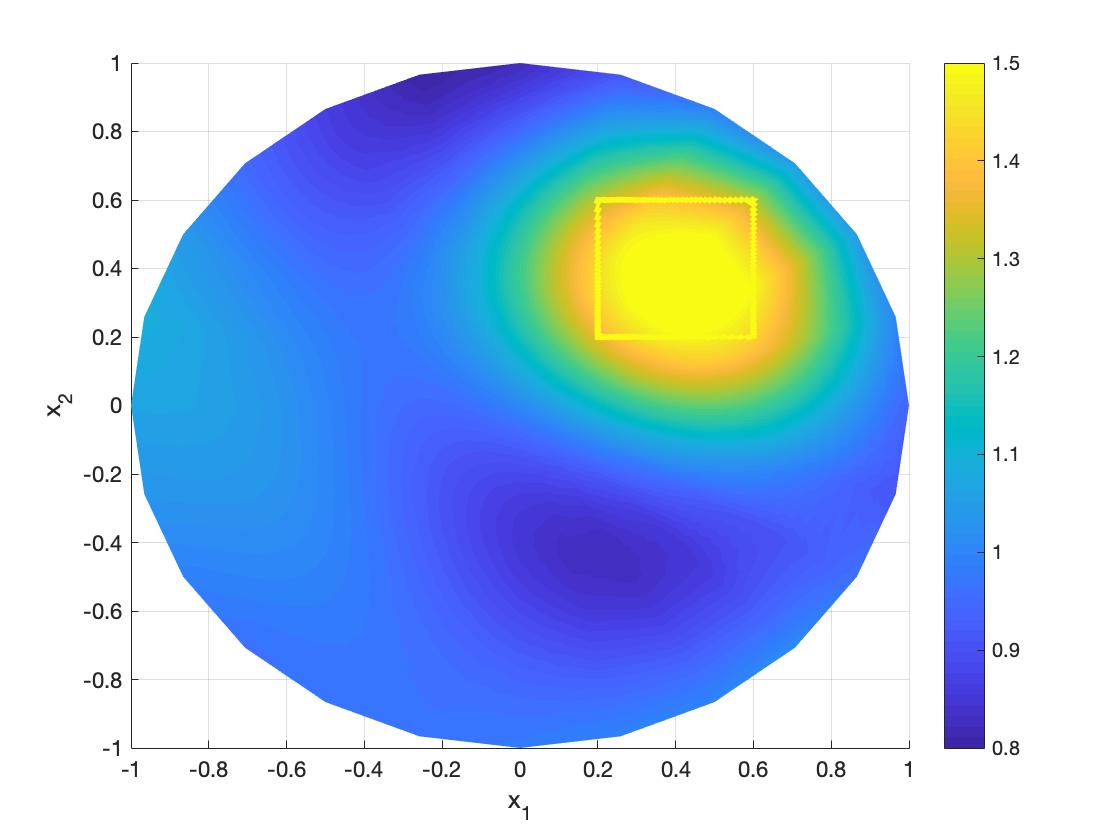}}
\subfloat[{approximation with $\delta=5\%$ }]{ \includegraphics[width=0.3\textwidth]{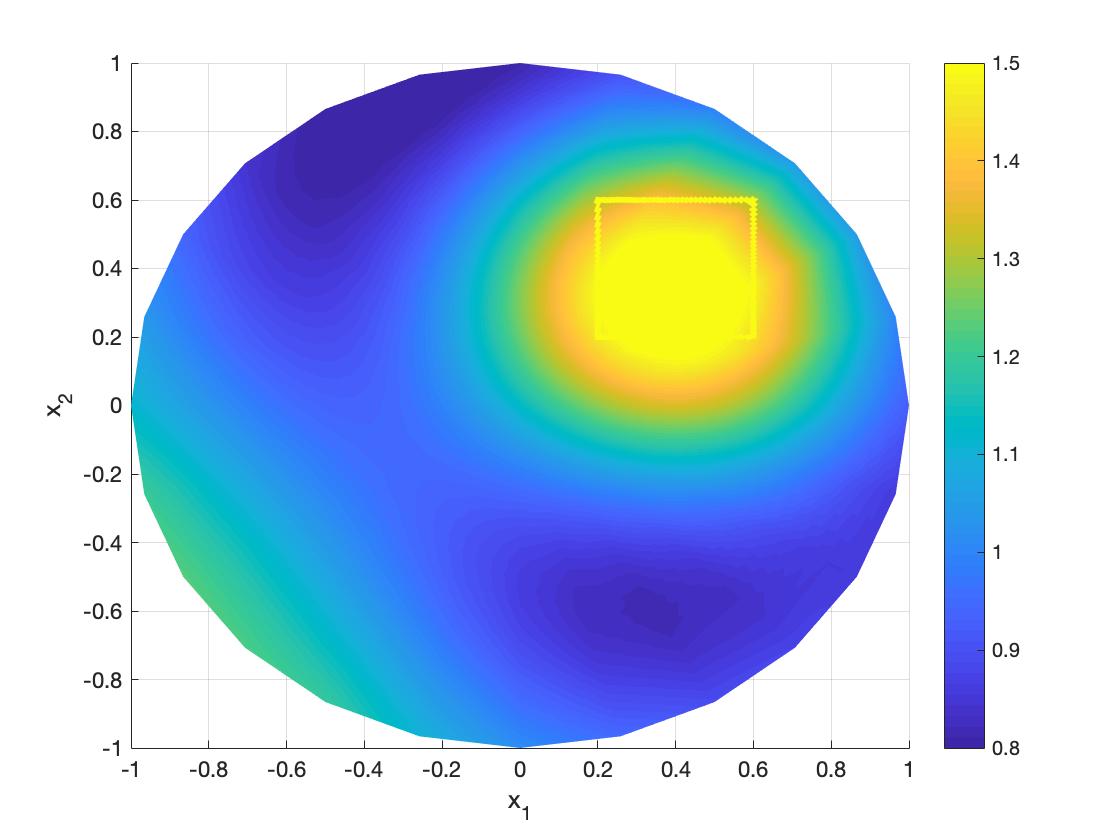}}
\subfloat[{approximation with $\delta=10\%$ }]{\includegraphics[width=0.3\textwidth]{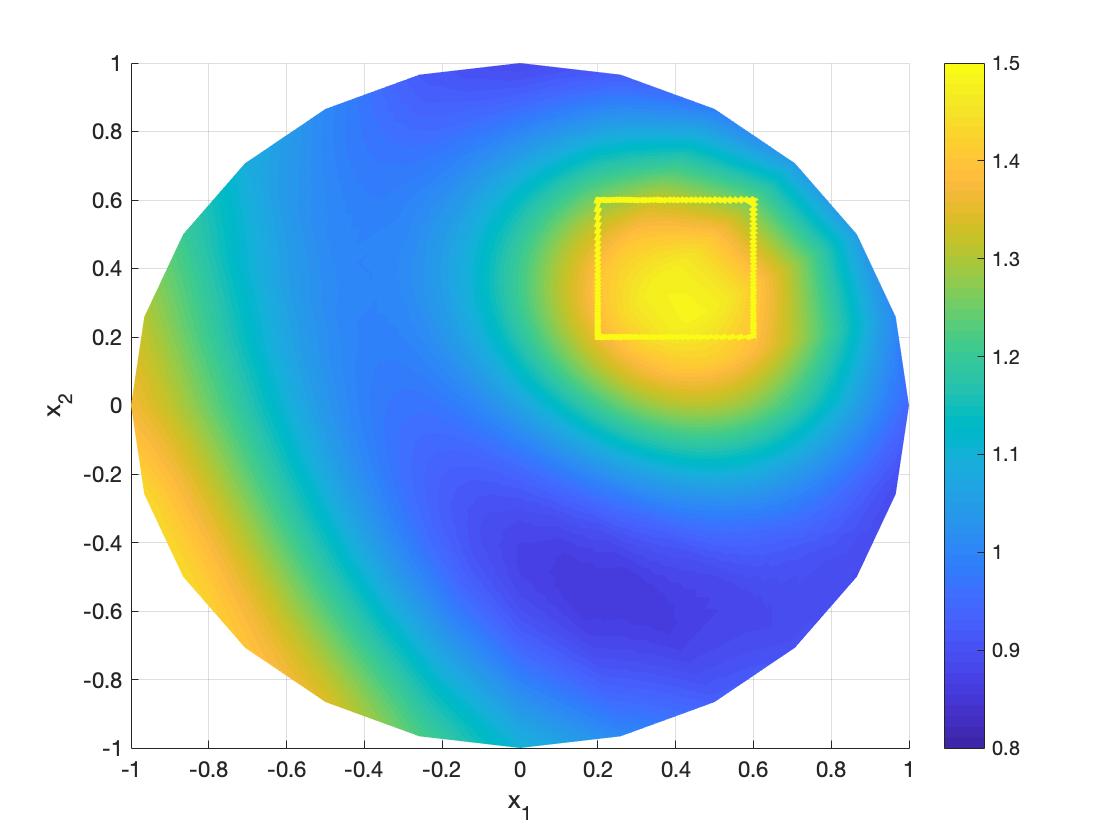}}\\
\subfloat[{relative error with $\delta=1\%$}]{\includegraphics[width=0.3\textwidth]{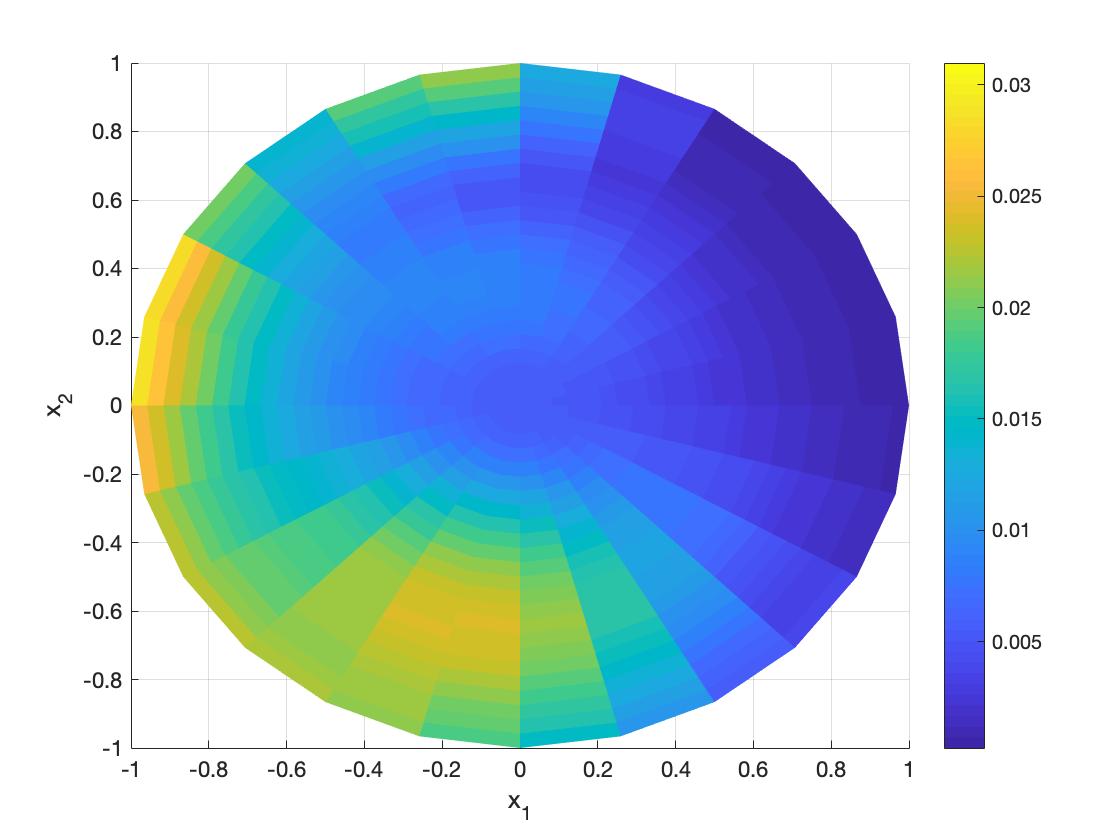}}
\subfloat[{relative error with $\delta=5\%$}]{ \includegraphics[width=0.3\textwidth]{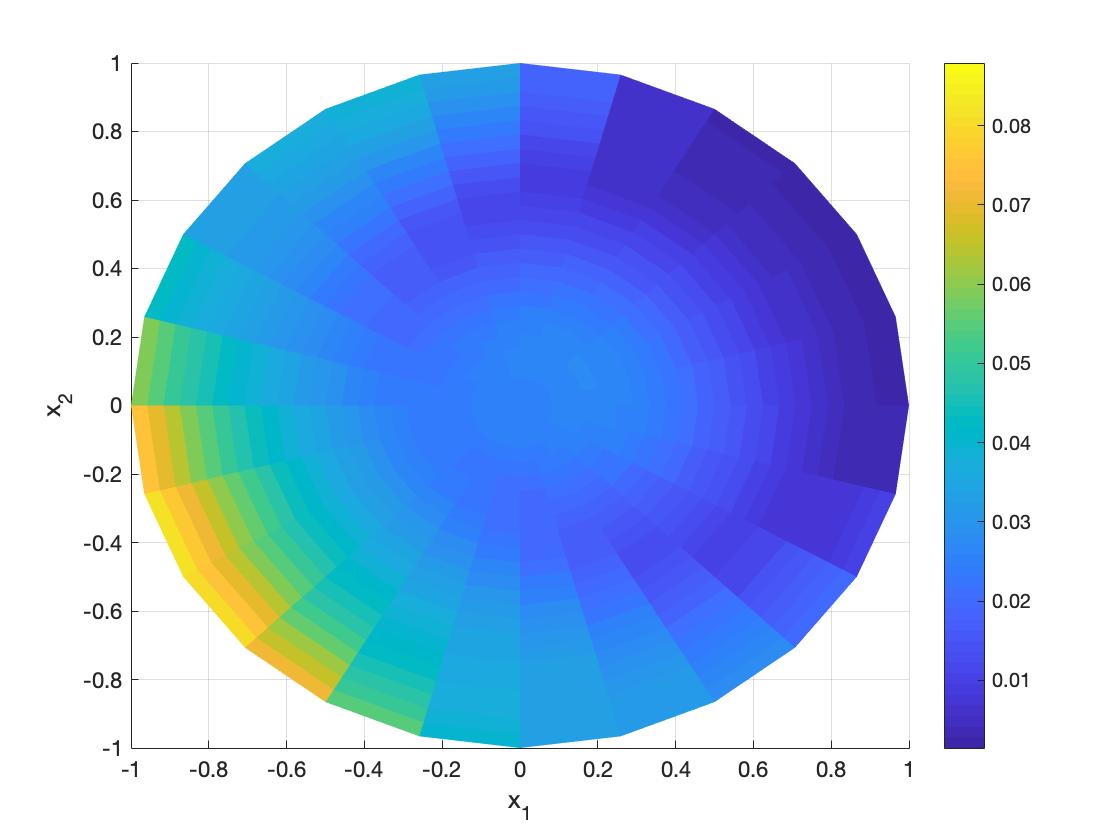}}
\subfloat[{relative error with $\delta=10\%$}]{\includegraphics[width=0.3\textwidth]{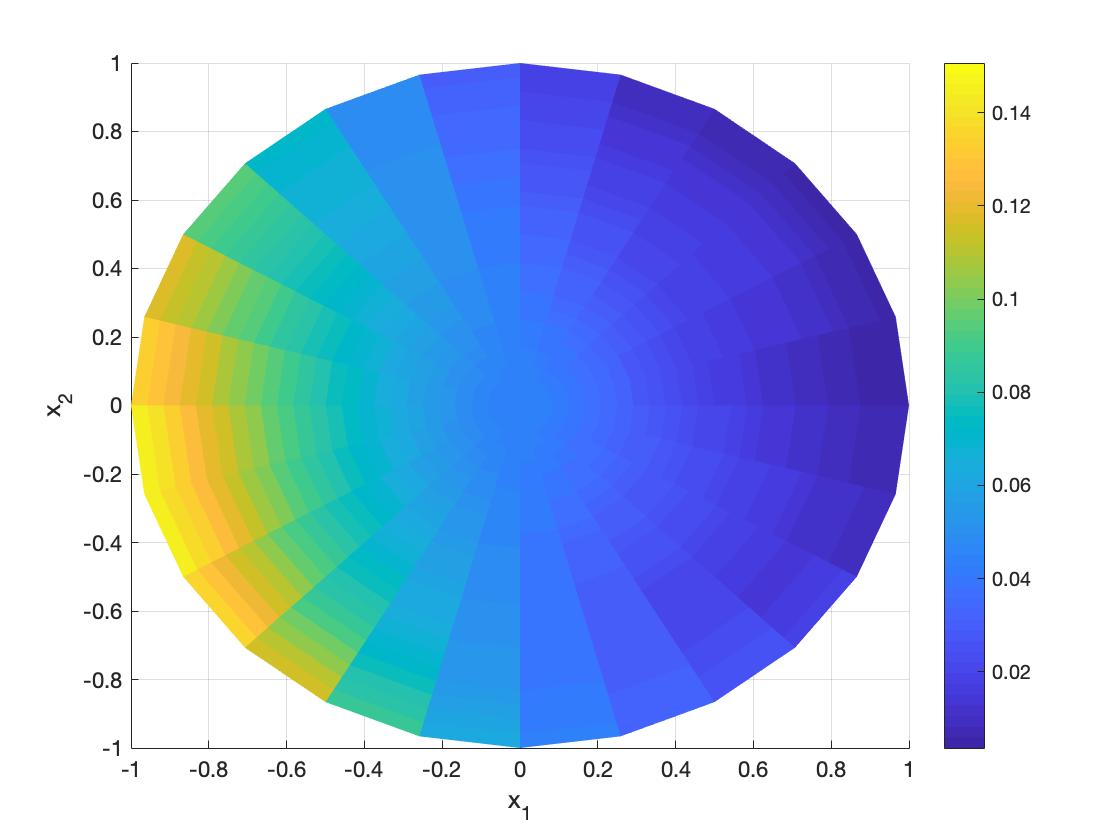}}
\caption{Numerical results with different noisy level for Example \ref{ex5}.}\label{ex5f2}
\end{figure}

\begin{figure}[htbp]
  \centering
\subfloat[{approximation with data on $\Gamma_\frac{\pi}{6}$}]{\includegraphics[width=0.3\textwidth]{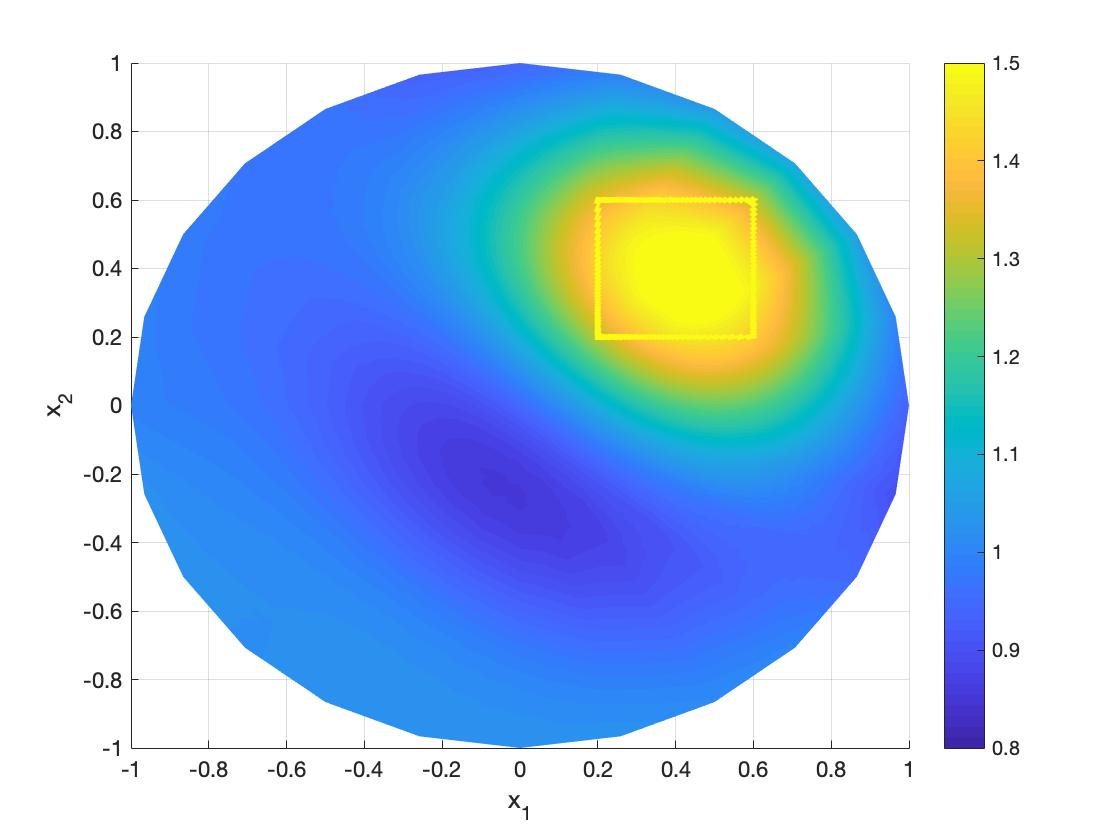}}
\subfloat[{approximation with data on $\Gamma_\frac{\pi}{4}$}]{ \includegraphics[width=0.3\textwidth]{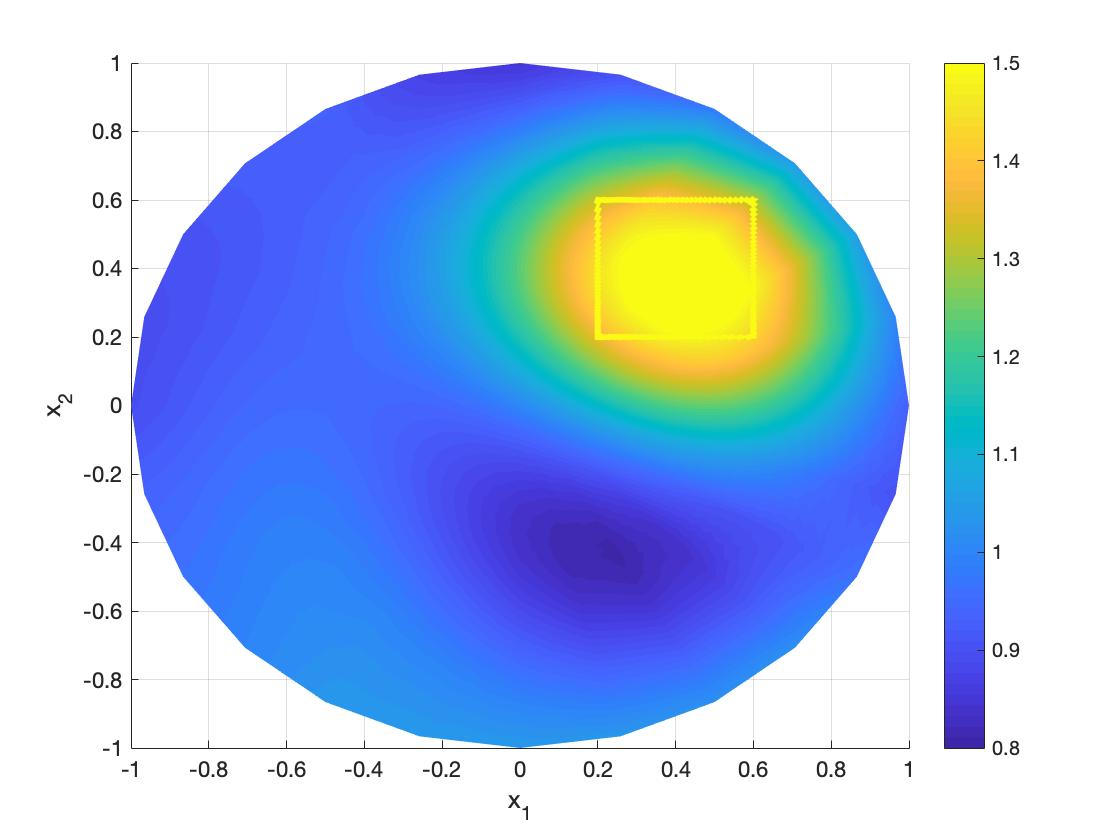}}
\subfloat[{approximation with data on $\Gamma_\frac{\pi}{2}$}]{\includegraphics[width=0.3\textwidth]{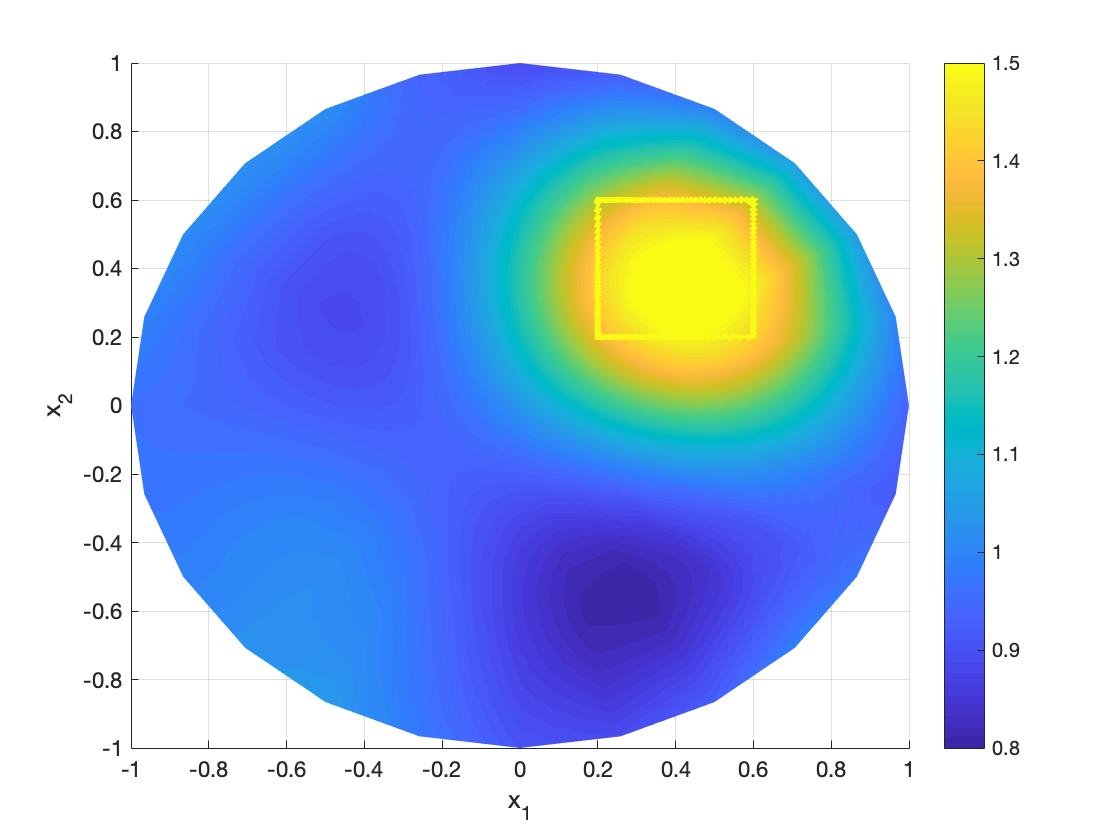}}\\
\subfloat[{relative error with data on $\Gamma_\frac{\pi}{6}$}]{\includegraphics[width=0.3\textwidth]{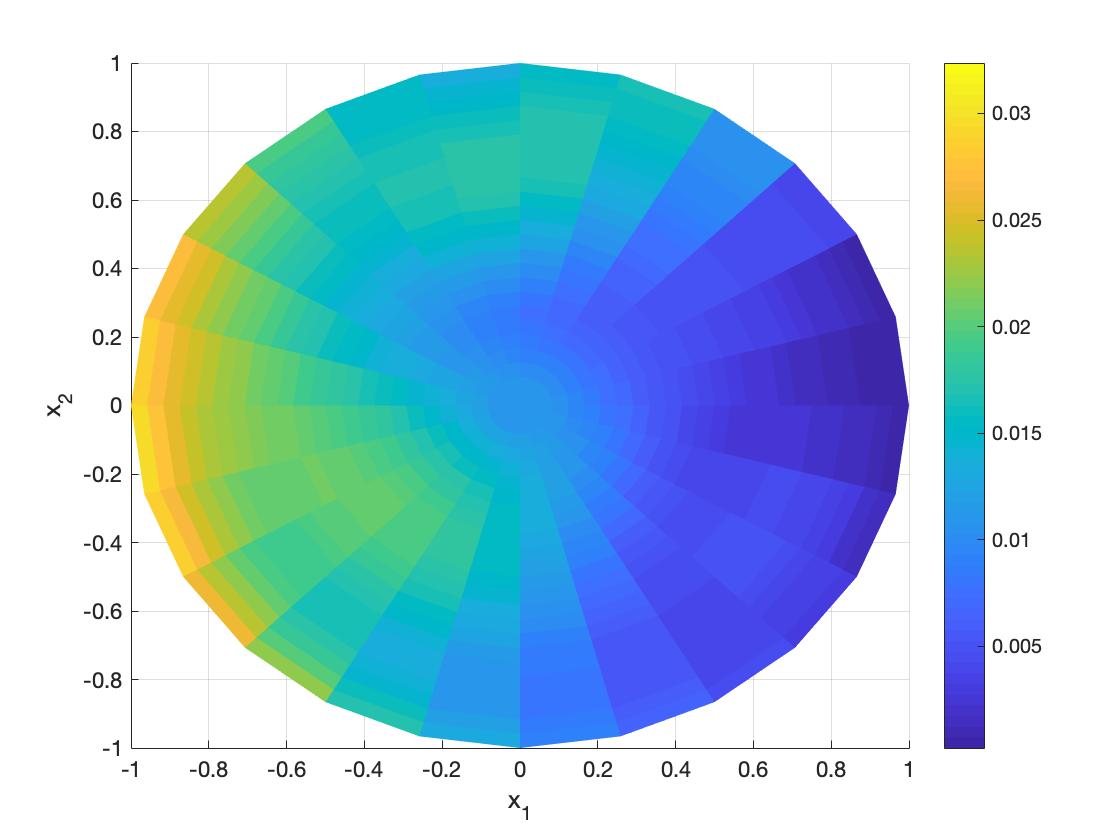}}
\subfloat[{relative error with data on $\Gamma_\frac{\pi}{4}$}]{ \includegraphics[width=0.3\textwidth]{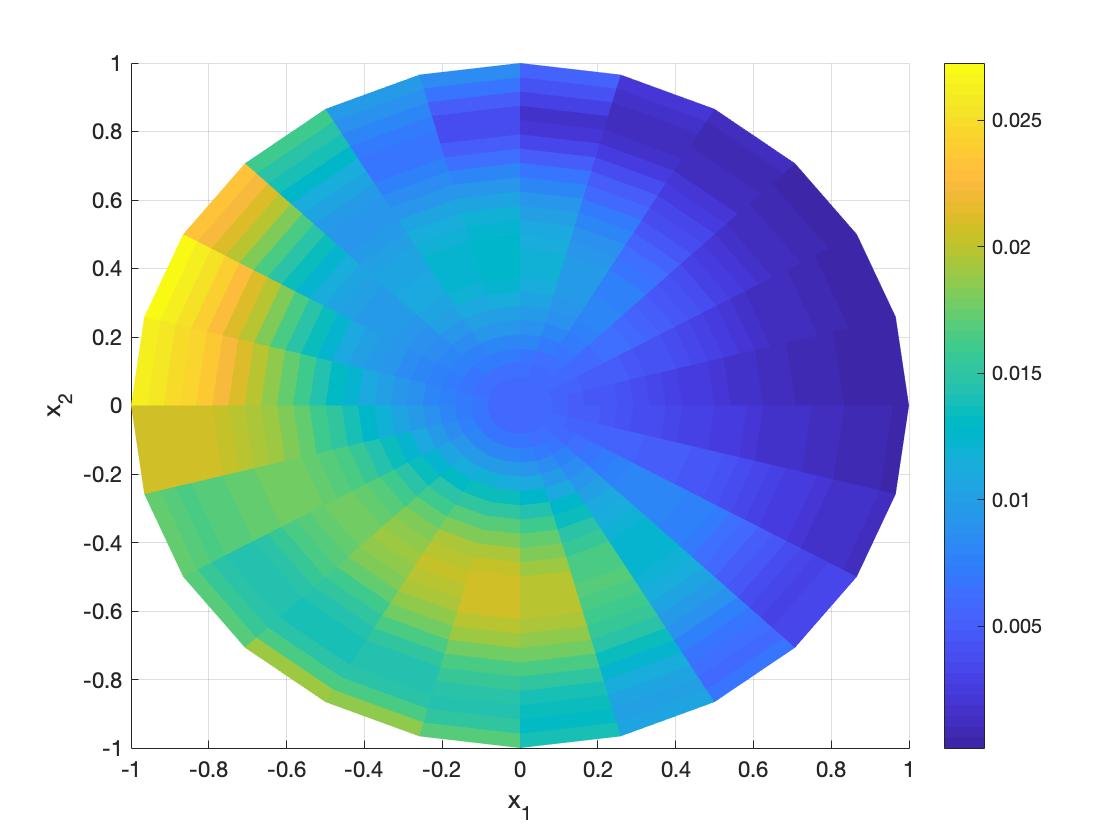}}
\subfloat[{relative error  with data on $\Gamma_\frac{\pi}{2}$}]{\includegraphics[width=0.3\textwidth]{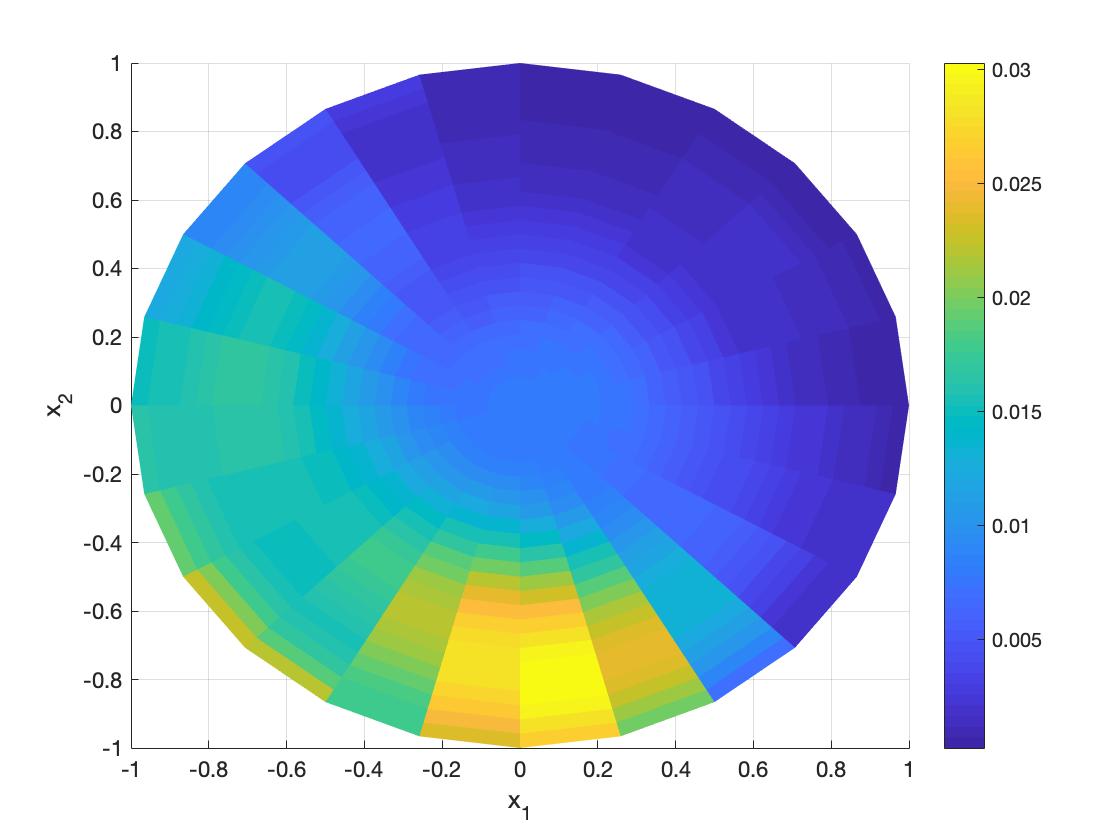}}
\caption{Numerical results with Cauchy data on different partial boundary for Example \ref{ex5}.}\label{ex5f4}
\end{figure}

\begin{example}\label{ex5}
Suppose that $G=\lbrace r^2\leq 1\rbrace$, where $r\triangleq \sqrt{x_1^2+x_2^2}$, $\tan \theta=x_2/x_1$, and $\Gamma_{\Theta}=\lbrace \theta\in\lbrack 0,\Theta\rbrack, r=1\rbrace$. Let
$
y_{0}(x_1,x_2)=
\begin{cases}
3,&0.2\leq x_1\leq 0.6,0.2\leq x_2\leq 0.6\\
1,&otherwise,
\end{cases},\quad g_1(t,x_1,x_2)=1.
$\end{example}
We show the numerical results when $\Theta=\frac\pi3$ with varies $\delta$ and results with different partial of boundary for $\Theta=\frac\pi6,\frac\pi4,\frac\pi2$ for $\delta=3\%$ in Figure \ref{ex5f2} and \ref{ex5f4}, respectively.

In conclusion, we do the numerical experiments for the Cauchy problem of stochastic parabolic equation for examples on both  1-D and 2-D  bounded domains by smooth, piecewise smooth and noncontinuous initial conditions, numerical results illustrate that the proposed method works well for all these examples. Although we haven't solved examples for  other irregular 2-D domains and 3-D domains, it can be seen from the kernel-based learning theory that the method should be satisfied, because there is no limitation of dimension $n$ and shape of domains in the kernel-based learning theory but only distances between source points and collection point be concerned. Moreover, since no iterated method used in the numerical computation, the computational cost is reduced. However, the convergence rate for the kernel-based learning theory is still open.

\section*{Acknowledgement}
The authors would like to thank Professor Kui Ren for providing valuable  suggestions. This work is partially supported by the NSFC (No. 12071061), the Science Fund for Distinguished Young Scholars of Sichuan Province (No. 2022JDJQ0035).

%\bibliographystyle{plain}
%\bibliography{bibcauchy}

\begin{thebibliography}{10}

\bibitem{barbu2}
V. Barbu and M. R{\"o}ckner.
\newblock Backward uniqueness of stochastic parabolic like equations driven by
  gaussian multiplicative noise.
\newblock {\em Stochastic Processes and their Applications}, 126(7):2163--2179,
  2016.

\bibitem{BK}
L. Beilina and M.V. Klibanov.
\newblock {\em Approximate global convergence and adaptivity for coefficient
  inverse problems}.
\newblock Springer Science \& Business Media, 2012.

\bibitem{2009}
R.C. Dalang, D. Khoshnevisan, C. M\"uller, D. Nualart, and Y.M.
  Xiao.
\newblock {\em A minicourse on stochastic partial differential equations},
  volume 1962.
\newblock Springer, 2009.

\bibitem{DD2022IP}
F.F. Dou and W. Du.
\newblock Determination of the solution of a stochastic parabolic equation by
  the terminal value.
\newblock {\em Inverse Problems}, 38(7):075010, 2022.

\bibitem{DH2012}
F.F. Dou and Y.C. Hon.
\newblock Kernel-based approximation for cauchy problem of the time-fractional
  diffusion equation.
\newblock {\em Engineering Analysis with Boundary Elements}, 36(9):1344--1352,
  2012.

\bibitem{Fu}
X. Fu and X. Liu.
\newblock A weighted identity for stochastic partial differential operators and
  its applications.
\newblock {\em Journal of Differential Equations}, 262(6):3551--3582, 2017.

\bibitem{1986An}
J.B. Walsh.
\newblock {\em An introduction to stochastic partial differential equations}.
\newblock Springer, 1986.

\bibitem{Glimm} J. Glimm. 
\newblock {\em    Nonlinear and stochastic phenomena:
the grand challenge for partial differential
equations.} 
\newblock  {\em SIAM Review}, \rm {\bf 33}, 626--643, 1991.

\bibitem{T1992}T. Takeshi.
 \newblock{\em Successive approximations to solutions of stochastic differential equations.}  
 \newblock {\em Journal of Differential Equations}, \rm{\bf 96}, 152-199, 1992.

\bibitem{KP1992}
P.E. Kloeden,and E. Platen.
\newblock {\em  Numerical Solution of Stochastic Differential Equations.}
\newblock Springer, 1992. 

\bibitem{GLWX2021IP}
Y. Gong, P. Li, X. Wang, and X. Xu.
\newblock Numerical solution of an inverse random source problem for the time
  fractional diffusion equation via phaselift.
\newblock {\em Inverse Problems}, 37(4):045001, 2021.

\bibitem{Hansen}
P.C. Hansen.
\newblock Numerical tools for analysis and solution of fredholm integral
  equations of the first kind.
\newblock {\em Inverse Problems}, 8(6):849-872, 1992.

\bibitem{HW2004}
Y.C. Hon and T. Wei.
\newblock A fundamental solution method for inverse heat conduction problem.
\newblock {\em Engineering Analysis with Boundary Elements}, 28(5):489--495,
  2004.

\bibitem{IK2000}
V.~Isakov and S.~Kindermann.
\newblock Identification of the diffusion coefficient in a one-dimensional
  parabolic equation.
\newblock {\em Inverse Problems}, 16(3):665--680, 2000.

\bibitem{KZ1998}
Y.L. Keung and J.~Zou.
\newblock Numerical identifications of parameters in parabolic systems.
\newblock {\em Inverse Problems}, 14(1):83--100, 1998.

\bibitem{KT2007}
M.V. Klibanov and A.V. Tikhonravov.
\newblock Estimates of initial conditions of parabolic equations and
  inequalities in infinite domains via lateral cauchy data.
\newblock {\em Journal of Differential Equations}, 237(1):198--224, 2007.

\bibitem{Klib2006}
M.V. Klibanov.
\newblock Estimates of initial conditions of parabolic equations and
  inequalities via lateral Cauchy data.
\newblock  {\em Inverse Problems}, 22(2):495, 2006.

\bibitem{Kotelenez} P.~Kotelenez. \newblock {\em Stochastic ordinary and
stochastic partial differential equations.
Transition from microscopic to macroscopic
equations}. 
\newblock  Springer, New York, 2008.

\bibitem{Ksurvey}
M.V. Klibanov.
\newblock Carleman estimates for global uniqueness, stability and numerical
  methods for coefficient inverse problems.
\newblock {\em Journal of Inverse and Ill-Posed Problems}, 21(4):477--560,
  2013.


\bibitem{KKLY2017}
Klibanov, Michael V.; Koshev, Nikolaj A.; Li, Jingzhi; Yagola, Anatoly G.
\newblock Numerical solution of an ill-posed Cauchy problem for a quasilinear parabolic equation using a Carleman weight function
\newblock {\em Journal of Inverse and Ill-Posed Problems}, 24(6):761--776, 2016.

\bibitem{LL}
H. Li and Q.~L{\"u}.
\newblock A quantitative boundary unique continuation for stochastic parabolic
  equations.
\newblock {\em Journal of Mathematical Analysis and Applications},
  402(2):518--526, 2013.

\bibitem{LYZ2009CPAA}
J. Li, M. Yamamoto, and J. Zou.
\newblock Conditional stability and numerical reconstruction of initial
  temperature.
\newblock {\em Communications on Pure \& Applied Analysis}, 8(1):361--382,
  2009.

\bibitem{Lu1}
Q.~L{\"u}.
\newblock Carleman estimate for stochastic parabolic equations and inverse
  stochastic parabolic problems.
\newblock {\em Inverse Problems}, 28(4):045008, 2012.

\bibitem{LZ2021}
Q.~L{\"u} and X.~Zhang. 
\newblock {\em Mathematical control theory for stochastic partial differential equations}. \newblock{Springer}, Cham,  2021. 

\bibitem{LZ2023}
Q.~L{\"u} and X.~Zhang.
\newblock Inverse problems for stochastic partial differential equations: Some
  progresses and open problems, 
  \newblock{\em Numerical Algebra, Control and Optimization}, 2023, DOI: 10.3934/naco.2023014.

\bibitem{MMRTS}
K.R. M{\"u}ller, S. Mika, K. Tsuda, G. Ratschnd K. Sch{\"o}lkopf.
\newblock An introduction to kernel-based learning algorithms.
\newblock In {\em IEEE Transactions on Neural Networks}, 12(2):181:201, 2001.
\bibitem{NHZ2020IP}
P. Niu, T. Helin, and Z. Zhang.
\newblock An inverse random source problem in a stochastic fractional diffusion
  equation.
\newblock {\em Inverse Problems}, 36(4):045002, 2020.

\bibitem{Schaback2007}
R. Schaback.
\newblock Convergence of unsymmetric kernel-based meshless collocation methods.
\newblock {\em SIAM Journal on Numerical Analysis}, 45(1):333--351, 2007.

\bibitem{Tang-Zhang1}
S. Tang and X.~Zhang.
\newblock Null controllability for forward and backward stochastic parabolic
  equations.
\newblock {\em SIAM Journal on Control and Optimization}, 48(4):2191--2216,
  2009.

\bibitem{wu}
B. Wu, Q. Chen, and Z. Wang.
\newblock Carleman estimates for a stochastic degenerate parabolic equation and
  applications to null controllability and an inverse random source problem.
\newblock {\em Inverse Problems}, 36(7):075014, 2020.

\bibitem{Y}
M. Yamamoto.
\newblock Carleman estimates for parabolic equations and applications.
\newblock {\em Inverse Problems}, 25(12):123013, 2009.

\bibitem{yin2014}
Z. Yin.
\newblock A quantitative internal unique continuation for stochastic parabolic
  equations.
\newblock {\em Mathematical Control and Related Fields}, 5(1):165--176, 2015.

\bibitem{YY2009}
G.~Yuan and M.~Yamamoto.
\newblock Lipschitz stability in the determination of the principal part of a
  parabolic equation.
\newblock {\em ESAIM: Control, Optimisation and Calculus of Variations},
  15(3):525--554, 2009.

\bibitem{yuan2017}
G. Yuan.
\newblock Conditional stability in determination of initial data for stochastic
  parabolic equations.
\newblock {\em Inverse Problems}, 33(3):035014, 2017.

\bibitem{Zhang2008}
X.~Zhang.
\newblock Unique continuation for stochastic parabolic equations.
\newblock {\em Differential Integral Equations}, 21:81--93, 2008.

\end{thebibliography}

\end{document}